\documentclass{article}
\usepackage{epsfig}
\usepackage{graphics}
\usepackage{color}
\usepackage{amsmath}
\usepackage{theorem}
\usepackage{amssymb}
\theoremheaderfont{\scshape}

\newtheorem{Proposition}{Proposition}
  \newtheorem{Remark}[Proposition]{Remark}
  \newtheorem{Corollary}[Proposition]{Corollary}
  \newtheorem{Lemma}[Proposition]{Lemma}
  
  \newtheorem{Theorem}[Proposition]{Theorem}

    \def\Tm{\begin{Theorem}\label}
    \def\Pp{\begin{Proposition}\label}
    \def\Rm{\begin{Remark}\label}
    \def\Lm{\begin{Lemma}\label}
    \def\Co{\begin{Corollary}\label}

    \def\eTm{\end{Theorem}}
    \def\ePp{\end{Proposition}}
    \def\eRm{\end{Remark}}
    \def\eLm{\end{Lemma}}
    \def\eCo{\end{Corollary}}

    \def\bfC{{\bf C}}
    
    \def\Qp{{\bf Q}_\perp}

    \def\dpp{{\bf D}_\perp}
    
    \def\Rp{{\bf R}_\perp}

\makeatletter
\@addtoreset{equation}{section}
\makeatother

    \def\Box{{\hfill\hbox{\enspace${\sqre}$}} \smallskip}
    \def\sqr#1#2{{\vcenter{\vbox{\hrule height .#2pt
                             \hbox{\vrule width .#2pt height#1pt \kern#1pt
                                   \vrule width .#2pt}
                             \hrule height .#2pt}}}}
    \def\sqre{\mathchoice\sqr54\sqr54\sqr{4.1}3\sqr{3.5}3}

\newcommand{\bT}[1]{\begin{theorem}\label{#1}}
\newcommand{\be}[1]{\begin{equation}\label{#1}}
\newcommand{\ba}[1]{\begin{eqnarray}\label{#1}}
\newcommand{\ee}{\end{equation}}
\newcommand{\ea}{\end{eqnarray}}
\newcommand{\bl}[1]{\begin{lemma}\label{#1}}
\newcommand{\bp}[1]{\begin{proposition}\label{#1}}
\newcommand{\br}[1]{\begin{remark}\label{#1}}
    
    \def\bchi{\mbox{\raisebox{.4ex}{\begin{Large}$\chi$\end{Large}}}}

    \def\bflam{{\boldsymbol \lambda}}
    
    \def\bfbet{{\boldsymbol \beta}}
    
    \def\bfW{\mathbf{W}}
    \def\bfdel{\boldsymbol\delta}

    \def\RR{\mathbb{R}}
    \def\DD{\mathbb{D}}
    
    \def\CC{\mathbb{C}}
    
    \def\NN{\mathbb{N}}
    \def\ZZ{\mathbb{Z}}
    
    \def\z{\noindent}
    \def\bff{{\bf f}}
    \def\bffz{{\bf F}_0}
    \def\bfd{\mathbf{D}}
    \def\bfdd{\mathbf{d}}

    \def\bfh{{\bf h}}
    \def\bfq{{\bf q}}
    \def\bfH{{\bf H}}\def\bfY{{\bf Y}}

    \def\bfl{{\bf l}}
    \def\bfgg{\mathbf{g}}
    \def\bfg{{\bf G}}
    
    \def\bfv{{\mathbf{v}}}
    \def\bfV{\mathbf{V}}
    \def\bog{{\bf g}}
        
    \def\calv{{\DD_\epsilon}}
    
    \def\calnb{{{\cal N}}}
    \def\lloc{{L}^1_{loc}}

    \def\lone{{L^1}}

    \def\ga{\hat{\Gamma}}
    \def\gc{\hat{\Gamma}_c}
    \def\hatc{\hat C}
    \def\bc{\hat{B}_c}
    \def\bcp{\hat{B}'_c}
    \def\hb{\hat{B}}
    
    \def\bfii{\mathbf{i}}
    \def\bfa{{\bf a}}
    \def\bfB{{\bf B}}
    \def\bfA{{\bf A}}
    \def\bfQ{{\bf Q}}
    \def\zpp{(z\oplus\gc(z))}
    \def\ppp{(p\oplus\gc(-p))}

        \def\betapp{(\beta\oplus\bc)}

    \def\bfT{{\bf T}}
    \def\bft{{\bf t}}
    \def\bfy{{\bf y}}

    \def\bfv{{\bf v}}
    \def\bfR{{\bf R}}
    
    \def\lap{{\cal L}}

    \def\bfk{{\bf k}}
    \def\bfj{{\bf j}}
    
    \def\bfm{{\bf m}}
    \def\bfdl{{\bf d_l}}

    \def\lap{{\cal L}}
    \def\lapi{{\cal L}^{-1}}
    \def\bor{{\cal B}}

    \def\bfk{{\bf k}}
    
    \def\bfm{{\bf m}}
    \def\bfdl{{\bf d_l}}
    
    \def\heav{{\cal H}}

    \def\bflam{{\boldsymbol \lambda}}
    
    \def\bfbet{{\boldsymbol \beta}}
    
    \def\bfW{\mathbf{W}}
    \def\bfdel{\boldsymbol\delta}
\begin{document}

\title{ On Borel summation and Stokes phenomena for rank one
  nonlinear systems of ODE's }
\author{Ovidiu Costin \thanks {Mathematics Department,
University of Chicago, 5734 University Avenue, Chicago
IL 60637; e-mail: costin\symbol{64}math.uchicago.edu}}
\date{ }

\maketitle

\hyphenation{trans-series}

         \section{Introduction}

  In this paper we study analytic (linear or) nonlinear systems of
  ordinary differential equations, at an irregular singularity of rank
  one, under nonresonance conditions.  It is shown that the formal
  asymptotic exponential series solutions ({\em transseries}
  solutions: countable linear combinations of formal power series
  multiplied by small exponentials) are Borel summable in a
  generalized sense along any direction in which the exponentials
  decay.  Conversely, any solution that decreases along some direction
  is the Borel sum of a transseries.
  
  The summation procedure introduced is an extension of Borel
  summation which is linear, multiplicative, commutes with
  differentiation and complex conjugation. The summation algorithm
  uses the formal solutions alone (and not the differential
  equation that they solve). Along singular (Stokes) directions, the
  functions reconstructed by summation  are shown to be given by Laplace integrals
  along special paths, a subset of \'Ecalle's median paths.

  The one-to-one correspondence established between actual solutions
  and generalized Borel sums of transseries is constant between Stokes
  lines and {changes} if a Stokes line is crossed (local Stokes
  phenomenon). We analyze the connection between local and classical
  Stokes phenomena.
 
  We study the analytic properties of the Borel (formal inverse
  Laplace) transform of the series contained in the transseries of the
  transseries and give a systematic description of their
  singularities. These Borel transforms satisfy a hierarchy of
  convolution equations, for which we give the general solution in a
  space of hyperfunctions.  In addition, we show that they are {\em
    resurgent functions} in the sense of \'Ecalle.

 The summation
  procedure is not unique; we classify all proper
  extensions of Borel summation to  transseries solutions of nonresonant
  systems.
   
  We find formulas connecting the
  different series contained in the  transseries among themselves (resurgence
  equations). Resurgence turns out to be closely linked to the local
  Stokes phenomenon.

  The connection to Berry's hyperasymptotics and applications to the
  classification of differential equations are briefly discussed.

\subsection{General setting}

 We  consider the differential system
                
\begin{eqnarray}
 \label{eqor1}
  \bfy'=\mathbf{f}(x,\bfy)  \qquad \bfy\in\CC^n              
   \end{eqnarray}

\z under the following {\em assumptions:}

(a1) The function $\mathbf{f}$ is analytic
at $(\infty,0)$.

(a2) Nonresonnance: the eigenvalues $\lambda_i$ of the linearization

\begin{eqnarray}
  \label{linearized}
  \hat{\Lambda}:=-\left(\frac{\partial f_i}{\partial
    y_j}(\infty,0)\right)_{i,j=1,2,\ldots n}
\end{eqnarray}

\z are linearly independent over $\ZZ$ (in particular nonzero) and
such that the Stokes lines are distinct (a somewhat less
restrictive condition is actually used, cf. \S\ref{nonres}).

\z {\em Normalization}. It is convenient
to prepare (\ref{eqor1}) in the following way. Pulling out the inhomogeneous 
and the linear terms (relevant to leading order
asymptotics) we get

\begin{eqnarray}\label{eqor}
{\bf y}'={\bf f}_0(x)-\hat\Lambda {\bf y}-
\frac{1}{x}\hat B {\bf y}+{\bf g}(x,{\bf y})
\end{eqnarray}

Under the assumptions (a1) and (a2), by means of normal form
calculations it is possible to arrange (\ref{eqor}) so that
(\cite{Wasow}, \cite{To1})

(n1)
$\hat\Lambda=\mbox{diag}(\lambda_i)$ and 

(n2) $\hat B=\mbox{diag}(\beta_i)$

\z  For convenience, we rescale $x$ and
reorder the components of $\bfy$ so that

(n3) $\lambda_1=1$, and, with $\phi_i=\arg(\lambda_i)$, we have
$\phi_i<\phi_j$ if $i<j$. To simplify notations, we formulate some of
our results relative to $\lambda_1$; they can be easily adapted to any
other eigenvalue.

  To unify the treatment we
make, by taking $\bfy=\bfy_1 x^{-N}$ for some $N>0$,

(n4) $\Re(\beta_j)<0,\ j=1,2,\ldots,n$.

\z (there is an asymmetry at this point: the opposite inequality cannot be
achieved, in general, as simply and without violating analyticity at
infinity).  Finally, through a transformation of the form
$\bfy\leftrightarrow\bfy-\sum_{k=1}^M\bfa_k x^{-k}$ we arrange that

(n5) $ \mathbf{f}_0=O(x^{-M-1})\mbox{ and }\mathbf{g}(x,\bfy)=
O(\bfy^2,x^{-M-1}\bfy) $. We choose $M>1+\max_i\Re(-\beta_i)$.

{\em Formal solutions.} In prepared form, given (a1)
and (a2), (\ref{eqor}) admits 
an $n$--parameter family of formal exponential series solutions  (transseries)

\begin{eqnarray}
  \label{eqformgen,n}
   \tilde{\bfy}=\tilde{\bfy}_0+\sum_{\mathbf{k}\ge 0; |\mathbf{k}|>0}C_1^{k_1}\cdots C_n^{k_n}
\mathrm{e}^{-(\bfk\cdot\bflam) x}x^{\bfk\cdot\bfm}\tilde{\bfy}_{\bfk}
\end{eqnarray}

\z (see \cite{Wasow}, \cite{Cope},\cite{Iwano}, and also \S~\ref{sec:For} below) where
$m_i=1-\lfloor\beta_i\rfloor$,
 ($\lfloor\cdot\rfloor=$ integer part), $\bfC\in\CC^n$ is an arbitrary vector
of constants, and $\tilde{\mathbf{y}}_\bfk=x^{-\bfk(\bfbet+\bfm)}
\sum_{l=0}^{\infty}\bfa_{\bfk;l} x^{-l}$ are formal power series.

When $x$ is large in some direction $d$ in $\CC$, an important role is
played by the subset of transseries which are at the same time {\em
  asymptotic} expressions\footnote{An asymptotic
expansion of a function carries immediate information about behavior
of the
function near the expansion point (in contrast to
antiasymptotic expansions, e.g. a convergent doubly infinite Laurent series)}:  When there are infinitely many
exponentials in (\ref{eqformgen,n}) we ask that for all $i$ with
$C_i\ne 0$ we have $|\mathrm{e}^{-\lambda_i x}|\ll 1$ for large $x$ in
the given direction $d$ in $\CC$. Formally, agreeing to omit the
terms with $C_i=0$, with $x$ in $d$,
any {\em ascending} chain $\Re(-\bfk_1\cdot\bflam x)\le
\Re(-\bfk_2\cdot\bflam x)\le\ldots$, $\bfk_i\ne\bfk_j$, in
(\ref{eqformgen,n}) must be {\em finite} (the terms of an asymptotic
transseries are well-ordered with respect to $''\ll''$).  Thus for $x$ 
in some direction $d$  we only consider those transseries
that satisfy the condition:

(c1) $\xi+\phi_i:=\arg(x)+\phi_i\in (-\pi/2,\pi/2)$ for all $i$ such that $C_i\ne
0$. In other words, $C_i\ne 0$ implies that $\lambda_i$ lies in a
half-plane centered on $\overline{d}$, the complex conjugate direction
to $d$.

 From now on, ${\bflam}
=(\lambda_{i_1},\ldots,\lambda_{i_{n_1}})$,
$\bfbet=(\beta_{i_1},\ldots,\beta_{i_{n_1}})$,
$\bfm=(m_{i_1},\ldots,m_{i_{n_1}})$ and
$\bfbet'=\bfbet+\bfm$ where the indices $i_1,\ldots,i_{n_1}$ satisfy
(c1).  

We will henceforth consider that (\ref{eqor}) is presented in prepared
form, and use the designation transseries only for those formal solutions
satisfying (c1).

The series $\tilde{\bfy}_0$ is a formal solution of (\ref{eqor})
while, for $\bfk\ne 0$, $\tilde{\bfy}_\bfk$ satisfy a hierarchy of
linear differential equations \cite{Wasow} (see also \S~\ref{sec:For}
for a brief exposition and notations). Generically all the series
$\tilde{\bfy}_k$ are factorially divergent and there is no immediate
way to uniquely associate actual functions to them. Neither can
$\tilde{\bfy}$ be viewed as a {\em classical} asymptotic expansion since
the $\tilde{\bfy}_\bfk$ are  beyond all orders of each
other (e.g., for $\bfk\ne 0$  and all $l\in\NN$, $\mathrm{e}^{-\bflam\cdot\bfk x}x^{\bfk\cdot\bfm}
\tilde{\bfy}_\bfk =o(x^{-l})$).

One question is therefore to understand the relation between these
(algorithmically obtained) formal solutions and the actual solutions of
(\ref{eqor}). In the present paper we show that a suitable
generalization of Borel summation provides a one-to-one correspondence
between transseries and actual solutions of (\ref{eqor}):

\begin{eqnarray}
  \label{symbol}
  {\bfy}\rightleftharpoons\tilde{\bfy}_0+\sum_{\mathbf{k}\ge 0; |\mathbf{k}|>0
}C_1^{k_1}\cdots C_n^{k_n}
\mathrm{e}^{-(\bfk\cdot\bflam) x}x^{\bfk\cdot\bfm}\tilde{\bfy}_{\bfk}
\ \ \ (x\rightarrow\infty, \arg(x)=\xi)
\end{eqnarray}

Given $\bfy$, the value of $C_i$ can change only when
$\xi+\arg(\lambda_i-\bfk\cdot\bflam)=0$, $k_i\in\NN\cup\{0\}$, i.e.
when crossing one of the (finitely many by (c1)) Stokes lines.  The
correspondence (\ref{symbol}) defines a summation method, in the sense
that it is an extension of  convergent summation which
preserves its basic properties: linearity, multiplicativity,
commutation with differentiation and with complex conjugation. These
properties are essential for obtaining true solutions out of
transseries for nonlinear differential equations.  Our procedure is
similar to the medianization proposed by \'Ecalle, but (due to the
structure of (\ref{eqor})) requires substantially fewer analytic
continuation paths. In addition we classify in the context of
(\ref{eqor}) all admissible summation methods (there is a
one-parameter family of them, preserving the properties of usual
summation).  Summation recovers from transseries actual solutions of
(\ref{eqor}) without resorting to (\ref{eqor}) in the process.  In
addition, the analysis reveals a rich analytic structure and formulas
linking the various $\tilde{\bfy}_\bfk$ among themselves (resurgence
relations). In \cite{Costin} we studied this problem under further
restrictions on the transseries (decay of the exponentials in a full
half-plane) and on the differential equation.  Removing those
restrictions creates difficulties that required a new approach. New
resurgence relations are found and in addition we provide a complete
description, needed in applications, of the singularity structure of
the Borel transforms of $\tilde{\bfy}_\bfk$.

\subsubsection{Notes on Borel summation}
\label{sec:bsum}

The following is a very brief description; for more details on classical
Borel summation see \cite{Hardy}, \cite{Borel} and for recent
developments see \cite{Balser} and especially \cite{Ecalle}.

If $\tilde{f}=\sum_{k=0}^{\infty} a_kx^{-k-r}$ is a formal series with
$\Re(r)>0$, its Borel transform is defined as the (still formal)
series
$\mathcal{B}{\tilde{f}}=\sum_{k=0}^{\infty}p^{k-1+r}/\Gamma(k+r+1)$,
the term-by-term inverse Laplace transform of $\tilde{f}$. If
$r\in\NN^+$ and $\tilde{f}$ \emph{converges} (to $f$), then
$\mathcal{B}\tilde{f}$ converges in $\CC$ to an analytic
function which is Laplace ($\lap$) transformable and
$\lap\bor\tilde{f} =f$.  A similar property holds more generally when
$\Re(r)>0$, with now $f$ and $\bor\tilde{f}$ ramified analytic
functions. Even when $\tilde{f}$ is divergent (not faster than
factorially), $\bor \tilde{f}$ may have a nonzero radius of
convergence and define a germ of an analytic function $F(p)$. If
$F(p)$ can be analytically continued along a ray $\arg(p)=\phi$, and
its growth is at most exponential, then
$f=\lap_{\phi}F=\int_{\mathrm{e}^{\mathrm{i}\phi}\RR^+}F(p)\mathrm{e}^{-xp}\mathrm{d}p$
defines a function with the property $f\sim\tilde{f}$ as
$x\rightarrow\infty$ with $\Re(x \mathrm{e}^{\mathrm{i}\phi})>0$.  In
general now $F(p)$ is singular (not only for $p=0$), and
$\mathcal{L}_\phi F$ (when it exists) will depend on $\phi$; the usual
convention is to choose $\phi$ so that
 \begin{equation}\label{convxp}xp\in\RR^+\end{equation}
  Thus, the Borel sum of $\tilde{f}$ in the
direction $x$, if it exists, is defined as
$\lap_{\phi(x)}\mathcal{B}f$ with $-\phi(x)=\xi:=\arg(x)$.  However, when
$\tilde{f}$ is a series with real coefficients, it is a common occurrence that
$F(p)$ is singular for $p\in\RR^+$ (because of conjugation symmetry),
and then the {\em classical} Borel sum of $\tilde{f}$ along the real axis
(the interesting direction in many cases) is undefined.

The difficulty is more serious than it may seem. Summation along paths
that avoid the singularities from above or from below give different
results and thus would lead to an ambiguous (or unnatural) procedure.
More importantly, a ``summation'' procedure using such paths would not
commute with complex-conjugation since the ``sum'' will be, in
general, complex for real $\tilde{f}$ and would thus fail to be a
(proper) summation method. Symmetry considerations suggest a first
step towards {\em averaging}: summation along the half-sum of the two paths
does commutes with complex conjugation. But this solves a problem only to
create another one. The half-sum process fails to commute with
multiplication (of series) and is thus not a summation method, either.

It turns out that there exist more sophisticated averages which have
all the required properties to define a summation procedure.  The
technique of averaging, as well as the fundamental concepts of
analyzable functions and transseries, were discovered and studied by
Ecalle in his constructive approach to the Dulac conjecture (see
\cite{Ecalle-book}, \cite{Ecalle} and \cite{Ecalle2}). The concept of
analyzable function (also discovered by Kruskal in the context of
surreal analysis) is regarded as a very comprehensive generalization
of analyticity/ quasianalyticity. The widely held belief is that
all functions of ``natural origin'' must be analyzable. In particular,
analyzable functions have uniquely associated transseries which are
generalized-Borel summable, after a finite number of transformations
\cite{Ecalle2}.  We show that, in the particular case of
(\ref{eqor}), decreasing solutions are analyzable.

There is a wide class of admissible, all-purpose averaging methods
(\cite{Ecalle3}). As yet there is no unique, natural average and the
problem in its full generality is highly nontrivial.  We obtain the
{\em balanced average} directly from the study of the general solution
of the inverse Laplace transform of (\ref{eqor}). Its potential
nonuniqueness is lifted, in our context, by imposing compatibility
with {\em hyperasymptotics} an important improvement in asymptotic
calculations proposed by M.  Berry (\cite{Berry}, \cite{Berry-hyp},
\cite{Berry-Howls}, \cite{Berry-gamma}).

\subsubsection{Nonresonance}\label{nonres} (1) $\lambda_i,\,i=1,...,n_1$ are
assumed $\ZZ$-linearly independent for any $d$. (2) Let $\theta\in[0,2\pi)$ and
$\tilde{\boldsymbol{\lambda}}=(\lambda_{i_1},...,\lambda_{i_p})$ where
$\left|\arg\lambda_{i_j}-\theta\right|\in(-\pi/2,\pi/2)$ (those
eigenvalues contained in the open half-plane $H_\theta$ centered along
$\mathrm{e}^{\mathrm{i}\theta}$).  We require that
for any $\theta$ the complex numbers in the set
$\{\tilde{\lambda}_{i}-\mathbf{k}
\cdot\tilde{\boldsymbol{\lambda}}\in H_\theta:\bfk\in\NN^p,\,
i=1,...,p\}$
(note: the set is \emph{finite}) have
{\em distinct} directions. These are the Stokes lines $d_{i;\bfk}$.

That the set of $\boldsymbol{\lambda}$ which satisfy (1) and (2) has full
measure follows from the fact that (1) and (2) follow from the condition:

\begin{eqnarray}\label{strongnonr}
  \Big(\bfm,\bfm'\in\ZZ^n ,\ 
\alpha\in\RR\   \mbox{and}\
(\bfm-\alpha\bfm')\cdot\boldsymbol{\lambda}=0
\Big)\ \Rightarrow \Big(\bfm=\alpha\bfm'\Big)
\end{eqnarray}

\z Indeed, if (\ref{strongnonr}) fails, one  of $\Re\lambda_j,
\Im\lambda_j$ is a rational function with rational coefficients of the
other $\Re\lambda_j$ and $\Im\lambda_j$, corresponding to a zero
measure set in $\RR^{2n}$.

\subsection{Further notations and conventions}
\label{sec:anset}

If $y_{1}$ and $y_2$ are inverse Laplace transformable functions, then in a
neighborhood of the origin
$\lap^{-1}(y_1y_2)=(\lap^{-1}y_1)*(\lap^{-1}y_2)$, where for $f,
g\in\lone$ convolution is given by

\begin{eqnarray}\label{defconv}
f*g:=p\mapsto\int_0^p f(s)g(p-s)\mathrm{d}s
\end{eqnarray}

\z We use the convention $\NN\ni 0$.  Let 

\begin{eqnarray}
  \label{defW}
  \mathcal{W}=\left\{p\in\CC:p\ne k\lambda_i\,,\forall
k\in\NN,i=1,2,\ldots,n\right\}
\end{eqnarray}

The directions $d_j=\{p:\arg(p)=\phi_j\}, j=1,2,\ldots,n$ (cf. (a2))
are the {\em Stokes lines} of $\tilde{\bfy}_0$ (note: sometimes known as {\em anti-}Stokes
lines!).  We construct over $\mathcal{W}$ a surface $\mathcal{R}$,
consisting of homotopy classes of smooth curves in $\mathcal{W}$
starting at the origin, moving away from it, and crossing at most one
Stokes line, at most once (see Fig. 1):

\begin{eqnarray}\label{defpaths}
{\cal R}:=\Big\{\gamma:(0,1)\mapsto \mathcal{W}:\ 
\gamma(0_+)=0;\ \frac{\mathrm{d}}{\mathrm{d}t}|\gamma(t)|>0;\
\arg(\gamma(t))\ \mbox{monotonic}\Big\}\cr \mbox{modulo homotopies}\hfill
\end{eqnarray} 

\z Define $\mathcal{R}_1\subset \mathcal{R}$ by (\ref{defpaths}) with
the supplementary restriction $\arg(\gamma)\in(\psi_n-2\pi,\psi_2)$
where $\psi_n=\max\{-\pi/2,\phi_n-2\pi\}$ and
$\psi_2=\min\{\pi/2,\phi_2\}$. $\mathcal{R}_1$ may be viewed as the
part of the covering $\mathcal{R}$, above a sector containing the real
axis. Similarly we let $\mathcal{R'}_1\subset \mathcal{R}_1$
with the restriction that the curves $\gamma$ do not cross
the Stokes directions $d_{i,\bfk}$ (cf. \S\ref{nonres}), other than $\RR^+$,
and we let $\psi_{\pm}=\pm\max
(\pm\arg \gamma)$ with $\gamma\in\mathcal{R'}_1$.

\begin{picture}(0,0)%
\epsfig{file=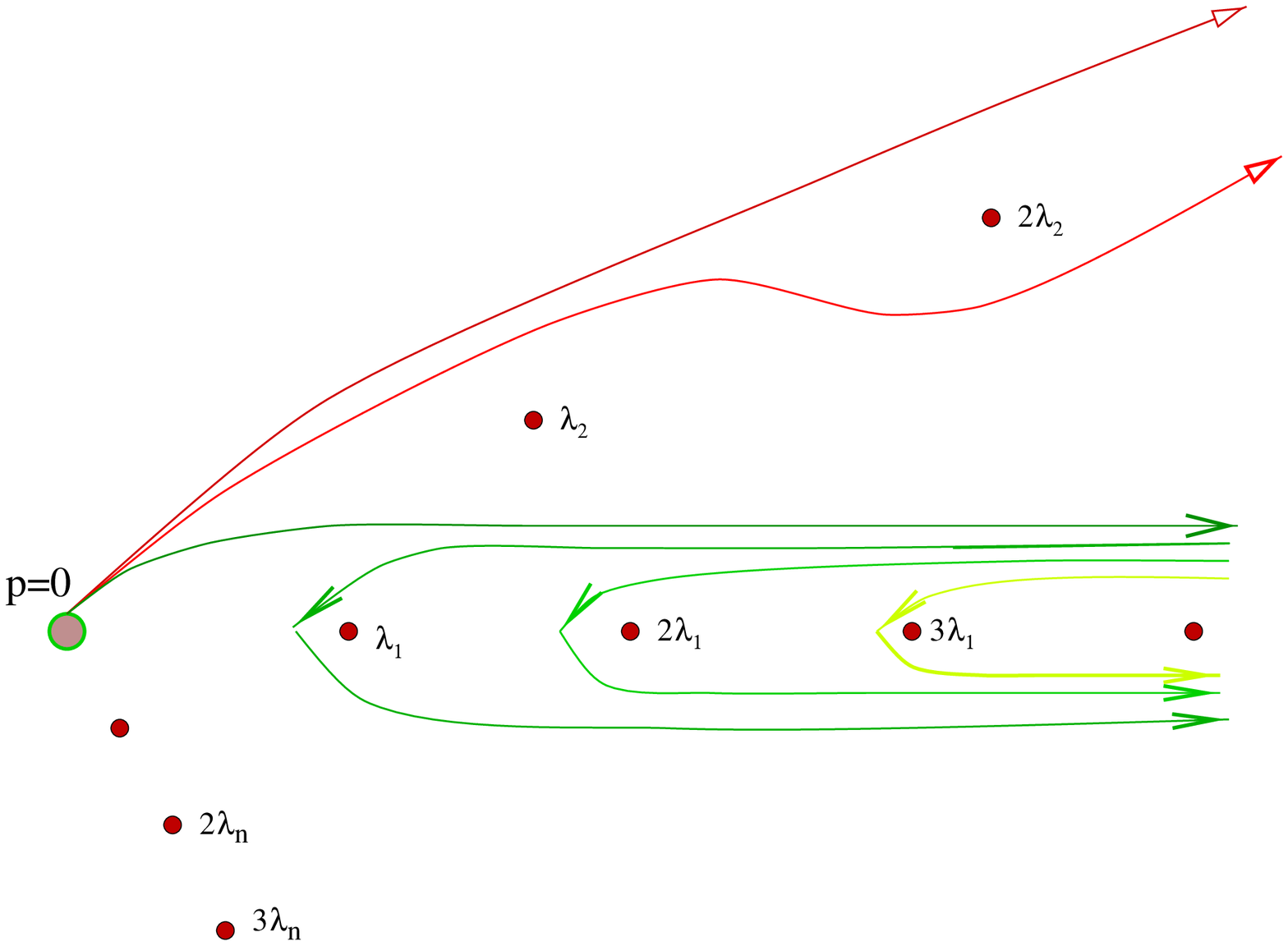, height=7cm}%
\end{picture}%
\setlength{\unitlength}{0.00033300in}%
\begingroup\makeatletter\ifx\SetFigFont\undefined%
\gdef\SetFigFont#1#2#3#4#5{%
  \reset@font\fontsize{#1}{#2pt}%
  \fontfamily{#3}\fontseries{#4}\fontshape{#5}%
  \selectfont}%
\fi\endgroup%
\begin{picture}(10890,7989)(4801,-9310)
\end{picture}

\smallskip
\centerline{{{Fig 1.} \emph{The paths
near $\lambda_2$ belong to $\mathcal{R}$.   }}}
\centerline{\em The paths
near $\lambda_1$ relate to the balanced average}

\bigskip

By $AC_\gamma(f)$ we denote the analytic
continuation of $f$ along a curve $\gamma$.
For the analytic continuations near a Stokes line $d_{i;\bfk}$ we use symbols
similar to \'Ecalle's: $f^-$ is the branch of $f$ along a path
$\gamma$ with $\arg(\gamma)<\phi_i$, while $f^{-j+}$ denotes the
branch along a path that crosses the Stokes line between $j\lambda_i$
and $(j+1)\lambda_i$ (see also \cite{Costin}). 

We use the notations $\mathcal{P}f$ for $\int_0^p f(s)\mathrm{d}s$ and
$\mathcal{P}_\gamma f$ if integration is along the curve $\gamma$.

We write $\bfk\succeq\bfk'$ if $k_i\ge k'_i$ for all $i$ and
$\bfk\succ\bfk'$ if $\bfk\succeq\bfk'$ and $\bfk\ne\bfk'$. The relation
$\succ$ is a well ordering on $\NN^{n_1}$.  We let $\mathbf{e}_j$ be
the unit vector in the $j^{\rm th}$ direction in $\NN^{n_1}$.

Formal expansions are denoted with a tilde, and capital letters
$\bfY,\bfV\ldots$ will usually denote Borel transforms or other
functions naturally
associated to Borel space. For notational convenience, we
will not however distinguish between the series
$\tilde{\bfY}_k=\bor\tilde{\bfy}_\bfk$, which in our case turn out to be
convergent, and the sums $\bfY_\bfk$ of these series  as germs
of ramified analytic functions.

By symmetry (renumbering the directions) it suffices to analyze the
singularity structure of $\bfY_0$ in $\mathcal{R}_1$ only. However,
(c1) breaks this symmetry for $\bfk\ne 0$ and the properties of these
$\bfY_\bfk$ 
will  be analyzed along some other directions as well.

$\bchi_A$ will denote the characteristic function of the set $A$.
We write $|{\bf f}|:=\max_{i}\{|f_i|\}$.
We have

\begin{eqnarray}\label{Taylor series}
{\bf g}(x,{\bfy})=\sum_{|{\bf l}|\ge 1}{\bf g}_{\bf l}(x) {\bf
y}^{\bf l}=\sum_{s\ge 0;|{\bf l}|\ge 1}{\bf g}_{s,\bf l}x^{-s}
{\bfy}^{\bf l} \ \ (|x|>x_0,|\bfy|<y_0)
\end{eqnarray}

\z where ${\bfy}^{\bf l}=y_1^{l_1}\cdots y_n^{l_n}$ and 
$|{\bf l}|=l_1+\cdots+l_n$. By construction ${\bf g}_{s,\bfl}=0$ 
if $|\bfl|=1$ and $s\le M$.

The formal inverse Laplace transform of $\bog(x,\bfy(x))$ (formal since $\bfy$ is still unrestricted)
is given by:

\begin{eqnarray}\label{lapdef}
{\cal L}^{-1}\left(\sum_{|\bf l|\ge 1}
{\bf y}(x)^{\bf l}\sum_{s\ge 0}{\bf g}_{s,\bf l}x^{-s}\right)
=\sum_{|\bf l|\ge 1}{{\bf G}}_{\bf l}*{\bf Y}^{*\bf l}+
\sum_{|\bf l|\ge 2}{\bf g}_{0,\bfl}{\bf Y}^{*\bf l}=:{{\calnb}}(\bf Y)
\cr&&
\end{eqnarray}

\z 
with $ {{\bf G}}_{\bf l}(p)=\sum_{s=1}^{\infty}{\bf g}_{s,\bf
  l}{p^{s-1}}/{s!}$ and $({{\bf
    G_{l}}*\bfY^{*\bfl}})_j:= \left({{\bf
    G}_{\bfl}}\right)_j*Y_1^{*l_1}*..*Y_{n}^{*l_n}$. By 
(n5), ${{\bf G}}^{(l)}_{1,\bfl}(0)=0\ \mbox{ if
  }|\bfl|=1$ and $l\le M$.  The inverse
Laplace transform of (\ref{eqor}) is the convolution equation:

\begin{eqnarray}\label{eqil}
-p{\bf Y}={\bf F}_0-\hat\Lambda {\bf Y}-
\hat B\mathcal{P}\bfY+{\cal  N}({\bf Y})
\end{eqnarray}

Let $\bfdd_\bfj(x):=\sum_{ \bf l\ge j} \binom{\bfl}{\bfj }
\bfgg_\bfl(x) \tilde{\bfy}_0^{\bfl -\bfj }$. Straightforward
calculation (see Appendix \S~\ref{sec:For}; cf. also \cite{Costin})
shows that the components $\tilde{\bfy}_\bfk$ of the transseries
satisfy the hierarchy of differential equations

\begin{eqnarray}
  \label{systemformv}
  &&\bfy'_\bfk+ \left(\hat\Lambda+\frac{1}{x}\left(\hat
  B+\bfk\cdot\bfm\right)-\bfk\cdot\bflam
\right)\bfy_k+\sum_{|\bfj|=1}\bfdd_{\bfj}(x)(\bfy_\bfk)^{\bfj}
=\bft_\bfk
\cr&&
\end{eqnarray}

\z where $\bft_\bfk=\bft_\bfk\left(\bfy_0,\left\{\bfy_{\bfk'}\right\}_{0\prec\bfk'\prec\bfk}\right)$ is a {\em polynomial} in
$\left\{\bfy_{\bfk'}\right\}_{0\prec\bfk'\prec\bfk}$ 
and in  $\{\bfdd_\bfj\}_{\bfj\le\bfk}$ (see (\ref{eqmygen})), with
$\bft(\bfy_0,\emptyset)=0$; $\bft_\bfk$ satisfies the homogeneity
relation

\begin{eqnarray}
  \label{homogeq2}\label{homogeq}
  \bft_\bfk\left(\bfy_0,\left\{C^{\bfk'}\bfy_{\bfk'}\right\}_{0\prec\bfk'\prec\bfk}\right)
=C^{\bfk}\bft_\bfk\left(\bfy_0,\left\{\bfy_{\bfk'}\right\}_{0\prec\bfk'\prec\bfk}\right)
\end{eqnarray}

\z Taking $\lapi$ in (\ref{systemformv}) we get, with $\bfd_\bfj=\sum_{
  \bf l\ge j} \binom{\bfl}{\bfj }\!\left[ \bfg_\bfl*\bfY_0^{*(\bfl
  -\bfj) }+\bfgg_{0,\bfl}*\bfY_0^{*(\bfl -\bfj)}\right]$,

\begin{eqnarray}
  \label{invlapvk}\label{eqMv}\label{invlapyk}
  &&\left(-p+\hat{\Lambda}-\bfk\cdot\bflam\right)\bfY_\bfk
+\left(\hat{B}+\bfk\cdot\bfm\right)
\mathcal{P}\bfY_\bfk+\sum_{|\bfj|=1}\bfd_\bfj*\bfY_\bfk^{*\bfj}
=\bfT_\bfk\cr&&
\end{eqnarray}

\z where $\mathbf{T}_\bfk$ is now a {\em convolution} polynomial, cf.
(\ref{defT}).


\subsection{Main results}
\label{sec:MR}

(a) {\em Analytic structure}.  \begin{Theorem}\label{AS} (i) $\bfY_0=\bor\tilde{\bfy}_0$
is analytic in $\mathcal{R}\cup \{0\}$.

The singularities of $\bfY_0$ (which are contained in the set
$\{l\lambda_j:l\in\NN^+,j=1,2,\ldots,n\}$) are described as follows.
For $l\in\NN^+$ and small $z$, using the notations explained
in \S\ref{sec:anset},

 \begin{multline}
   \label{SY0}
   \bfY_0^{\pm}(z+l\lambda_j)=\pm\Big[(\pm S_j)^l(\ln z)^{0,1}
   \bfY_{l\mathbf{e}_j}(z)\Big]^{(lm_j)}+\bfB_{lj}(z)=\cr
   \Big[z^{l\beta_j'-1}(\ln
   z)^{0,1}\,\bfA_{lj}(z)\Big]^{(lm_j)}+\bfB_{lj}(z) \ (l=1,2,\ldots)
 \end{multline}

 \z where the power of $\ln z$ is one iff
 $l\beta_j\in\ZZ$, and $\bfA_{lj},\bfB_{lj}$ are analytic for small
 $z$. The functions $\bfY_\bfk$ are,  exceptionally, analytic at
 $p=l\lambda_j$, $l\in\NN^+$, iff,

\begin{eqnarray}\label{defSj}
S_j=r_j\Gamma(\beta'_j)\left(\bfA_{1,j}\right)_j(0)=0
\end{eqnarray}

\z where $r_j=1-\mathrm{e}^{2\pi i(\beta'_j-1)}$
if $l\beta_j\notin\ZZ$ and $r_j=-2\pi i$
otherwise. The $S_j$ are Stokes constants,
see Theorem~\ref{Stokestr}.

(ii) $\bfY_\bfk=\bor{\tilde{\bfy}}_\bfk$, $|\bfk|>1$, are analytic in
$\mathcal{R}\backslash
\{-\bfk'\cdot{\bflam}+\lambda_i:\bfk'\le\bfk,1\le i\le n\}$.  For
$l\in\NN$ and $p$ 
 near $l\lambda_j$, $j=1,2,\ldots,n$ there exist $\bfA=\bfA_{\bfk jl}$ and
$\bfB=\bfB_{\bfk jl}$ analytic at zero so that ($z$ is as above)

 \begin{multline}
   \label{SYK}
  \bfY_\bfk^{\pm}(z+l\lambda_j)=
\pm\Big[(\pm S_j)^l\binom{k_j+l}{
     l}(\ln z)^{0,1}
   \bfY_{\bfk+l\mathbf{e}_j}(z)\Big]^{(lm_j)}+\bfB_{\bfk lj}(z)=\cr
  \Big[z^{\bfk\cdot{\bfbet}'+l\beta_j'-1}(\ln z)^{0,1}\,\bfA_{\bfk l
    j}(z)\Big]^{(lm_j)}+\bfB_{\bfk l j}(z)\ (l=0,1,2,\ldots)
       \end{multline}

\z where the power of $\ln z$ is $0$ iff $l=0$
or 
$\bfk\cdot{\bfbet}+l\beta_j-1\notin\ZZ$ and $\bfA_{\bfk 0 j}=\mathbf{e}_j/\Gamma(\beta'_j)$.
Near $p\in\{\lambda_i-\bfk'\cdot{\bflam}:0\prec\bfk'\le\bfk\}$, (where
$\bfY_0$ is analytic) $\bfY_\bfk,\,\bfk\ne 0$ have convergent
Puiseux series.

\end{Theorem} 

{\sc Remark:} The fact that the singular part of
$\bfY_\bfk(p+l\lambda_j)$ in (\ref{SY0}) and (\ref{SYK}) is a {\em
  multiple} of $\bfY_{\bfk+l\mathbf{e}_j}(p)$ is the effect
of  {\em resurgence} and  provides a way of
determining
the $\mathbf{Y}_\bfk$ given $\mathbf{Y}_0$ 
provided the $S_j$ are nonzero. Since, generically,
the $S_j$ are nonzero this is a  surprising upshot:  given one formal
solution, (generically) an $n$ parameter family of solutions can be
constructed out of it, without using (\ref{eqor}) in the process;
the differential equation itself is then recoverable (\cite{Costin3}).

\bigskip

By Theorem~\ref{AS} the Borel transforms
$\bfY_\bfk={\bor}\tilde{\bfy}_\bfk$ define germs of ramified analytic
functions and are continuable on the surface $\mathcal{R}$. In order
to be able to take Laplace transforms we need to define
$\bor\tilde{\bfy}_\bfk$ along any direction $d$ in $\mathcal{S}$. If
$d\ne d_{j,\bfk}$ then $\bfY_\bfk$ are analytically continuable along $d$ and
the continuations turn out to have all the properties that we need.
But along Stokes lines $d_{j,\bfk}$ analytic continuation is impossible:
in general the functions $\bfY_\bfk$ have an infinite array of
branch points (\ref{SYK}). In addition, while both $\bfY_\bfk^+$ and
$\bfY_\bfk^-$ turn out to be Laplace transformable (in distributions)
along $d_{j,\bfk}$, $\lap\bfY_\bfk^+$ and $\lap\bfY_\bfk^-$ are generically
{\em different}. Neither the upper nor the lower continuation would
give rise to a definition of Borel summation which commutes with
complex-conjugation, as discussed in the introduction.  Other analytic
continuations along paths $\gamma$ that {\em cross} $d_{j,\bfk}$ have even
worse problems, namely that $AC_\gamma(\bfY_\bfk*\bfY_\bfk)\ne
AC_\gamma(\bfY_\bfk)*_\gamma AC_\gamma(\bfY_\bfk)$, (see
\cite{Costin}; for notations only, see also \S\ref{sec:twolemmas}). As $\bor$ transforms differential equations into
convolution equations, the implication is that with such a $\gamma$,
$\lap AC_\gamma(\bfY_\bfk)$ would {\em not} be, in general, solutions
of their differential equations. Individual analytic continuations are
thus inadequate for solving (\ref{eqor}), but some
{\em averages} of analytic continuations do satisfy all the
requirements.  Let $\alpha=1/2+\mathrm{i}\sigma$ with $\sigma\in\RR$
and $\bor\tilde{\bfy}_\bfk$ be extended along $d_{j,\bfk}$ by
the  weighted average of analytic continuations

\begin{eqnarray}
\label{defmed}
\bor_\alpha\tilde{\bfy}_\bfk=\bfY_\bfk^{\alpha}=
\bfY_\bfk^++\sum_{j=1}^{\infty}\alpha^j\left(\bfY_\bfk^{-}
-\bfY_\bfk^{-({j-1})+}\right)
\end{eqnarray}

\begin{Remark}\label{RF}  Relation (\ref{defmed}) gives the most general reality preserving,
linear operator mapping formal power series solutions of
(\ref{eqor}) to solutions of
  (\ref{eqil}) in distributions (more precisely
in $\mathcal{D}'_{m,\nu}$; see \S\ref{sec:NC}).
\end{Remark}

\z This remark follows easily from Proposition~\ref{SRY0} and 
Theorem~\ref{RE} below.

The choice $\alpha=1/2$ has special properties; we call
$\bor_{\frac{1}{2}}\tilde{\bfy}_\bfk=\bfY^{ba}_\bfk$ the balanced average
of $\bfY_\bfk$. For this choice the expression (\ref{defmed})
 coincides with the one in which $+$ and $-$ are interchanged (Proposition~\ref{medianization}), accounting for the
reality-preserving property. Clearly, if $\bfY_\bfk$ is analytic along
$d_{j,\bfk}$, then the terms in the infinite sum vanish and
$\bfY_\bfk^{\alpha}=\bfY_\bfk$; we also let $\bfY_\bfk^{\alpha}=\bfY_\bfk$ if
$d\ne d_{j,\bfk}$, where again $\bfY_\bfk$ is analytic. It follows from
(\ref{defmed}) and Theorem~\ref{CEQ} below that the Laplace integral
of $\mathbf{Y}^{\alpha}_\bfk$ along $\RR^+$ can be deformed into contours as
those depicted in Fig. 1, with weight $-(-\alpha)^{k}$ for a contour
turning around $k\lambda_1$.

 In addition to symmetry (the balanced average equals the half sum
of the upper and lower continuations on $(0,2\lambda_j)$,
\cite{Costin3}),  an
asymptotic property uniquely picks $C=1/2$. Namely, for $C=1/2$
alone are the $\lap\bor\tilde{\bfy}_\bfk$ always {\em summable to the
  least term} cf. \S~\ref{sec:aver}.

(b) {\em Connection with (\ref{eqor}) and (\ref{homogeq})}.
Generalized Borel summation coincides with the usual Borel summation
when the transseries consists of only one term, the first series, when
that series is classically Borel summable. This is clear from
theorem~\ref{CEQ} (ii) below. Furthermore, generalized summation is a
map from a class of formal series to functions which is linear,
multiplicative, commutes with differentiation and complex conjugation
(cf. \S~\ref{sec:aver}, \S~\ref{sec:NC}) so it is a summation
procedure, which furthermore, establishes along every direction a one
to one correspondence between transseries and decaying actual solutions of
(\ref{eqor}) cf. \S~\ref{sec:aver}, Proposition~\ref{medianization}
and Theorem~\ref{CEQ} below.

For
clarity we state the results for $x\in S_x$, a sector in the
right half plane containing
$\lambda_1=1$ in which (c1) holds and for $p$ in the associated domain
$\mathcal{R'}_1$, but $\lambda_1$ plays no special role as discussed
in the introduction.

\begin{Theorem}\label{CEQ} (i) The branches of $(\bfY_\bfk)_\gamma$ in
  $\mathcal{R'}_1$ ($\mathcal{R}_1$ if $\bfk=0$)
have limits in a $C^*$-algebra of distributions,
$\mathcal{D}'_{m,\nu}(\RR^+)\subset\mathcal{D}'$ (cf. \S~\ref{sec:NC}).
Their Laplace transforms in $\mathcal{D}'_{m,\nu}(\RR^+)$
$\lap(\bfY_\bfk)_\gamma$ exist simultaneously and with
$x\in\mathcal{S}_x$ and for any $\delta>0$ there is a constant $K$ and
an $x_1$ large enough, so that for $\Re(x)>x_1$ we have
$\left|\lap(\bfY_\bfk\right)_\gamma(x)|\le K\delta^{|\bfk|}$.

In addition, $\mathbf{Y}_\bfk(p\mathrm{e}^{\mathrm{i}\phi})$ are
continuous in $\phi$ with respect to the $\mathcal{D}'_{m,\nu}$
topology, (separately) on $(\psi_-,0]$ and $[0,\psi_+)$ . 

If
$m>\max_i(m_i)$ and $l<\min_i |\lambda_i|$ then
$\mathbf{Y}_0(p\mathrm{e}^{\mathrm{i}\phi})$ is continuous in
$\phi\in[0,2\pi]\backslash\{\phi_i:i\le n\}$ in the
$\mathcal{D}'_{m,\nu}(\RR^+,l)$ topology and has (at most) jump
discontinuities for $\phi=\phi_i$. For each $\bfk$, $|\bfk|\ge 1$ and
any $K$ there is an $l>0$ and an $m$ such that
$\mathbf{Y}_\bfk(p\mathrm{e}^{\mathrm{i}\phi})$ are continuous in
$\phi\in[0,2\pi]\backslash\{\phi_i; -\bfk'\cdot\bflam+\lambda_i:i\le n
,\bfk'\le\bfk\}$ in the $\mathcal{D}'_{m,\nu}((0,K),l)$ topology and
have (at most) jump discontinuities on the boundary.

(ii) The sum (\ref{defmed}) converges in $\mathcal{D}'_{m,\nu}$ (and
coincides with the analytic continuation of $\bfY_\bfk$ when
$\bfY_\bfk$ is analytic along $\RR^+$). For any $\delta$ there is a
large enough $x_1$ {\em independent of $\bfk$} so that
$\bfY^{ba}_\bfk(p)$ with $p\in\mathcal{R}'_1$ are Laplace transformable
in $\mathcal{D}'_{m,\nu}$
for $\Re(xp)>x_1$ and furthermore $|(\lap\bfY^{ba}_\bfk)(x)|\le
\delta^{|\bfk|}$. In addition, if $d\ne\RR^+$, then for large $\nu$,
$\bfY_\bfk\in L^1_\nu(d)$.

The functions $\lap\bfY_\bfk^{ba}$ are analytic for
$\Re(xp)>x_1$. For any $\bfC\in\CC^{n_1}$ there is an $x_1(\bfC)$
large enough so that the sum

\begin{eqnarray}
  \label{soleqn}
  \bfy=\lap\bfY_0^{ba}+\sum_{|\bfk|> 0}\bfC^{\bfk}\mathrm{e}^{-\bfk\cdot\bflam
    x}x^{-\bfk\cdot\bfbet}\lap\bfY_\bfk^{ba}
\end{eqnarray}

\z converges uniformly for $\Re(xp)>x_1(\bfC)$, and $\bfy$ is a solution
of (\ref{eqor}). When the direction
of $p$ is not the real axis then, by definition, $\bfY^{ba}_\bfk=\bfY_\bfk$,
$\mathcal{L}$ is the usual Laplace transform 
and (\ref{soleqn}) becomes

\begin{eqnarray}
  \label{soleqnpm}
  \bfy=\lap\bfY_0+\sum_{|\bfk|> 0}\bfC^{\bfk}\mathrm{e}^{-\bfk\cdot\bflam
    x}x^{-\bfk\cdot\bfbet}\lap\bfY_\bfk
\end{eqnarray}

In addition, $\lap\bfY_\bfk^{ba}\sim \tilde{\bfy}_\bfk$ for large $x$
in the half plane $\Re(xp)>x_1$, for all $\bfk$, uniformly.

iii) More generally, for any $\alpha$ and any 
solution $\bfy$ of (\ref{eqor}) such that $\bfy\sim \tilde{\bfy}_0$ for
large $x$ along a ray
in $S_x$ there exists a constant
vector $\bfC=\bfC_{\alpha;\bfy}$ so that

\begin{eqnarray}
  \label{soleqnpa}
  \bfy=\lap\bor_\alpha\tilde{\mathbf{y}}_0+\sum_{|\bfk|> 0}\bfC^{\bfk}\mathrm{e}^{-\bfk\cdot\bflam
    x}x^{-\bfk\cdot\bfbet}\lap\bor_\alpha\tilde{\mathbf{y}}_\bfk
\end{eqnarray}

Given $\alpha$ the representation (\ref{soleqnpa}) of $\bfy$ is unique
(see also \S~\ref{sec:bsum} above for the convention on the direction
of Laplace integration).

\end{Theorem}

Of special interest are
the cases $\alpha=1/2$, discussed above, and also  $\alpha=0,1$ which give:

\begin{eqnarray}
  \label{soleqnp}
  \bfy=\lap\bfY_0^{\pm}+\sum_{|\bfk|> 0}\bfC^{\bfk}\mathrm{e}^{-\bfk\cdot\bflam
    x}x^{-\bfk\cdot\bfbet}\lap\bfY_\bfk^{\pm}
\end{eqnarray}

{\em (c) Resurgence properties; local Stokes phenomenon}. 

It turns out that the formal series $\tilde{\bfy}_\bfk$ are connected
among each-other via their Borel transforms. Resurgence
formulas link $\bfY_\bfk$ to analytic continuations of
$\bfY_{\bfk'}$ with $\bfk'\prec\bfk$, in a way that, generically,
$\bfY_0$ contains enough information to compute all $\bfY_\bfk$.

Various resurgence properties have been observed in different
contexts, and the term resurgence
has been used with slightly different interpretations.  In the
hyperasymptotic theory of M. Berry, it was discovered that the first
asymptotic series reappears in various shapes in the process of
computing higher terms of the expansions.  J. \'Ecalle, in his
comprehensive theory of analyzable functions, has obtained a general
resurgence principle, the {\em bridge equation} \cite{Ecalle}. The common
denominator of resurgence is the reappearance of ``earlier'' terms in the
formulas of ``later'' ones.  It turns out that,
for our problem, resurgence is fundamentally linked to the {\em
  Stokes phenomenon}. In the following formulas we make the convention
 $\bfY_\bfk(p-j)=0$ for $p<j$ as an element
of $\mathcal{D}'_{m,\nu}(\RR^+)$. We again state
the results 
is stated  for $p\in\mathcal{R}'_1$ and
$x\in S_x$ 
but  hold in any sector where (c1) is valid.

\begin{Theorem}\label{RE}
i) For all $\bfk$ and $\Re(p)>j,\Im(p)>0$ as well
as in $\mathcal{D}'_{m,\nu}$ we have

\begin{eqnarray}
  \label{mainresur}
  \bfY_{\bfk}^{\pm j\mp}(p)-\bfY_{\bfk}^{\pm (j-1) \mp}(p) = (\pm S_1)^j\binom{k_1+j}{j}
  \left(\bfY^\mp_{\bfk+j\mathbf{e}_1}(p-j)\right)^{(mj)}
\end{eqnarray}

\z and also, 

\begin{eqnarray}
  \label{thirdresu}
  \mathbf{Y}_{\bfk}^\pm=\bfY_\bfk^{\mp}+\sum_{j\ge 1} \binom{j+k_1}{ k_1}(\pm S_1)^{j}(\mathbf{Y}^\mp_{\bfk+j\mathbf{e}_1}(p-j))^{(mj)}
\end{eqnarray}

ii) {\em Local Stokes transition.}

\z Consider the expression of a fixed solution $\bfy$ of (\ref{eqor})
as a Borel summed transseries (\ref{soleqn}). As $\arg(x)$ varies,
(\ref{soleqn}) changes only through $\bfC$, and that change occurs
when Stokes lines are crossed (cf. \S\ref{nonres}; the Stokes lines
of $\bfY_0$ are the directions of $\lambda_i$). We have, in the neighborhood
of $\RR^+$, with $S_1$
defined in (\ref{defSj}):

\begin{eqnarray}
  \label{microsto}
  \bfC(\xi)=\left\{\begin{array}{ccc} &\bfC^-=\bfC(-0)\qquad\mbox{for
    $\xi<0$}\\ &\bfC^0=\bfC(-0)+\frac{1}{2}S_1\mathbf{e}_1\qquad\mbox{for
    $\xi=0$}\\ &\bfC^+=\bfC(-0)+S_1\mathbf{e}_1\qquad\mbox{for
    $\xi>0$}\\ &\end{array}\right.
\end{eqnarray}

\end{Theorem}

\z {\em (d) Classical Stokes phenomena and local Stokes transitions}.
Again we formulate the result below for $\lambda_1$ but with
straightforward adjustments it holds relative to any other eigenvalue.
Let $\mathbf{C}$ be of the form
$C_1\mathbf{e}_1$.  Along the imaginary axis, condition
(c1) fails. The positive and negative imaginary are the {\em
  antistokes} lines corresponding to $\lambda_1=1$ (note: sometimes
called Stokes lines!).  If we choose paths in the right half plane
approaching the positive/negative imaginary axis in such a way that
$|x^{-\beta_1-l}\mathrm{e}^{-x}|\rightarrow K\ne 0$ along them, where
$l+\beta\in(0,M)$, then $\bfy\sim 
C^{\pm} x^{-l-\beta_1}\mathrm{e}^{-x}+\tilde{\bfy}_0$ for large $x$ and the term
multiplied by $K$ is now the {\em leading} behavior of $\bfy$. The
particular choice of $K$ and $l$ within this range is rather
arbitrary, the main point being that along such special curves, the
constant $\bfC$ is definable in terms of {\em classical} asymptotics.
Within the right half plane, it is only near the imaginary axis that
this happens, since otherwise the exponential term is smaller than all
terms of $\tilde{\bfy}_0$. On the other hand Borel summation makes
possible the definition of $\bfC$ throughout the right half plane, and we now
address the issue of the relation between classical asymptotics and
exponential asymptotics.

\begin{Theorem}\label{Stokestr} Let $\gamma^{\pm}$ be two paths in the right half plane,
near the positive/ negative imaginary axis such that
$|x^{-\beta_1+1}e^{-x\lambda_1}|\rightarrow 1$ as $x\rightarrow\infty$
along $\gamma^{\pm}$. Consider the solution $\bfy$ of (\ref{eqor})
given in (\ref{soleqn}) with $\mathbf{C}=C\mathbf{e}_1$ and where
the path of integration is $p\in\RR^+$. Then 

\begin{gather}\label{classicS}
\bfy=
(C\pm\frac{1}{2}S_1)\mathbf{e}_1 x^{-\beta_1+1}e^{-x\lambda_1}(1+o(1))
\end{gather}

\z for large $x$ along $\gamma^{\pm}$, where $S_1$ is the same as in
(\ref{defSj}), (\ref{microsto}).

\end{Theorem}

Classical asymptotics loses track of the value of $\bfC$ along any
ray other than the imaginary directions, as the terms multiplied by
$\bfC$ will be hidden ``beyond all orders'' of the classically divergent series
$\tilde{\bfy}_0$. In contrast to the classical picture, we see that
through generalized Borel summation the constant $\bfC$ is precisely
defined throughout the positive half-plane and the question of where
the change in $\bfC$ occurs is well defined.  

Formula (\ref{microsto}) is the exponential asymptotic expression of
the Stokes phenomenon. It shows that the constant jumps as the Stokes
line is crossed, (\ref{microsto}), as originally predicted by Stokes
himself \cite{Stokes}. Subsequently, the original ideas of Stokes,
based on optimal truncation of series were greatly refined by M.
Berry, leading to his theory of hyperasymptotic expansions and a
description of Stokes transitions for saddle integrals \cite{Berry}.

If more than one component of $\bfC$ is nonzero, then in general there
is no direction along which $\bfC$ can be defined through {\em
  classical} asymptotics. Part of the difficulty of studying nonlinear
Stokes phenomena using classical tools stems from this fact.

Relation (\ref{microsto}) expresses the evolution of $\bfC$ and the
presence of a Stokes phenomenon beyond all orders of Poincar\'e
asymptotics.

\section{Proofs and further results}
\label{sec:proofs}

The layout of the proofs is as follows. We study the formal inverse
Laplace transforms of (\ref{eqor}) and (\ref{systemformv})
in a $C^*$-algebra of distributions that we introduce.
 Using a fixed point principle we find the general solution
of these convolution equations, and then study their properties. We
then show that generalized Laplace transforms of these distributions
exist and have all the required properties. The resurgence formulas
are obtained by comparing different expressions for the same solution
of (\ref{eqor}) near a Stokes line.


Since the proofs rely to a large extent on the detailed study of the
convolution equations (\ref{eqil}), (\ref{eqMv}), we start with a few
heuristic remarks. In a convolution equation such as
(\ref{eqil}), the term $(-p+\hat{\Lambda})\bfY$  plays a role similar to that of the highest derivative term
in a differential equation. To illustrate this, assume a solution
$\bfY$ is already given on an interval, say $(0,a)$, and we wish to
extend it to $(0,a+\epsilon)$. We look for such a solution in the form
$\bfY+ \boldsymbol\delta$, where we take $\bfY=0$ on
$(a,a+\epsilon)$ and $\boldsymbol\delta=0$ on $(0,a)$. If $\epsilon$
is small, then $\boldsymbol\delta*\boldsymbol\delta=0$ and the
equation in $\boldsymbol\delta$ is linear inhomogeneous. The terms
that involve integrals of $\boldsymbol\delta$ are of order
$O(\epsilon)\|\boldsymbol\delta\|$ as $\epsilon\rightarrow 0$, so that the
dominant terms are the forcing term, together with 
$(-p+\hat{\Lambda})\boldsymbol\delta$ provided the coefficient is invertible. If in addition
the forcing is non-singular then 
$\boldsymbol\delta$ can be found, e.g., by a convergent $\epsilon$
expansion; this is the analog of an ordinary 
point of a differential equation.  $\boldsymbol \delta$ can be  singular if either
$\Delta_p=\mbox{det\,}(p-\hat{\Lambda})=0$
 or the forcing is singular. To
understand
the qualitative behavior  near a zero of $\Delta_p$ one has to keep
(at least) one more term, the leading term among those previously
discarded (i.e. the second term in the notation (\ref{eqMv})).  In
this approximation, $\boldsymbol\delta$ satisfies  a differential
equation.

In our problem, there are $n$ roots of $\Delta_p$ but because of the
nonlinearity and nonlocality of the equations, a singularity generates
(when convolved with itself) a whole array of singularities affecting
$\boldsymbol\delta$ through the forcing term.  Through convolution, a
nonintegrable singularity produces further singularities of even lower
regularity. We introduce a distribution space,
$\mathcal{D}'_{m,\nu,l}$, whose degree of regularity decreases with
the distance from the origin, at a ``convolution-like'' rate; these
distributions form a convolution Banach algebra (cf.
\S~\ref{sec:NC}).

Technically, the proofs rely on suitable fixed-point theorems in
spaces having some of the properties we want to prove, notably in
terms of regularity and behavior at infinity.  This is combined with a
local analysis near noninvertibility points of the dominant term,
which is done by treating the convolution equations as a 
perturbation of the approximating differential equation mentioned
above, which splits the singularity thus making again possible the use
of fixed point theorems. This analysis is used in order to find the
resurgence properties, which in turn are used to prove (using among
others Lemma~\ref{cuteCauchy} below) the sharper results on global
analyticity and structure of singularities.

 We start by introducing some useful functional spaces
and derive specific fixed point theorems.

\subsection{Technical constructions and results}\label{sec:Tec}

\subsubsection{Focusing spaces and algebras}

We say that a family of norms $\|\|_{\nu}$ depending on a parameter
$\nu\in\RR^+ $ is {\bf focusing} if for any $f$ with $\|f\|_{\nu_0}<\infty$

\begin{eqnarray}
  \label{focus-pocus}
  \|f\|_\nu\downarrow 0 \mbox{ as }\nu\uparrow\infty
\end{eqnarray}

Let $\mathcal{E}$ be a linear space and $\{\|\|_\nu\}$ a family of
norms satisfying (\ref{focus-pocus}).  For each $\nu$ we define a
Banach space $\mathcal{B}_\nu$ as the completion of
$\{f\in\mathcal{E}:\|f\|_{\nu}<\infty\}$. Enlarging $\mathcal{E}$ if
needed, we may assume that $\mathcal{B}_\nu\subset\mathcal{E}$. For
$\alpha<\beta$, (\ref{focus-pocus}) shows that the identity is an
embedding of $\mathcal{B}_\alpha$ in $\mathcal{B}_\beta$. Let
$\mathcal{F}\subset\mathcal{E}$ be the projective limit of the
$\mathcal{B}_\nu$.  That is to say

\begin{eqnarray}
  \label{focuproj}
  \mathcal{F}:=\bigcup_{\nu>0}\mathcal{B}_\nu
\end{eqnarray}

\z is endowed with the topology  in which a sequence is convergent if it
converges in {\em some} $\mathcal{B}_\nu$. We call $\mathcal{F}$
a {\bf focusing space}.

Consider now the case when
$\left(\mathcal{B}_{\nu},+,*,\|\|_\nu\right)$ are commutative Banach
algebras. Then $\mathcal{F}$ inherits a structure of a commutative
algebra, in which $*$ (``convolution'') is continuous. We say that
$\left(\mathcal{F},*,\|\|_\nu\right)$ is a {\bf focusing algebra}.

\subsubsection{Examples}
\label{sec:NC}

Let $K\in\RR^+$ and
$\mathcal{S}=\mathcal{S}_{K,\alpha_1,\alpha_2}=\{p:\arg(p)\in[\alpha_1,\alpha_2]\subset
(-\pi/2,\pi/2),|p|\le K\}$ (or a finite union of such sectors)
and $\mathcal{V}$ be a small neighborhood
of the origin. $\overline{\mathcal{V}}$ will be the closure of
$\mathcal{V}$, cut along the negative axis, and together with these
 upper and
lower cuts.

 {\bf (1)}. $\lone_\nu(\mathcal{K})$. Let
 $\mathcal{K}=\mathcal{S}_{K,\phi,\phi}$.
 The space
$\lone_\nu(\mathcal{K})$ with
 convolution (\ref{defconv}) is a commutative
Banach algebra under each  of the (equivalent) norms

\begin{eqnarray}
  \label{norm00}
  \|f\|_\nu=\int_0^K \mathrm{e}^{-\nu t}|f(t\exp(\mathrm{i}\phi))|\mathrm{d}t
\end{eqnarray}

\z Indeed, with $F(s):=f(s\mathrm{e}^{\mathrm{i}\phi})$ and $G(s):=g(s\mathrm{e}^{\mathrm{i}\phi})$ we have:

\begin{multline}\label{prfRem}
\int_0^{K}\mathrm{d}t\mathrm{e}^{-\nu t}\left|\int_0^{t}\mathrm{d}sF(s)G(t-s)\right|\le
\int_0^{K}\mathrm{d}t\mathrm{e}^{-\nu t}\int_0^{t}\mathrm{d}s|F(s)G(t-s)|\cr
=\int_0^{K}\int_{0}^{K-v}\mathrm{e}^{-\nu (u+v)}|F(v)||G(u)|\mathrm{d}u\mathrm{d}v\cr\le
\int_0^{K}\int_{0}^{K}\mathrm{e}^{-\nu (u+v)}|F(v)||G(u)|
\mathrm{d}u\mathrm{d}v=\|f\|_{\nu }\|g\|_{\nu }
\end{multline}

\z By dominated convergence
$\|f\|_\nu\downarrow 0$ as $\nu\uparrow\infty$
and thus $\lone(\mathcal{K})$ is a focusing algebra.

{\bf (2)} If $K=\infty$ in example $(1)$, then the norms (\ref{norm00})
are not equivalent anymore for different $\nu$, but convolution is
still continuous in (\ref{norm00}) and  the projective limit of
the $L^1_\nu(\RR^+ \mathrm{e}^{\mathrm{i}\phi})$, $\mathcal{F}(\RR^+
\mathrm{e}^{\mathrm{i}\phi})\subset L^1_{loc}(\RR^+ \mathrm{e}^{\mathrm{i}\phi})$, is a focusing algebra.

{\bf (3a)} $\mathcal{T}_\beta(\mathcal{S}\cup\overline{\mathcal{V}})$.
For $\Re(\beta)> 0$ and $\phi_1\ne\phi_2$,
this space is given by $\{f:f(p)=p^{\beta}F(p)\}$,
where $F$ is analytic in the interior
of $\mathcal{S}\cup\mathcal{V}$ and continuous in
its closure. We take the family of (equivalent) norms

\begin{eqnarray}
  \label{normF1}
  \|f\|_{\nu,\beta}=K\sup_{s\in
   \mathcal{S}\cup\overline{\mathcal{V}}}\left|\mathrm{e}^{-\nu p}f(p)\right|
\end{eqnarray}

\z It is clear that convergence of $f$ in $\|\|_{\nu,\beta}$ implies
uniform convergence of $F$ on compact sets in
$\mathcal{S}\cup\mathcal{V}$ (for $p$ near zero, this follows from
Cauchy's formula).  $\mathcal{T}_{\beta}$ are thus Banach spaces and
 focusing spaces by (\ref{normF1}). The spaces
$\{\mathcal{T}_{\beta}\}_{\beta}$ are isomorphic to each-other.
Taking $s=pt$ in (\ref{defconv}) we find that

\begin{gather}
  \label{defconv31}
 p^{-\beta_1-\beta_2-1} (f_1*f_2)(p)=
\int_0^1 t^{\beta_1}F_1(pt) (1-t)^{\beta_2}
F_2(p(1-t))\mathrm{d}t=F(p)
\end{gather}

\z where $F$ is manifestly analytic, and that the application
\begin{eqnarray}
  \label{convodom1}
  (\cdot *\cdot):\mathcal{T}_{\beta_1}\times
\mathcal{T}_{\beta_2}\mapsto \mathcal{T}_{\beta_1+\beta_2+1}
\end{eqnarray}

\z is continuous:

\begin{eqnarray}
  \label{normconvo1}
  &&\|f_1*f_2\|_{\nu,\beta_1+\beta_2+1}=
  K\sup_p\left|\mathrm{e}^{-\nu
    p}\int_0^ps^{\beta_1}F_1(s)(p-s)^{\beta_2}F_2(p-s)\mathrm{d}s\right|\cr&&\le 
K^{-1}\sup_p\int_0^p\left|K\mathrm{e}^{-\nu s}s^{\beta_1}F_1(s)K\mathrm{e}^{-\nu (p-s)}
(p-s)^{\beta_2}F_2(p-s)\right|\mathrm{d}|s|\cr&&
\le \|f_1\|_{\nu,\beta_1}\|f_2\|_{\nu,\beta_2}
\end{eqnarray}

A natural generalization of $\mathcal{T}_\beta$ is obtained taking
$\beta_1,\ldots,\beta_N\in\CC$ with positive real parts, no two of
them differing by an integer. If $f_\beta=\sum_{i=1}^k p^{\beta_i}A_i(p)$
with $A_i$ analytic, then $f_\beta\equiv 0$ iff $A_i\equiv 0$ for all $i$
(e.g., by a  Puiseux series argument).  It is then
natural to identify the space $\mathcal{T}_{\{\beta_1,\ldots,\beta_k\}}$
of functions of the form $f_\beta$ with
$\oplus_{i=1}^k \mathcal{T}_{\beta_i}$. Convolution with analytic
functions is defined on $\mathcal{T}_{\{\beta_1,\ldots,\beta_k\}}$ while
convolution of two functions in $\mathcal{T}_{\{\beta_1,\ldots,\beta_k\}}$
takes values in $\mathcal{T}_{\{\beta_i+\beta_j \,\mbox{mod 1}\} }$.
We write $\mathcal{T}_{\{\cdot\}}$ when the concrete
values of $\beta_1,\ldots,\beta_k$  do not matter.

{\bf (3b)} A particular case of the preceding example is
$\mathcal{A}_{z,l}(\mathcal{S\cup \mathcal{V}})$ consisting of
analytic functions in the interior of $\mathcal{S}\cup\mathcal{V}$,
continuous on its closure, and vanishing at the origin together with
the first $l$ derivatives. $\mathcal{A}_{z,l}$ can be identified with
$\mathcal{T}_l$.

{\bf (4)} $\mathcal{D}'_{m,\nu}$, the ``staircase distributions''.
{\em Proofs} of the properties stated in this paragraph and more details are
given in \S~\ref{starca}.  Let $\mathcal{D}(0,x)$ be the test
functions on $(0,x)$ and $\mathcal{D}=\mathcal{D}(0,\infty)$. Let
$\mathcal{D}'_m\subset\mathcal{D}'$ be the distributions $f$ for which
$f=F_k^{(km)}$ on $\mathcal{D}(0,k+1)$ with $F_k\in\lone[0,k+1]$.
There is a uniquely associated ``staircase decomposition'', a sequence
$\left\{\Delta_i(f)\right\}_{i\in\NN}=\left\{\Delta_i\right\}_{i\in\NN}$
such that $\Delta_i\in\lone(\RR^+)$,
$\Delta_i=\Delta_i\bchi_{[i,i+1]}$ and

\begin{eqnarray}
  \label{stdec0}
  f=\sum_{i=0}^{\infty}\Delta_i^{(mi)}
\end{eqnarray}

\z Convolution is well defined on $\mathcal{D}'_m$ by

\begin{eqnarray}
  \label{formulaforconv}
   \label{stdeccomv}
  f*\tilde{f}:=(F_k*\tilde{F}_k)^{(2km)}\ \ \mbox{in } \mathcal{D}'(0,k+1)
\end{eqnarray}

\z and $(\mathcal{D}'_m,+,*)$ is a commutative algebra. We define, for
$f\in\mathcal{D}'_{m}$,

\begin{eqnarray}
  \label{normdistr00}
   \|f\|_{\nu ,m}:=c_m\sum_{i=0}^{\infty}\nu ^{im}\|\Delta_i\|_\nu 
\end{eqnarray}

\z where $c_m$ is defined in Lemma 39,  and $\|\Delta\|_\nu$ is
computed from (\ref{norm00}) with $K=\infty$.  Then
(\ref{normdistr00}) is a norm on $\mathcal{D}'_m$ and
$\mathcal{D}'_{m,\nu}=\left(\mathcal{D}'_m,+,*,\|\|_{m,\nu}\right)$ is
a Banach algebra.  With respect to the family of norms $\|\|_{m,\nu}$,
the projective limit of the $\mathcal{D}'_{m,\nu}$, $\mathcal{F}_{m}$
is a focusing algebra.

For any  $f\in L^1_{\nu_0}(\RR^+)$ there is a constant
$C(\nu,\nu_0)$ such that $f\in \mathcal{D}'_{m,\nu}$ for all
$\nu>\nu_0$ and

\begin{eqnarray}
  \label{majornorm}
  \|f\|_{\mathcal{D}'_{m,\nu}}\le C(\nu_0,\nu)\|f\|_{ L^1_{\nu_0}}
\end{eqnarray}

\z and formula (\ref{formulaforconv}) is equivalent to (\ref{defconv}),
in this case.

The operator $f(p)\mapsto pf(p)$ is continuous from $\mathcal{D}'_{m,\nu}$
to $\mathcal{D}'_{m,\nu+\delta}$ for any $\delta>0$. 

For $a\in\RR^+$,
$\mathcal{D}'_{m,\nu}(a,\infty)=\{f\in\mathcal{D}'_{m,\nu}:
\Delta_i(x)=0$ for  $x<a\}$ is a closed ideal in $\mathcal{D}'_{m,\nu}$
(isomorphic to the restriction $\mathcal{D}'_{m,\nu}(a,\infty)$ of
$\mathcal{D}'_{m,\nu}$ to $\mathcal{D}(a,\infty)$). The restrictions
$\mathcal{D}'_{m,\nu}(a,b)$ of $\mathcal{D}'_{m,\nu}$
to $\mathcal{D}(a,b)$ are for $0<a<b<\infty$ Banach spaces 
 with respect to the norm
(\ref{normdistr00}) restricted to $(a,b)$.

The functions in
$\mathcal{D}\left(\RR^+\backslash\NN\right)$ 
are dense in $\mathcal{D}'_{m,\nu}$, with respect to the norm
(\ref{normdistr00}) (Lemma~\ref{imbeddi}).

If we choose a different interval length $l>0$ instead of $l=1$ in the
partition associated to (\ref{stdec0}), we then write
$\mathcal{D}'_{m,\nu}(l)$. Obviously, dilation gives a natural
isomorphism between these structures. If
$d=\{t\mathrm{e}^{\mathrm{i}\phi}:t\in\RR^+\}$ is any ray,
$\mathcal{D}'_{m,\nu}(d)$ and
$\mathcal{F}_{m;\phi}$ are defined in an analogous way and have
the same properties as their real counterpart.

Laplace transforms are naturally defined in
$\mathcal{D}'_{m,\nu}$.

\begin{Lemma}\label{existe} Laplace transform extends continuously
from $\mathcal{D}(\RR^+\backslash\NN)$ to $\mathcal{D}'_{m,\nu}(\RR^+)$ by the formula

\begin{eqnarray}
  \label{laptra}
  (\lap f)(x):=\sum_{k=0}^{\infty}x^{mk}\int_0^{\infty}\mathrm{e}^{-sx}\Delta_k(s)\mathrm{d}s
\end{eqnarray}

\z In particular, with $f,g, h'\in\mathcal{D}'_{m,\nu }$ we have
\begin{eqnarray}
  \label{lapcomuta}
 && \lap(f*g)=\lap(f)\lap(g)\cr &&
 \lap(h')=x\lap(h)-h(0)\cr &&
\lap(pf)=-(\lap(f))'
\end{eqnarray}

For $x\in S_\nu= \{x:\Re(x)>\nu\} $ the sum (\ref{laptra}) converges
absolutely.  Laplace transform is,  for fixed $x\in S_\nu$, a continuous
functional (of norm less than one) on $\mathcal{D}'_{m,\nu}$.

$(\lap f)(x)$ is analytic in $S_\nu$.

Furthermore, $\mathcal{L}$  is {\em injective} in $\mathcal{D}'_{m,\nu}$. 
\end{Lemma}

\z The proof is given in \S~\ref{sec:LT}.
We conclude this section with a few   remarks.
\begin{Remark}\label{substiconv}\label{monotprop} Let $\mathcal{U}$ be one
 of the spaces considered in the examples and $\nu$ be large:

i) if  $g$ is analytic, (and if $g\in\mathcal{U}$ in Examples (2) and (4))
then $L_g:=f\mapsto f*g$ 
is a bounded operator and $\|L_g\|=O(\nu^{-1})$
($\mathcal{P}$ is such an operator, since $\mathcal{P}f=f*1$);

ii) replacing $*$ by $*_\phi$ defined as $f*_\phi g:=\mathrm{e}^{\mathrm{i}\phi}(f*g)$ leads
to an isomorphic structure;

iii) if $g\in\mathcal{U}$ is a real valued nonnegative function and
$|f|\le g$ ($|\mathcal{P}^{km}f|\le \mathcal{P}^{km}g$ for all $k$, on
$(0,k+1)$ if the space is $\mathcal{D}'_{m,\nu}$) then
$\|f\|\le\|g\|$.

\end{Remark}

(i) In $L^1_\nu$ and $\mathcal{D}'_{m,\nu}$ this follows from
the continuity of convolution and (\ref{majornorm}). In
the examples (3), the natural inclusion 
$\mathcal{T}_{\beta+\NN}\subset\mathcal{T}_{\beta}$ together
with (\ref{convodom1}) and (\ref{normconvo1})
makes convolution with an analytic function
 continuous in $\mathcal{T}_\beta$ and the claim follows from
the estimate
$\|f*g\|\le\max|g|\, \|f\|\, \|\mathcal{P}|\mathrm{e}^{\nu p}|\|$.

(ii) The isomorphism is given by $f\mapsto \mathrm{e}^{-\mathrm{i}\phi}f$.

(iii) Since $\mathcal{P}$ is positivity preserving, writing $|f|\le g
$ as $ -g\le (\Re,\Im) f\le g$ the property is obvious when $f,g$ are
functions, while for $\mathcal{D}'_{m,\nu}$ it follows from equation
(\ref{defDelI}) below.

\Box

\subsubsection{Vectorial convolution and focusing spaces}

We endow $\mathcal{B}_{\nu}^n$ with a Banach space structure by
identifying it with $\mathcal{B}_\nu\oplus
\cdots\oplus\mathcal{B}_\nu$ ($n$ times). The projective limit of the
$\mathcal{B}_\nu$, $\mathcal{F}^n$ is, clearly, a focusing space.  We
define a convolution on
$(\cdot\,*\,\cdot):\mathcal{F}^n\mapsto\mathcal{F}$ ({\em not}
$\mathcal{F}^n\mapsto \mathcal{F}^n$) by

\begin{eqnarray}
  \label{convect2}
  \bfV*\mathbf{W}:=\sum_{i,j=1}^n V_i*W_j 
\end{eqnarray}

\z
We write $\bfV^{*\bfl}:=V_1^{*l_1}*V_2^{*l_2}*\cdots*V_n^{*l_n}$
with the conventions $V^{*1}=V$ and that the factors 
with $l_i=0$ are omitted.

\subsubsection{A fixed point property}

\begin{Lemma}\label{gennonsense} Let $\mathcal{F}$ be a focusing space and
$\mathcal{N}$ be a (linear or nonlinear) operator defined on
$\mathcal{F}$.
Equivalently, in view of (\ref{focus-pocus}), let $\mathcal{N}$
be defined on 
$\bigcup_{\nu>\nu_0}B_\nu(\delta)$ with
$B_\nu(\delta)=\{f:\|f\|_\nu\le\delta\}$ for some $\delta>0$.
Assume $\mathcal{N}$ is
{\bf eventually contractive} in the following sense. There exist
$\nu_0,\,\epsilon>0$, $\alpha<1$, so that if $\nu\ge\nu_0$ and
$\|f\|_\nu + \|g\|_\nu\le \epsilon$ then

\begin{eqnarray}
  \label{condop}
  \|\mathcal{N}(f+g)
-\mathcal{N}(g)\|_\nu\le \alpha\|g\|_\nu
\end{eqnarray}

\z Then $\mathcal{N}$  has a unique
fixed point $f_0\in\mathcal{F}$.

If $\mathcal{N}$ depends continuously (in the strong topology)
on a parameter
$\phi$ for  $\nu>\nu_0$  and if
the constants $\nu_0,\alpha$  and $\epsilon$ above
do not depend on $\phi$, then 
the fixed point $f_\phi$ is also continuous in $\phi$. Furthermore,
$\lim_{\nu\rightarrow\infty}\sup_{\phi}\|f_\phi\|_{\nu}=0$.

\end{Lemma}

The proof is straightforward. To show existence, take $\nu>\nu_0$
large enough so that, by (\ref{focus-pocus})
$\|\mathcal{N}(0)\|_\nu<(1-\alpha)\epsilon$. Then the closed ball
$B_\nu(\epsilon)$ is mapped by $\mathcal{N}$ into itself
for any $\phi$ by (\ref{condop}) and $\mathcal{N}$ is contractive
there.  The fixed point obtained, for instance, as the limit of the
(convergent, uniformly in $\phi$) iteration
$\phi_{n+1}=\mathcal{N}(\phi_n); \phi_0=0$ is continuous in $\phi$
since $\mathcal{N}$ is. By construction 
$\|f_0\|_\nu\le \epsilon$, for all $\phi$.

For uniqueness, let $f_0$ and $f_1$ be fixed points of
$\mathcal{N}$;  by (\ref{focus-pocus}) there is a $\nu>\nu_0$ so that
$\|f_{0,1}\|_\nu<\epsilon$. Then by (\ref{condop}),
$\|f_0-f_1\|_\nu\le\alpha\|f_0-f_1\|_\nu$ and thus $f_0=f_1$.  \Box.

Let $\nu_0>0$ and let $\{\mathbf{M}\}_{\bfl}$ be a sequence of linear
operators $\mathbf{M}_\bfl:\mathcal{B}_\nu\mapsto\mathcal{B}_\nu^n$
for $\nu\ge\nu_0$.
Assume that for some $\kappa$ and all $\bfl$, $|\bfl|\ge 1$,

\begin{eqnarray}
  \label{condgG}
   &&\|\mathbf{M}_{\bfl}\|_\nu\le C_\nu\kappa^{|\bfl|}\ \mbox{ and }
 C_{1,\nu}:=\lim_{\nu\rightarrow\infty}\max_{|\bfl|=1}\|\mathbf{M}_{\bfl}\|_\nu=0
\end{eqnarray}

\z Let ${\bf F}_0\in\mathcal{F}^n$ and
$\mathbf{M}:\mathcal{F}^n\mapsto\mathcal{F}^n$ be defined by

\begin{eqnarray}\label{Mdef}
\mathbf{M}(\mathbf{Y}):=\mathbf{F}_0+
\sum_{|\bf l|\ge 1}\mathbf{M}_{\bf l}\left({\bf Y}^{*\bf l}\right)
\end{eqnarray}

\begin{Lemma}\label{infico}

 $\mathbf{M}$ satisfies the assumptions of  Lemma~\ref{gennonsense}.
$\mathbf{M}$   has therefore a unique
fixed point in $\mathcal{F}^n$.

\end{Lemma}

\z {\em Proof}. We first need the following estimate.

\begin{Remark}\label{4}  Let $\mathbf{V},\mathbf{W}\in\mathcal{F}^n$. For $|\bfl|>0$ 
and any $\nu$ we have,
with $\|\|=\|\|_\nu$,
 
\begin{eqnarray}\label{estconvn}
\|\mathbf{W}_\bfl\|:=\|({\bf V}+{\bf W})^{*\bfl}-\mathbf{V}^{*\bfl}\|\le |\bfl|\left(
\|\mathbf{V}\|+\|\mathbf{W}\|\right)^{|\bfl|-1}\|\mathbf{W}\|
\end{eqnarray}

\end{Remark}

This inequality is obtained by induction on $\bfl$, with respect to
$\prec$.  For $|\bfl|=1$, (\ref{estconvn}) is trivial. Assume
(\ref{estconvn}) holds for all $\bfl\prec\bfl_1$; without loss of
generality we may consider that $\bfl_1 =\bfl_0+\mathbf{e}_1$. We have:

\begin{multline*}
  \|({\bf V}+{\bf W})^{*\bfl_1}-\mathbf{V}^{*\bfl_1}\|= \|({\bf
    V}+{\bf W})^{*\bfl_0}*({\bf V}_1+{\bf
    W}_1)-\mathbf{V}^{*\bfl_1}\|\cr
  =\|(\mathbf{V}^{*\bfl_0}+\!
\mathbf{W}_{\bfl_0})*(V_1+W_1)-\!\mathbf{V}^{*\bfl_1}\|=
  \|\mathbf{V}^{*\bfl_0}\!*\!W_1+\mathbf{W}_{\bfl_0}\!*\!V_1+
\mathbf{W}_{\bfl_0}\!*\!W_1\|\cr
  \le
  \|\mathbf{V}\|^{|\bfl_0|}\|\mathbf{W}\|+
\|\mathbf{W}_{\bfl_0}\|\|\mathbf{V}\|+
  \|\mathbf{W}_{\bfl_0}\|\|\mathbf{W}\|\cr\le \|\mathbf{W}\|
  \left(\|\mathbf{V}\|^{|\bfl_0|}+
|\bfl_0|(\|\mathbf{V}\|+\|\mathbf{W}\|)^{|\bfl_0|}\right)\le
    \|\mathbf{W}\|(|\bfl_0|+1)(\|\mathbf{V}\|+\|\mathbf{W}\|)^{|\bfl_0|}
\end{multline*}

\z and (\ref{estconvn}) is proven.

For the sum in $\mathbf{M}$ to converge in $\|\|_\nu $ it suffices to
choose $\nu $ such that $\|\bfV\|_\nu <\kappa^{-1}$. Let
$\epsilon<\kappa^{-1}$ and $\bfV,\mathbf{W}$ be such that
$\|\bfV\|+\|\mathbf{W}\|<\epsilon$.  We have

\begin{multline}
  \label{difN}
 \left\|\mathbf{M}(\bfV+\mathbf{W})-\mathbf{M}(\bfV)\right\|_\nu
\le \left(nC_{1,\nu}+C_\nu\sum_{|\bf l|\ge 2}
|\bfl|(\kappa\epsilon)^{|\bfl|-1}\right)\|\mathbf{W}\|_\nu\cr
\le \left(nC_{1,\nu}+n2^nC_\nu\kappa\epsilon+nC_\nu\frac{(2-\kappa\epsilon)(\kappa\epsilon)^{1+2(n-1)}}
{(1-\kappa\epsilon)^{n+1}}\right)\|\mathbf{W}\|_\nu=K_\nu\|\mathbf{W}\|_\nu
\end{multline}

\z where we separated out the terms with $|\bfl|=2, l_i=0$ or $1$ and
for the rest of the terms used the identity $\sum_{l_i\ge
  2}|\bfl|\mathrm{e}^{-\gamma|\bfl|}=- \frac{\mathrm{d}}{\mathrm{d}\gamma}\left(\sum_{l\ge
  2}\mathrm{e}^{-\gamma l}\right)^n$.  We see that, in fact,
$\lim_{\nu\rightarrow\infty}K_\nu=0$. \Box

\subsubsection{Two lemmas on analytic structure}\label{sec:twolemmas}

\begin{Lemma}\label{cuteCauchy} Let $f$ be analytic in the unit disc cut along the
positive axis and let $0<g(x)\in C^1[0,1]$. Assume that 
$\lim_{\epsilon\downarrow
  0}f(x\pm \mathrm{i}\epsilon g(x)) = f^{\pm}(x)$ in $\lone[0,1]$ and 

\begin{eqnarray}
  \label{condipo}
 f^{+}(x )-f^{-}(x)=f_\delta(x)=x^rA(x )
\end{eqnarray}

\z with $\Re(r)>-1$, where $A(\xi)$ extends to an analytic function for
$|\xi |<a\le 1$. Then there exists a function
$B$ analytic in $|\xi|<a$ so that

\begin{eqnarray}
  \label{concluelempo}
  &&f(\xi )=\frac{1}{1-\mathrm{e}^{2\pi i r}}\xi ^rA(\xi )+B(\xi )\ \ \
  (r\notin\NN)\cr&&
 f(\xi )=\frac{i}{2\pi}\ln(\xi )\xi^r A(\xi )+B(\xi ) \ \ \
  (r\in\NN)\cr&&
\end{eqnarray}

\z If $ f^{+}(x )-f^{-}(x)$ is a linear combination
$\sum_{i=1}^N x^{r_i}A_i(x )$ (under the same assumptions
on $r_i$ and $A_i$), then $f$ is given by 
the corresponding superposition of terms of the 
form (\ref{concluelempo}).

\end{Lemma}

The proof is given in \S~\ref{sec:A0}.

In the following, $\gamma:\RR^+\mapsto\CC$ will denote smooth curves
in $\mathcal{R}'_1$, $\gamma_k$ denotes a curve that crosses through
the interval $(k,k+1)$, $\gamma_{\epsilon}=\Re(\gamma)+i
\epsilon\Im(\gamma)$ (cf.\S~\ref{sec:anset}).
Let $a,b\in(0,\pi/2)$ and
$\mathcal{S}_0=\{p:\arg(p)\in (\psi_-,0)\cup(0,\psi_+)\}$.
Let $f$ be a function analytic in $\mathcal{R}'_1$ so that
$f\circ\gamma_\epsilon\in\mathcal{D}'_{m,\nu}$ has limits in $
\mathcal{D}'_{m,\nu}$ as $\epsilon\downarrow 0$. We denote
the space of such functions by $\mathcal{D}'_{m,\nu}(\mathcal{R}_1)$.
Let $ f ^{\pm}_0=f^{\pm}$, $F_{j}=\mathcal{P}^{(mj)}f$, and for $j>0$
\begin{gather}
  \label{decom0}
 f_{j}^+(z-j)=F_j^{-j+}(z)-F_j^{-(j-1)+}(z);\  f_j^-(z-j)=F_j^{+j-}(z)-F_j^{+(j-1)-}(z)
\end{gather}

\z  By construction the $ f_j^+$
are in $\lone[0,1-\epsilon)$, analytic in a sectorial neighborhood of
$z=0$ and can be extended analytically for $\Im(z)>0,\Re(z)>0$ (this
last property motivates the choice $+$ for superscript, while the
right shift is chosen in view of our application).  Also by
construction $ f_j^{\pm}(z)=0$ for $z<0$; it is convenient to extend
$ f_j^{\pm}$ by zero throughout $\Re(z)<0$.  We have the
``telescopic'' decomposition

\begin{eqnarray}
  \label{defindecom2}
 f^{\mp\,j\,\pm}(z)=\sum_{i=0}^{j}( f_i^{\pm})^{(mi)}(z-i)
\end{eqnarray}

Relation (\ref{defindecom2}) holds in
$\mathcal{D}'$ along the real axis, {\em and} as an equality of analytic
functions for $\Re(z)>j,\Im(z)>0$. 
 For instance, for 
$f=(z-2)^{-1}\in\mathcal{D}'_{1\nu}$ we have $ f_1^{+}=0$,
$ f_2^{+}=-2\pi i z$ for $\Re(z)>0$, and $f_2^{+}=0$ otherwise.

Conversely, a decomposition of the form (\ref{defindecom2}) together
with analyticity in $\mathcal{S}_0$ implies analyticity in
$\mathcal{R}'_1$. More precisely, assume that
$f(t\exp(\mathrm{i}\phi))\in\mathcal{D}'_{m,\nu}(\RR^+)$ for $\phi\in
(\psi_-,0)\cup(0,\psi_+)$ and that $\lim_{\phi\downarrow
  0}f(\cdot\,\exp(\pm \mathrm{i}\phi)=f^{\pm}$ in $\mathcal{D}'_{m,\nu}$. Assume in
addition that there exists the decomposition

\begin{eqnarray}
  \label{decom}
  f^{\pm}=\sum_{k=0}^{\infty}
\left[ f ^{\mp}_k(p-k)\bchi_{[k,\infty]}\right]^{(mk)}
\end{eqnarray}

\z  where for each
$k$, $ f ^{\mp}_k\in\mathcal{D}'_{m,\nu}(\RR^+)$ (note: the $ f ^{\mp}_k\bchi$
are uniquely determined, cf. Remark~\ref{density}). Assume in addition
that $ f ^{\mp}_k$
extend analytically to $\mathcal{S}_0^{\mp}$ in the following sense:
there exist $g_k^{\mp}$ analytic in $\mathcal{S}_0^{\mp}$, with
$g_k^{\mp}(\mp t\exp(\mathrm{i}\phi))\in \mathcal{D}'_{m,\nu}(\RR^+)$ and such
that $\lim_{\phi\downarrow 0}g_k^{\mp}(\mp t\exp(\mathrm{i}\phi))= f_k^{\mp}$ in
$\mathcal{D}'_{m,\nu}(\RR^+)$. 
\begin{Lemma}\label{STR1} (i) Under the above conditions, $f$ extends analytically
to $\mathcal{R}'_1$.

(ii) If for small argument $ f_k\in\mathcal{T}_\beta$ then 

\begin{eqnarray}
  \label{ASTR1}
  f^{\pm}(p+k)= \frac{1}{1-\mathrm{e}^{2\pi \mathrm{i}\beta}}(f^{\mp}_k(p))^{(mk)}+a(p)
\mbox{ or }\frac{i}{2\pi}( f ^{\mp}_k(p)\ln p)^{(mk)}+a(p)
\end{eqnarray}

\z according to whether $\beta\notin\ZZ$ or $\beta\in\ZZ$
respectively, where $a$ is analytic at zero. As in Lemma
~\ref{cuteCauchy}, if $ f_k=\sum_{i\le j} f_{ki}$
with $ f_{ki}\in\mathcal{T}_{\beta_i}$ then
$f^{\pm}(p+k)$ is the corresponding superposition 
of terms of the form (\ref{ASTR1}).

\end{Lemma}

The proof is given in \S~\ref{sec:A0}.

\subsubsection{Convolutions, analyticity and averaging}
\label{sec:aver}

Define $\mathcal{F}(\mathcal{R}'_1)$
as the functions in $\mathcal{D}'_{m,\nu}(\mathcal{R}'_1)$ 
such that,
in the decomposition (\ref{defindecom2}) $\|( f_{j}(p-j))^{(mj)}\|\le
K_f(\nu)^j$ where $\lim_{\nu\rightarrow\infty}K_f(\nu)=0$.

If $\gamma$ is a straight line in $\mathcal{R}'_1$
then  $AC_\gamma(f*g)=AC_\gamma(
f)_\gamma*AC_\gamma(g)$ (if $\gamma$
is {\em not} equivalent to a straight line, this equality
is generally {\em false} \cite{Costin}).
Convolution does however commute with suitable averages of analytic
continuations, as Proposition~\ref{medianpropo} below shows.  In view of symmetry, we only need to look at the
properties of the $+$ decomposition.

For $\alpha\in\CC$,  consider the operator
$\mathcal{A}_\alpha:\mathcal{F}(\mathcal{R}'_1)\mapsto \mathcal{F}_m(\RR^+)$
given by

\begin{eqnarray}
  \label{generalmed}
  \mathcal{A}_\alpha(f):=f^{[\alpha]}(p)=\sum_{i=0}^{\infty}\alpha^i\left( f ^+_i(p-j)\right)^{(mi)}
\end{eqnarray}

\z In our assumptions, convergence is ensured in
$\mathcal{D}'_{m,\nu}$ for large enough $\nu$. An important case is
the balanced average, $\alpha=1/2$. The operator $\mathcal{A}_{\frac{1}{2}}$ is
similar to \'Ecalle's medianization, and is designed to substitute for
analytic continuation along the singularity line $\RR^+$ in a way
compatible with the $*-$algebra structure. As mentioned before, it
can be shown that under our assumptions on (\ref{eqor}), only for the
choice $\alpha=1/2$ is the difference between the
$f=\lap\mathcal{A}_\alpha F$ and the optimally truncated asymptotic
series of $f$ always of the order of magnitude of the least term of
the series \cite{CK2}.  Borel summability techniques and
hyperasymptotic methods \cite{Berry}, \cite{Berry-hyp} give, whenever
they both apply, the same association between transseries and actual
functions.

\begin{Proposition}\label{medianpropo} i) If $f,g\in\mathcal{F}(\mathcal{R}'_1)$ then $f*g$ defined
for small argument by (\ref{defconv}) extends
to a function in $\mathcal{F}(\mathcal{R}'_1)$. We have

\begin{eqnarray}
  \label{defindecom3}
(f*g)_j^{\pm}=\sum_{s=0}^j f_s^{\pm}*g_{j-s}^{\pm}
\end{eqnarray}

\z and $K_{f*g}(\nu)\le K_f(\nu)+K_g(\nu)$. If $h$ is analytic and
bounded in the right half plane and $f\in\mathcal{F}(\mathcal{R}'_1)$
then $hf\in\mathcal{F}(\mathcal{R}'_1)$.

ii) If $h$ is analytic in $\mathcal{R}'_1\cup\RR^+$ and
$f,g\in\mathcal{F}(\mathcal{R}'_1)$, $a,b\in\CC$,  then

\begin{eqnarray}
  \label{assertmed}
  \mathcal{A}_\alpha(af+bg)&=&a\mathcal{A}_\alpha(f)+b
  \mathcal{A}_\alpha(g)
\cr \mathcal{A}_\alpha(hf)&=&h\mathcal{A}_\alpha(f)\cr
  \mathcal{A}_\alpha(1)&=&1\cr
  \mathcal{A}_\alpha(f*g)&=&\mathcal{A}_\alpha(f)*\mathcal{A}_\alpha(g)
 \end{eqnarray}

\z If $\mathbf{M}$ satisfies the hypothesis of Lemma~\ref{gennonsense} and in
addition
$\mathbf{M}_{\bfl}(\mathcal{F}(\mathcal{R}'_1))\subset\mathcal{F}(\mathcal{R}'_1)$ and $\mathcal{A}_\alpha
\mathbf{M}_{\bfl}=\mathbf{M}_{\bfl}\mathcal{A}_\alpha$ 
then 

\begin{eqnarray}
  \label{assertmed3}
  \mathcal{A}_\alpha \mathbf{M}=\mathbf{M}\mathcal{A}_\alpha 
\end{eqnarray}

\z In particular, if $\bfY$ is a fixed point of $\mathcal{M}$ then so
is $\mathcal{A}_\alpha\bfY$. An example is the case
$\mathbf{M}_{\bfl}(\bfY)= \mathbf{G}_{\bfl}*\bfY^{*\bfl}$ with
$\mathbf{G}_{\bfl}$ analytic in $\mathcal{R}'_1\cup\RR^+$ and such that
for some $\kappa$ and all $\bfl$ we have $|\mathbf{G}_{\bfl}(p)|\le\exp(\kappa
p)$.

\end{Proposition}

\z The proof and further details are given in
Appendix~{\ref{sec:A0}}. 

Let now
$\mathcal{F}_r(\mathcal{R}'_1)\subset\mathcal{F}(\mathcal{R}'_1)$
consist in functions $f$ whose only singularities are regular, in the
sense that the elements $f_{j}$ (cf. (\ref{decom0})) are of the form
$(\sum_{i=1}^{N_j}p^{a_{i,j}}A_{i,j})^{(mj)}$ where $A_{i,j}$ are
analytic near $p=0$.

\begin{Remark}\label{preservestruct} $\mathcal{F}_r(\mathcal{R}'_1)$ is stable with respect to
convolution.

\end{Remark}

\z By construction, (Proposition~\ref{medianpropo})
and (\ref{defindecom3}), for small $p$,
$(f_{1}*f_2)_j $ is a sum of terms of the form

$$p^{a_1}A_1*p^{a_2}A_2$$

\z and the proof follows without any difficulty from
(\ref{defconv31}).

\subsection{Main proofs}

\begin{Proposition}\label{9}\label{estimd} i)
For any $\kappa>\max\{x_0,y_0^{-1}\}$, cf. (\ref{Taylor
  series}), there is a constant $K>0$ such that for all $\bfl\succ 0$

\begin{eqnarray}\label{exponesti}
\sup_{p\in\CC}\mathrm{e}^{-\kappa |p|}|{\bf G}_{\bf l}(p)|\le K \kappa^{|\bf l|} 
\end{eqnarray}

\z  (cf. (\ref{eqil})).  

ii) Let $\mathcal{F}$ be one of the focusing algebras in
\S~\ref{sec:NC} and $\bfY\in\mathcal{F}$.  Let

$$ \bfd_\bfj= \sum_{ \bf l\ge j} \binom{\bfl}{\bfj
  }\left[ \bfg_\bfl*\bfY^{*(\bfl -\bfj) }+\bfgg_{0,\bfl}*\bfY^{*(\bfl
  -\bfj)}\right]$$

\z Then for large $\nu$ and some $\kappa_1>0$, $\|\mathbf{D}_\bfj\|\le
\kappa_1^{|\bfj|}$ while for $|\mathbf{j}|=1$,
$\|\mathbf{D}_\bfj\|=
o(\nu^{-M})$.
\end{Proposition}

{\em Proof.}

\z  (i) From the analyticity of $\mathbf{g}$ it follows 
 that

\begin{eqnarray}\label{unifesti}
|{\bf g}_{m,\bf l}|<\mbox {const. } \kappa^{m+|\bfl|}
\end{eqnarray}

\z where the constant is independent on $m$ and $\bfl$.  
Then, by (\ref{lapdef}),

\begin{equation}\label{estimmodg}
|{\bf G}_{\bf l}(p)|<\mathrm{const.}\ \kappa^{|\bfl|+1}\frac{\mathrm{e}^{\kappa| p|}-1}{\kappa |p|}<
\mathrm{const.}\ \kappa^{|\bfl|+1}\mathrm{e}^{\kappa |p|}
\end{equation}

\z The last claim is a direct consequence of (n5).

\bigskip

(ii)
Note first that $\sum_{|\bfl|=l}1\le l^n\le 2^{nl}$ (since
$l_i\le l$). Also,  $\binom{l_i}{j_i}\le 2^{l_i}$  so that
$\binom{\bfl}{\bfj}\le 2^{n|\bfl|}$. By (n5) we have  $\mathbf{g}_{0,\bfl}=0$ if
$|\bfl|\le 1$. Choosing $\delta<2^{2nl+1}\kappa$  and $\nu$ such that
$\|\bfY\|\le\delta$,  $\|\mathbf{D}_\bfj\|$ is  estimated
by

$$\sum_{l\ge j}2^{2nl}\delta^{l-j}\kappa_l^l$$

\z where $j=|\bfj|,\kappa_l=\|\mathbf{G}_\bfl\|+|\mathbf{g}_{0,\bfl}|\le 2\kappa$ if $l>1$ and $\kappa_l=o(\nu^{-M})$
for $l=1$, by (i), and the result follows. 

\Box

Without loss of generality, we analyze (\ref{eqil}) in
a neighborhood of $d_1$, the Stokes line corresponding
to $\lambda_1=1$. (For the equations (\ref{eqMv})
we will need, in addition, to study a direction where
$p-\lambda_i+\mathbf{k}\cdot \boldsymbol\lambda=0$.)

\z  Let $\epsilon$ 
and $c_0$ be small  and positive, $\mathcal{V}=\{p:|p|<1\}$,

\begin{equation}
  \label{defSC}
\mathcal{S}_c=\{p:\arg(p)\in [\psi_n-2\pi+c,-c]\cup[c,\psi_2-c]\}
\end{equation}

\z $\mathcal{S}_0=\cup_{0<c<c_0}\mathcal{S}_c$, $\mathcal{S}_c^{\pm}=
\mathcal{S}_c\cap\{p:\pm\arg(p)>0\}$, and let
$\mathcal{S}'_0,\mathcal{S}'_c,{\mathcal{S}_c^{\pm}}'$ be defined
correspondingly, with $\psi_-$ and $\psi_+$ replacing $\psi_n-2\pi$
and $\psi_2$, respectively.

\begin{Lemma}\label{analyticase} i) For any ray $d\subset \mathcal{S}_0$, there is a
  unique solution of (\ref{eqil}) in $\mathcal{D}'_m({d},l)$ for any
  $m\in\NN$, $l\in(0,1)$.  This solution, $\bfY_0$, is in fact
  analytic in $\mathcal{S}_0\cup\mathcal{V}$ and, for large enough
  $\nu(d)$, $\bfY_0\in L^1_\nu(d)$.

ii) The function $\bfY_0(\cdot\,\mathrm{e}^{\mathrm{i}\phi})$ is continuous in
$\phi\in(\psi_n,0]$ and (separately in)
$\phi\in[0,\psi_2)$ in the
$\mathcal{D}'_{\nu,m}(\RR^+,l)$ topology and
$\sup_{\phi\in[\psi_1,\psi_2]}\|\bfY_0(p\mathrm{e}^{\mathrm{i}\phi})\|_{m,\nu}
\rightarrow 0$ as $\nu\rightarrow\infty$.

iii) The description (\ref{SY0}) holds for $j=1,l=1$.

iv) For $a>1$, there is a one-parameter family
of solutions of (\ref{eqil}) in $\mathcal{D}'_m(0,a)$.

\end{Lemma}

\z {\em Note:} the hyperfunctions
$\bfY_0(p\mathrm{e}^{\pm 0i})=\bfY_0^{\pm}$ are different, in general.

\begin{Corollary}\label{prim} For $k\in\NN\cup\{0\}$, the function
$\mathcal{P}^{mk}\bfY_0(p\mathrm{e}^{\mathrm{i}\phi})$ is  continuous in
$S_k^{-}:=\{p:0\le|p|<k+1;\arg(p)\in (\psi_n,0]\}$ and in
$S_k^{+}:=\{p:0\le|p|<k+1;\arg(p)\in [0,\psi_2)\}$ (and, of course,
analytic in  $S_0$).  \end{Corollary}

{\em Proof of Lemma~\ref{analyticase}.}

We write  (\ref{eqil}) in the form:

\begin{eqnarray}\label{eqilm}
\bfY=
\left(\hat{\Lambda}-p\right)^{-1}\left({\bf F}_0-
\hat B\mathcal{P}\bfY+{\cal  N}({\bf Y})\right)=\mathcal{M}(\bfY)
\end{eqnarray}

\z Let $d_K$ be an initial segment of $d$ of length $K<\infty$.  As the matrix
$\hat{\Lambda}-p$ is invertible in $\mathcal{S}_c$, it is easy to see
that the operator $\mathcal{M}$ in (\ref{eqilm}) satisfies the
conditions of Lemma~\ref{infico}, in the spaces $ L^1_\nu(d)$,
$\mathcal{A}_{z,l}(\mathcal{S})$ ($l\le M$) (cf.  Examples (1) through
(4) in \S\ref{sec:NC} and Remark~\ref{9}).  The conditions are in
addition also satisfied in $\mathcal{D}'_{m,\nu}(d)$; due to the
special structure of this space, the proof the boundedness of the
operator $U=(\Lambda-p)^{-1}$ is more delicate, and is given  in
Lemma~\ref{cinftycase}.

Thus, $\mathcal{M}$ has a unique
fixed point in each of these spaces. The obvious inclusions between
these spaces complete the proof of part (i). For the rest of 
Lemma~\ref{analyticase} we need more results.

\begin{Proposition}\label{IN0}
The properties stated in  Lemma~\ref{analyticase} 
hold in $\mathcal{S}_0\cap\{p:|p|<1+\epsilon\}$

\end{Proposition}

The proof  is given in \S~\ref{sec:A1}. \Box

\begin{Proposition}\label{contidis}
Let $\bfW_0$ be a solution of (\ref{eqil}) in $\mathcal{D}'_m(0,a)$
with $a>1$. For $b\ge a$ there exists a unique solution of (\ref{eqil})
in $\mathcal{D}'_m(0,b)$, which
agrees
with $\bfW_0$ on $\mathcal{D}(0,a)$.

\end{Proposition}

\z We use the decomposition $\mathcal{D}'_m(0,b)=
\mathcal{D}'_m(0,a)\bigoplus\mathcal{D}'_m(a,b)$, $(a<b\le\infty)$ We
identify $\bfW_0$ with an element of $\mathcal{D}'_m(0,b)$ by
extending it with zero and define $\mathcal{M}_1$ on
$\mathcal{F}_m(a,b)=\cup_{\nu>\nu_0}\mathcal{D}'_{m,\nu}(a,b)$ by

\begin{eqnarray}
  \label{decom2}
  \mathcal{M}_1(\bfW_1)=\mathcal{M}(\bfW_0+\bfW_1)-\mathcal{M}(\bfW_0)
\end{eqnarray}

\z and  (\ref{eqilm}) becomes

\begin{eqnarray}
  \label{eqilm1}
  \bfW_1=\mathcal{M}_1(\bfW_1)
\end{eqnarray}

\z By Lemma~\ref{gennonsense}, $\mathcal{M}$ is eventually contractive
and by (\ref{decom2}), clearly, so is $\mathcal{M}_1$. By
Lemma~\ref{gennonsense} then, $\mathcal{M}_1$ has a unique fixed point
in $\mathcal{F}_{m}(a,b)$ ($b\le\infty$).  In view of
the inclusions between $\mathcal{F}_{m}(a,\infty)$ and
$\mathcal{F}_{m}(a,b)$, the proof is complete. \Box.

{\em Proof of Lemma~\ref{analyticase} (ii)}

Let $\phi\in[0,\epsilon]$ with $\epsilon$ small. We regard the space
$\mathcal{F}_{m}$ as fixed and $\mathcal{M}_1$ as being dependent on
$\phi$ through $p$ and $*_\phi$ (cf. Remark~\ref{substiconv}).
Firstly, $\bfW_0(\phi)$ is continuous in $\phi\in[0,\epsilon]$ and
$\|\bfW_0(p\mathrm{e}^{\mathrm{i}\phi})\|_\nu=O(1/\nu)$ uniformly in $\phi$ as follows
from Proposition~\ref{IN0}.  Then the infinite sum in the definition
of $\mathcal{M}_1$ is uniformly convergent in
$\mathcal{D}'_{m,\nu}(a,\infty)$ for $\nu$ large enough by
(\ref{exponesti}).  The operator
$U=(p\mathrm{e}^{\mathrm{i}\phi}-\hat{\Lambda})^{-1}=\mathrm{e}^{-\mathrm{i}\phi}(p-\mathrm{e}^{-\mathrm{i}\phi}\hat{\Lambda})^{-1}$
is strongly continuous in $\phi$ in $\mathcal{F}_{m}(a,\infty)$,
$a>1$, cf. Lemma~\ref{cinftycase}. By Remark~\ref{substiconv},
$\mathcal{M}_1$ is $\phi$-continuous and Lemma~\ref{gennonsense}
applies.

 \Box

\z We now study the convolution equations associated to the
higher terms in the transseries, (\ref{eqMv}).

\begin{Lemma}\label{higherterms} i) Given the vector of constants $\bfC\in\CC^{n_1}$
  and  and in addition
  given $\bfY_0$ for $d=\RR^+$, there is a unique solution of
  (\ref{invlapvk}) in $\mathcal{F}_{m}$ with the (singular) initial
  condition

\begin{eqnarray}
  \label{condinit}
  \bfY_{\mathbf{e}_j}(p)= C_j
  \Gamma(\beta_j')^{-1}p^{\beta_j'}(\mathbf{e}_j+o(1)) \qquad(p\rightarrow
  0,\ 
j=1,2,\ldots,n_1)
\end{eqnarray} 

\z The general solution 
of (\ref{eqsvec0}), (\ref{invlapvk}) is 

\begin{eqnarray}
  \label{gesovk}
\mathbf{C}^\bfk\mathbf{Y}_\bfk, \ \mathbf{C}\in\CC^{n_1}
\end{eqnarray}

\z where $\mathbf{Y}_\bfk$ is the solution for
$\mathbf{C}=(1,1,\ldots,1)$.

ii) In a neighborhood
of $p=0$ we have

$$\bfY_\bfk(p)=p^{\bfk\bfbet'-1}\bfA_\bfk(p)$$

\z with $\bfA_\bfk$ analytic near the origin.

iii) The functions $\bfY_\bfk$, $\bfk\succeq 0$ are analytic in
$\mathcal{R}'_1$ and 
$\bfY_\bfk(x\mathrm{e}^{\mathrm{i}\phi})$ are continuous in $\phi$ with respect to the
$\mathcal{D}'_{m,\nu}$ topology for $\phi\in(\psi_-,0]$
and for $\phi\in[0,\psi_+)$.

iv) Each $\bfY_\bfk$ is in $\mathcal{F}(\mathcal{R}'_1)$
(cf. \S~\ref{sec:aver}).
  Furthermore, there is a constant
$K$ and a function $\delta(\nu)$ such that
$\lim_{\nu\rightarrow\infty}\delta(\nu)=0$ and in the decomposition
(\ref{defindecom2}) of $\mathbf{Y}_\bfk$ we have
$\mathbf{Y}_{\bfk;j}\in\mathcal{T}_{(\bfk+j\mathbf{e}_1)\bfbet'-1}$ and

\begin{eqnarray}
  \label{normVK}
&&\sup_{\phi,\bfk}\delta(\nu)^{-|\bfk|}\|\mathbf{Y}_{\bfk}\|_{\mathcal{D}'_{m,\nu}(\mathrm{e}^{\mathrm{i}\phi}\RR^+)}<K\cr&&
  \sup_{\phi,\bfk,j}\delta(\nu)^{-|\bfk|+j}\|\mathbf{Y}_{\bfk;j}\|_{\mathcal{D}'_{m,\nu}(\mathrm{e}^{\mathrm{i}\phi}\RR^+)}<K
\end{eqnarray}

\z where $\phi$ runs in $(\psi_-,\psi_+)$. 

v) The functions $\mathbf{Y}_\bfk(\cdot\,\mathrm{e}^{\mathrm{i}\phi})$,
$\phi\in(\psi_-,\psi_+)$, are simultaneously Laplace transformable in
$\mathcal{D}'_{m,\nu}(\RR^+)$ and their Laplace transforms are
solutions of (\ref{systemformv}).  The expression

\begin{eqnarray}
  \label{solupperlower}
 \bfy^{\pm} =\lap\bfY_0^{\pm}+\sum_{\bfk\succ 0}x^{\bfm\cdot\bfk}\bfC^{\bfk}\mathrm{e}^{-\bfk\cdot\bflam x}
  \lap \bfY^{\pm}_\bfk
\end{eqnarray}

\z 
is uniformly convergent for large enough $x(\mathbf{C})$ in some open
sector. In addition, (\ref{solupperlower}) is a solution
of $(\ref{eqor})$ and $\lap\bfY^{\pm}_\bfk\sim\tilde{\bfy}_\bfk$
for large $x$ in the half plane $\Re(xp)>0$. 
 
\end{Lemma}

Without loss of  generality, we analyze (\ref{invlapvk}) in
$\mathcal{T}_{\{\cdot\}}({\mathcal{S}_c^+}')$, $\mathcal{D}'_{m,\nu}(d)$,
with $d\in{\mathcal{S}_c^+}'$, and in $\mathcal{D}'_{m,\nu}(\RR^+)$.  We
denote all the corresponding norms by $\|\|_\nu$.

\begin{Remark}\label{stableconv} Assume that for $\bfk'\prec\bfk$ we have
$\bfY_{\bfk'}\in\mathcal{T}_{\bfk'\bfbet'-1}$. Then, in (\ref{defT}),
we have
$\bfT_\bfk(\bfY_0,\{\bfY_{\bfk'}\})\in\mathcal{T}_{\bfk\bfbet'-1}$.
\end{Remark}

This follows immediately from Equations (\ref{defconv31}) and
(\ref{convodom1}) and from the homogeneity of $\mathcal{T}$ implicit
in the sum $\sum_{(\mathbf{i}_{mp};\bfk)}$ (the notation is explained
after (\ref{eqmygen})). \Box

 For $|\bfk|>1$ we take
$\bfW_\bfk:= \bfY_\bfk$ and $\mathbf{R}_\bfk:=\mathbf{T}_\bfk$ and
write (\ref{invlapvk}) as
\begin{equation}\label{eqabstractm}
(1+J_\bfk)\bfW_\bfk=\hat Q_\bfk^{-1} \bfR_\bfk
\end{equation}

\z with $\hat Q_\bfk:=(-\hat\Lambda+p+\bfk\cdot\bflam)$ (notice that
for $|\bfk|>1$ and $p\in\mathcal{S}_0'$
we have $\mathrm{det}\,\hat Q_\bfk(p)\ne 0 $).

\begin{gather}\label{defjm}
(J_\bfk\bfW)(p):=\hat Q_\bfk^{-1}\left(\left(\hat
B+\bfm\cdot\bfk\right)\int_0^p\bfW(s)\mathrm{d}s-\sum_{j=1}^n
\int_0^p W_j(s)\bfd_j(p-s)\mathrm{d}s\right)\cr
\end{gather}

The case $|\bfk|=1$ is special in that $p=0$ is a singularity.  The
corresponding statements of Lemma~\ref{higherterms} for $|p|<\epsilon
$ with $\epsilon$ small are proven in Proposition~\ref{localresult}:
let $\mathbf{W}^0_\bfk$ be the functions provided there.  For the
analytic part of Lemma~\ref{higherterms} we need to show that
$\mathbf{W}^0_\bfk$ extend analytically to solutions of (\ref{eqMv}).
To unify the treatment we derive equations of the form
(\ref{eqabstractm}) for these continuations.  Let
$\delta\in{\mathcal{S}_c^+}',|\delta|<\epsilon$.  Using Proposition~\ref{localresult} and standard analyticity
arguments, it suffices to show that $\mathbf{W}^0_\bfk$ extend
analytically in any sector ${\subset\mathcal{S}_c^+}'$ centered at $\delta$,
and that the corresponding convolution equation is satisfied
along the ray $d_\delta\ni\delta$.  We let
$\mathbf{a}=\bfW^0_\bfk(\delta)$, and with $\bfW^1_\bfk=\bfW^0_\bfk$
for $|p|<|\delta|$ along $d_\delta$ and zero otherwise, we write
$\bfY_\bfk(p)=\bfW^1_\bfk(p)+\mathbf{a}+\bfW_\bfk(p-\delta)$. For
$p\in d_\delta$ we find that $\bfW_\bfk$ must satisfy
(\ref{eqabstractm}) where $\hat
Q_\bfk:=(-\hat\Lambda+p+\bfk\cdot\bflam+\delta)$ and ${\bf R}_\bfk(p)$
is given by

$$
\left(\bfm\cdot\bfk+\hat B\right)
  \left(\int_0^{\delta}\bfW^0_\bfk(s)\mathrm{d}s-\mathbf{a}p\right)
+\sum_{j=1}^n 
\int_0^{\delta}(\bfW^0_\bfk)_j(s)\bfd_{\mathbf{e}_j}(p-s)\mathrm{d}s
-\mathbf{a}\mathbf{D}_{\mathbf{e}_j}(p)
$$

\z ${\bf R}_\bfk(p)$ is manifestly analytic in ${\mathcal{S}^+}'$.
Since $\bfW_\bfk(p)=\bfW^0_\bfk(p+\delta)-\mathbf{a}$ is already a
solution of (\ref{eqabstractm}) for small $p$, and in this case the
left side of (\ref{eqabstractm}) vanishes for $p=0$, we have ${\bf
  R}_\bfk\in\mathcal{T}_1$, for $|\bfk|=1$.

Combined with Remark~\ref{stableconv} and induction on $\bfk$, the
following result completes the proof of Lemma~\ref{higherterms}, parts
(i) and (ii).  

\begin{Proposition}\label{Uniformnorm}
 
i) For large $\nu$ and  constants $K_1$ and $K_2(\nu)$ independent of
$\bfk$, with $K_2(\nu)=O(\nu^{-1})$ we have $\|Q_\bfk^{-1}\|\le\frac{K_1}{|\bfk|}$
and
\begin{eqnarray}
  \label{normQk}
 \|J_\bfk\|\le
  K_2(\nu)
\end{eqnarray}

ii)  For large $\nu$, the operators $(1+J_\bfk)$ defined in $\mathcal{D}'_{m,\nu}$, and
also in
$\mathcal{T}_{\bfk\bfbet'-1}$ for $|\bfk|>1$ and
in $\mathcal{T}_1$ for $|\bfk|=1$ are simultaneously invertible.
Given $\bfY_0$ and $\bfC$, the $\bfW_\bfk, \,|\bfk|\ge 1$ are
uniquely determined. For any $\delta>0$ there is a large enough  $\nu$,
so that

\begin{equation}\label{estimunifk}
\|\bfW_\bfk\|\le\delta^{|\bfk|},\ k=0,1,..
\end{equation}

\z (in the $\mathcal{D}'_{m,\nu}$ topology, (\ref{estimunifk}) 
hold uniformly in $\phi\in[\psi_-+\epsilon,0]$ and
$\phi\in[0,\psi_+-\epsilon]$ for any small $\epsilon>0$).

\end{Proposition}

\Box

{\em Proof.} 

(i) 
For $\mathcal{T}_{\bfk\bfbet'-1}$, this follows immediately from
Remark~\ref{substiconv}, and the constants will depend on the
parameter $c$ in $\mathcal{S}_c^+\cap\mathcal{S}_{\RR}$.  Given $\bfY_0$, the estimates
(\ref{normQk}) are also true in $\mathcal{D}'_{m,\nu;\phi}$ this time
{\em uniformly} in $\phi$, down to $\phi=0$. The proof of
(\ref{normQk}) in this case is given in Appendix~\ref{sec:A1}, in
Lemma~\ref{cinftycase}, from which the continuity of $J_\bfk$ in
$\phi$ also follows.

(ii) From (\ref{eqabstractm}) and (i) we get, for some $K$ and $j\ge
1$ $\|\bfW_\bfk\|\le K\|\bfR_\bfk\|$.  We first show inductively that
the $\bfW_\bfk$ are bounded. Choosing a suitably large $\nu(\epsilon)$
we can make $\max_{|\bfk|\le 1}\|\bfW_\bfk\|_\nu\le \epsilon$ for any
positive $\epsilon$ (uniformly in $\phi$).  We show by induction that
$\|\bfW_\bfk\|_\nu\le \epsilon$ for all $k$.  Using (\ref{estimunifk}),
(\ref{eqMv}), (\ref{combineq}), Proposition~\ref{estimd} and the crude
estimate $\binom{a}{ b}\le 2^a$ we get

\begin{equation}\label{finesti1}\|\bfW_\bfk\|_\nu\le K\|\bfR_\bfk\|_\nu\le
\sum_{\bfl\le\bfk}\kappa_1^{|\bfl|}
 \epsilon^{|\bfk|}\sum_{(\bfii_{mp})}1\le \epsilon^{|\bfk|}
\sum_{s=0}^{|\bfk|} \kappa_1^s
2^{n_1(|\bfk|+s)}2^{s+n_1}\le (C_1 \epsilon)^{|\bfk|}
\end{equation}

\z where $C_1$ does not depend on $\epsilon,\bfk$. Choosing $\epsilon$
so that $\epsilon<C_1^{-2}$ we have, for $|\bfk|\ge 2$ $(C_1
\epsilon)^{|\bfk|}<\epsilon$ completing the induction step . But as we
now know that
$\|\bfW_\bfk\|_\nu\le \epsilon$,  the same 
inequalities (\ref{finesti1}) show that in fact
$\|\bfW_\bfk\|_\nu\le (C_1\epsilon)^{|\bfk|}$. Choosing $\epsilon$
small enough, the first part of Proposition~\ref{Uniformnorm}, (ii)
follows.
\Box

\begin{Proposition}\label{SRY0}
i) Let $\bfY_0$ be given by Lemma~\ref{analyticase}. We
have, in $\mathcal{D}'_{m,\nu}(\RR^+)$

\begin{eqnarray}
  \label{firstdecY}
  \mathbf{Y}_0^{\pm}=\mathbf{Y}_0^{\mp}+\sum_{k=1}^{\infty}(\pm
  S_1)^k\left(\bfY^{\mp}_{k\mathbf{e}_1}(p-k)\bchi_{[k,\infty]}\right)
  ^{(mk)}
\end{eqnarray}

\z (cf. (\ref{defSj})) and  $\mathbf{Y}_0$ is analytic in $\mathcal{R}_1$.

ii) The
general solution of (\ref{eqil}) in $\mathcal{D}'_{m,\nu}(\RR^+)$
is 
\begin{eqnarray}
  \label{SGEQIL}
  \bfY_0^++\sum_{k=1}^{\infty}{C}^k\left(\bfY^+_{k\mathbf{e}_1}(p-k)\bchi_{[k,\infty]}\right)
^{(mk)}
\end{eqnarray}

\z with arbitrary $C$ (a similar statement
holds with $\bfY_0^-$ replacing $\bfY_0^+$). 

\end{Proposition}

{\em Proof}.

\z We start with (ii). Assuming first (\ref{SGEQIL}) is indeed a
solution of (\ref{eqil}), to see that there are no others, it suffices
by Proposition~\ref{contidis} to check that (\ref{SGEQIL}) is the
general solution on $[0,1+\epsilon)$.  The latter part is immediate
from Remark~\ref{coinci} and Proposition~\ref{localresult} below.  Now
$\bfY_0^+$ is a solution of (\ref{eqil}), by Lemma~\ref{analyticase}
(ii); the sum (\ref{SGEQIL}) is convergent in $\mathcal{D}'_{m,\nu}$
by (\ref{finesti1}). Since
$\left(\bfY_{k\mathbf{e}_1}(p-k)\bchi_{[k,\infty]}\right)
^{(mk)}\in\mathcal{D}'_{m,\nu}(k,\infty)$, to show that (\ref{SGEQIL})
is a solution on $\RR^+$, we check inductively on $j$ that
$\bfH_j=\sum_{k=0}^{j}\left(\bfY_{k\mathbf{e}_1}(p-k)\bchi_{[k,\infty]}\right)
^{(mk)}$ solves (\ref{eqil}) in $\mathcal{D}'[0,j+1)$.  Assuming this
for $j'<j$ and looking for a solution on $[0,j+1)$ in the form
$\bfH_{j+1}=\bfH_j+\left(\tilde{\bfY}_{j}(p-j)\bchi_{[j,\infty]}\right)
^{(mj)}$ we obtain, by a straightforward calculation, using the
induction hypothesis and (\ref{u4})

\begin{multline}
    \label{eqykj}
     (\hat\Lambda-p) \tilde{\bfY}_j^{(mj)}+
\hat B\mathcal{P}\tilde{\bfY}_j^{(mj)}-\sum_{j=1}^n
(\tilde{\bfY}_j)_j^{(mj)}*\bfd_{j\mathbf{e}_j}\cr =
\sum_{|\bfl|> 1}\bfdl *
\sum_{\Sigma s=j}*\prod_{i=1}^n*\prod_{j=1}^{l_i}(\tilde{\bfY}_{s_{i,j}})^{(s_{i,j})}_i=:
\bfR_j(p)
  \end{multline}

  \z which integrated $(mj)$ times is exactly the equation for
  $\bfY_{j\mathbf{e}_j}$, cf. also \S~\ref{sec:For}. The claim now
  follows from Lemma~\ref{higherterms}, (iii).  For (i), we note as
  before that $\bfY_0^{\pm}$ are indeed solutions of (\ref{eqil}).
  Applying (ii), we only need to identify $C$ for which purpose we
  compare the left side with the right side on $(1,1+\epsilon)$, where
  all the terms except for $k=1$ vanish and Remark~\ref{identif} below
  applies.  Lemma~\ref{STR1} completes the proof.

\Box

{\em Proof of Lemma~\ref{higherterms}, (iv)}.
 We 
 let  $\bfY_{j}=\bfY_{j\mathbf{e}_j}$.  By Lemma~\ref{STR1},
  $\bfd_{j\mathbf{e}_j}\in\mathcal{F}(\mathcal{R}'_1)$ since

\begin{multline}\label{exprdj}
  \bfd_\bfj^+=\sum_{\bfl\ge\bfj}\binom{\bfl}{\bfj}\bfg_{\bfl}*(\bfY_0^{+})^{*(\bfl-\bfj)}=
  \sum_{\bfl\ge\bfj}\bfg_{\bfl}*\Bigg[\bfY_0^{-}\cr+\sum_{s\ge
    1} S^s(\bfY_s^{-}(p-s))^{(ms)}\Bigg]^{*(\bfl-\bfj)}=
  \sum_{\bfk\succeq 0}\left(\sum_{s\ge 1}
  S^s(\bfY_s^{-})^{(ms)}\right)^{*\bfk}\bfQ_{\bfk\bfj}\cr =
  \sum_{l=0}^{\infty} S^l(\bfd_{\bfj; l}^-)^{(ml)}
\end{multline}

\z where  $\bfQ_{\bfk\bfj}=\sum_{\bfl\ge\bfj+\bfk}
\binom{\bfl}{\bfj}\binom{\bfl-\bfj}{\bfk}
\bfg_\bfl(\bfY_0^-)^{*(\bfl-\bfk-\bfj)}$ and

\begin{eqnarray}\label{definition dj}
&&(\bfd_{\bfj; l}^-)^{(ml)}=
\sum_{\begin{subarray}\bfk\in\NN^n\cr|\bfk|\in[0,l]\end{subarray}}\bfQ_{\bfk\bfj}*
    \sum_{(i_{rs}:l)} 
\sideset{^*}{}\prod_{r=1}^n\sideset{^*}{}
\prod_{s=1}^{k_r}\left(\bfY^{-}_{i_{rs}}\right)^{(mi_{rs})}_r\cr&&=
\left(\sum_{\begin{subarray}\bfk\in\NN^n\cr|\bfk|\in[0,l]\end{subarray}}\bfQ_{\bfk\bfj}*
    \sum_{(i_{rs}:l)} 
\sideset{^*}{}\prod_{r=1}^n\sideset{^*}{}
\prod_{s=1}^{k_r}\left(\bfY^{-}_{i_{rs}}\right)_r\right)^{(ml)}\cr&&
 \end{eqnarray}

 \z and the notation $\sum_{(i_{rs}:l)}$ is explained after Eq.
 (\ref{eqmygen})); in particular, (1) there are only finitely many
 terms in (\ref{definition dj}) and (2) by homogeneity,
 $\mathbf{D}_{\bfj;l}\in\mathcal{T}_{l\beta'_1-1}$. In addition, it
 follows, as in (\ref{finesti1}), (noting that we only need the
 $\bfj$ with $|\bfj|=1$) that $\|\mathbf{Q}_{\bfk;\bfj}\|\le
 K(4\kappa)^{|\bfk|}$ and $\|\mathbf{D}_{\bfj;l}\|\le K
 (2^{2n+3}\kappa C_1\epsilon)^l$.  If we look for $\bfY_\bfk^{+}$ in
 the form $\mathbf{Y}_\bfk^-+\sum_{l=1}^{\infty}(\bfY^{-}_{\bfk;
   l}(p-l))^{(ml)}$ then the equation for $\bfY^-_{\bfk; l}$, $l\ge
 1$, reads

\begin{multline}
  \label{eqcomplk}
  \left(-p+\hat{\Lambda}-\bfk\cdot\bflam-l\right)\bfY^-_{\bfk;l}+
\left(\hat{B}+\bfm\cdot\bfk+m_1 l\right)
\mathcal{P}\bfY^-_{\bfk;l}\cr+\sum_{|\bfj|=1}
\bfd^{-}_{\bfj;0}*\left(\bfY^-_{\bfk;l}\right)^\bfj
=\left(\bfT_\bfk(\bfY^-_{\cdot;l})\right)-
\sum_{s=1}^{l-1}\sum_{|\bfj|=1}
\bfd_{\bfj;s}*\left(\bfY^{-}_{\bfk;l-s}\right)^{\bfj}
\end{multline}

\z By induction, exactly as in (iii), it follows that
$\mathbf{Y}_{\bfk;l}(\cdot\,\mathrm{e}^{\mathrm{i}\phi})\in\mathcal{T}_{\bfk\bfbet'+l\beta'_1-1}$,
are $\phi-$ continuous in $\mathcal{D}'_{m,\nu}(\RR^+)$ and that, with $\nu$
large enough independent of $\bfk,l,\phi$,
$\|\mathbf{Y}_{\bfk;l}\|_{\mathcal{D}'_{m,\nu}}\le
K\delta^{|\bfk|+l}$. Analyticity in $\mathcal{R}'_1$ follows now from
Lemma~\ref{STR1}.

(v) Laplace transformability as well as the fact that $\bfy_\bfk$
solve (\ref{systemformv}) follow immediately from (\ref{estimunifk})
and Lemma~\ref{existe}. Uniform convergence follows from
(\ref{normVK}) and Lemma~\ref{existe}. Let
$\bfy_\bfk^{\pm}=\lap\bfY_\bfk^{\pm}$. Now since $\bfy_0^{\pm}+\sum
x^{\bfm\cdot\bfk}\mathbf{C}^{\bfk}\mathrm{e}^{-\bfk\cdot\bflam x} \bfy_k^{\pm}$
formally solves (\ref{eqor}) (by the very construction of
(\ref{systemformv}), see \S\ref{sec:For}) and is a uniformly
convergent function series, the conclusion follows together with
the fact that $\lap\bfY_\bfk\sim \tilde{\bfy}_\bfk$  since by (ii),
and (iv) $\lap\bfY_\bfk$ have power series
asymptotics which by construction must be formal solutions of
(\ref{systemformv}).

 \Box 

 {\em Proof of Theorem~\ref{AS}, (ii)} \z Modulo relabeling of the
 spatial directions and rescaling, the proofs of the properties of the
 solution $\mathbf{Y}_0$ of (\ref{eqil}) are valid in a sector
 centered around any eigenvalue $\lambda_i$.  The above proofs of the
 analytic properties work with virtually no change along any direction
 so that $p+\bfk'\cdot{\bflam}-\lambda_i\ne 0$ and the same is true
 with respect to the $\mathcal{D}'_{m,\nu}$ properties of
 $\mathbf{Y}_\bfk$ restricted to {\em compact sets} in $\CC$.  The
 analysis in $\mathcal{D}'_{m,\nu}(e^{i\phi}\RR^+)$ (i.e., along the
 infinite ray) and the associated {\em norm} estimates are {\em not}
 valid along directions so that $\Re(p-\lambda_i+\bfk'\cdot{\bflam})\le 0$
 (because of possible presence of small denominators in
 $(1+J_\bfk)^{-1}$, the only difference in this case, but a
 significant one). However,  we do not take Laplace
 transforms along such directions (cf. (c1)) and the associated
infinite ray norms are not
 needed for our purposes.  For the analytic properties in the
 neighborhood of $-\bfk'\cdot{\bflam}+\lambda_i$, see
 Remark~\ref{wrongdirection}.  \Box

\subsubsection{Higher resurgence relations}
\label{sec:str-sing}

\begin{Proposition}\label{asymptrick} i)  Let $\bfy_1$ and $\bfy_2$ be solutions of (\ref{eqor})
so that $\bfy_{1,2}\sim\tilde{\bfy}_0$ for large $x$ in an
open sector $S$ (or in some direction
$d$);
then $\bfy_1-\bfy_2=\sum_{j}C_j\mathrm{e}^{-\lambda_{i_j}
  x}x^{-\beta_{i_j}}(\mathbf{e}_{i_j}+o(1))$ for some constants $C_j$, where the indices run over
  the eigenvalues $\lambda_{i_j}$ with the property $\Re(\lambda_{i_j}
  x)>0$ in $S$ (or $d$). 
If $\bfy_1-\bfy_2=o(\mathrm{e}^{-\lambda_{i_j}
    x}x^{-\beta_{i_j}})$ for all $j$, then $\bfy_1=\bfy_2$.
 
ii) Let $\bfy_1$ and $\bfy_2$ be solutions of
(\ref{eqor}) and assume that $\bfy_1-\bfy_2$ has differentiable
asymptotics of the form 
$\mathbf{K}a\exp(-ax)x^b(1+o(1))$ with  $\Re(ax)>0$ and $\mathbf{K}\ne 0$, for large $x$.
Then $a=\lambda_i$ for some $i$.

iii) Let $\mathbf{U}_\bfk\in\mathcal{T}_{\{\cdot\}}$ for all $\bfk$,
$|\bfk|>1$.  Assume in addition that for
large $\nu$ there is a function $\delta(\nu)$ vanishing as
$\nu\rightarrow\infty$ such that

\begin{gather}
\label{cd1}
\sup_{\bfk}\delta^{-|\bfk|}\int_{d}\left|\mathbf{U}_{\bfk}(p)\mathrm{e}^{-\nu p}\right|\mathrm{d}|p|<K<\infty
\end{gather}

\z Then, if $\bfy_1,\bfy_2$ are solutions of (\ref{eqor}) in $S$ where in addition

\begin{eqnarray}
  \label{lapcond}
  \bfy_1-\bfy_2=\sum_{|\bfk|>1}
\mathrm{e}^{-\bflam\cdot\bfk x}x^{\bfm\cdot\bfk}\int_d\mathbf{U}_\bfk(p)
\exp(-xp)\mathrm{d}p
\end{eqnarray}

\z where $\bflam,x$ are as in (c1), then $\bfy_1=\bfy_2$, and
$\mathbf{U}_\bfk=0$ for all $\bfk$, $|\bfk|>1$.

\end{Proposition}

{\em Proof}. (i) is a classical result (see \cite{Iwano}
for the general treatment and \cite{Wasow} for a brief presentation of
special cases and further references).  However, what
is actually needed for
our purposes can  be reduced to the more familiar {\em linear}
asymptotic theory in the following way. Let $d$ be a direction in the
complex plane and let $\bfy_0$, $\bfy_1$ be solutions of (\ref{eqor})
such that $\bfy_{0,1}\sim\tilde{\bfy}_0$ for large $x$ along $d$.
Then, by (n5), $\bfy_{0,1}(x)=O(x^{-M})$ and for any $j$,
${\bf g}^{(\mathbf{e}_j)}(x,\bfy_{0,1}(x))=O(x^{-M})$.  If
$\boldsymbol\delta=\bfy_{1}-\bfy_{0}$ then by hypothesis
$\boldsymbol\delta(x)=o(x^{-l})$ along $d$, for all $l$. The function
$\boldsymbol\delta$ is locally analytic and satisfies the  equation

\begin{multline}
  \label{eqdel}
  \bfdel'=-\hat{\Lambda}\bfdel-\frac{1}{x}\hat{B}\bfdel
+\sum_{|\bfk|=1}{\bf g}^{(\bfk)}(x,\bfy_0)\bfdel^{\bfk}
+\sum_{|\bfk|>1}{\bf g}^{(\bfk)}(x,\bfy_0)\bfdel^{\bfk}=\cr
-\hat{\Lambda}\bfdel-\frac{1}{x}\hat{B}\bfdel
+\frac{1}{x^M}\sum_{j=1}^n(\bfdel)_j\bfh_{\mathbf{e}_j}(x)
\end{multline}

\z where $\bfh_\bfk(x)$ are bounded along $d$. Obviously, because of
the link between $\bfdel$ and $\bfh_\bfk$, the $\bfdel$ we started
with might be the only solution of (\ref{eqdel}) which is also a
difference of solutions of (\ref{eqor}).  The asymptotic
characterization we need holds nevertheless for {\em all} decaying
solutions solutions of (\ref{eqdel}): since no two eigenvalues are
equal, there exists by the well-known linear asymptotic theory
\cite{Wasow} a fundamental set $\{\bfdel_i\}_{1\le i\le n}$ of
solutions of (\ref{eqdel}) such that $\bfdel_i\sim \mathrm{e}^{-\lambda_i x}
x^{-\beta_i}(\mathbf{e}_i+o(1))$.  Thus $\bfdel=\sum_{i=1}^n
C_i\bfdel_i=\sum_{i=1}^n C_i \mathrm{e}^{-\lambda_i x}
x^{-\beta_i}(\mathbf{e}_i+o(1))$. Since $\Re(-\beta_i)>0$ and the
$\lambda_i$ are distinct we must have $C_i=0$ for all $i$ for which
$\Re(-\lambda_i x)\ge 0$, otherwise $|\bfdel(x)|$ would be unbounded
for large $x$; the first part of (i) is
proven. If on the other hand $\bfdel=o(\mathrm{e}^{-\lambda_{i_j}
  x}x^{-\beta_{i_j}})$ for all $j$, again because the $\lambda_i$ are
independent, it follows that $C_i=0$ for all $i=1,2,\ldots,n$, thus
$\bfdel=0$.

(ii) is now obvious.

For (iii), note first that by (\ref{cd1}) and (c1) the RHS of
(\ref{lapcond}) converges uniformly for large $x$ in some open sector.
In addition, by an arbitrarily small change in $\xi=\arg(x)$, we can
make the set $\{\Re(x\lambda_i)\}_i$ $\ZZ$-independent (the existence
of $\bfk(\xi)\ne 0$ s.t. $\Re(\mathrm{e}^{\mathrm{i}\xi}\bfk\cdot\bflam)=0$ for $\xi$ in an
interval of would imply the existence of a {\em common} $\bfk$ for a
set of $\xi$ with an accumulation point, giving $\bfk\bflam=0$). We
choose such a $\xi$.  Assume now there exist $\bfk$ so that
$\mathbf{U}_\bfk\ne 0$; among them let ${\bfk_0}$ have the least
$\Re(x\bfk\cdot\bflam)$.  By (\ref{cd1}) for large $x$,
$\bfy_1-\bfy_2\sim \mathrm{e}^{-\bflam\cdot\bfk_0 x} x^{\bfm\cdot\bfk_0}\lap_\phi
\mathbf{U}_{{\bfk_0}}(1+o(1))$. Because
$\mathbf{U}_{{\bfk_0}}\in\mathcal{T}_{\{\cdot \}}$, and by (\ref{cd1}),
$\lap_\phi \mathbf{U}_{\bfk_0}$ has a differentiable power series
asymptotics which is the term-by-term Laplace transform of the Puiseux
series at the origin of $\mathbf{U}_{{\bfk_0}}$, and thus non-zero. This
contradicts (i) because with
$|{\bfk_0}|>1$ we have $\bflam\cdot{\bfk_0} \ne\lambda_j$ for all $j$ 
($\ZZ-$independence).
Thus $\mathbf{U}_\bfk=0$ for all $\bfk$.

\Box

We let $\bfC\in(\CC\backslash\{0\})^{n_1}$ be an arbitrary
 constant vector. 

\z For $x$ large enough, $\bfy^+$ defined in
(\ref{solupperlower})  is a solution of
(\ref{eqor}) in an open sector containing $x$. We now use
Lemma~\ref{higherterms} to write (\ref{solupperlower}) in terms of
functions analytic in the lower half plane:

\begin{multline}
  \label{solupper-}
  \bfy^+ =\lap\bfY_0^-+\sum_{j=1}^{\infty} x^{mj}\mathrm{e}^{-jx}
  \lap\bfY_{0;j}^-+ \sum_{\bfk\succ 0,j\ge
    0}x^{\bfm\cdot\bfk+mj}\bfC^{\bfk}\mathrm{e}^{-(\bfk\cdot\bflam+j\lambda_1)
    x} \lap \bfY^-_{\bfk;j}\cr =\lap\bfY_0^-+\sum_{\bfk\succ
    0}x^{\bfm\bfk}\mathrm{e}^{-(\bfk\cdot\bflam)
    x}\sum_{\bfk';j:\bfk'+j\mathbf{e}_1=\bfk}\bfC^{\bfk'}\lap
  \bfY^-_{\bfk';j}
\end{multline}

\z where, by (\ref{firstdecY}) we have
$\bfY_{0;j}=S_1^j\bfY^-_{j\mathbf{e}_1}$.
On the other hand, the expression 

\begin{eqnarray}
  \label{sollower}
 \bfy^- =\lap\bfY_0^-+\sum_{\bfk\succ 0}x^{\bfm\cdot\bfk}\tilde{\bfC}^{\bfk}\mathrm{e}^{-\bfk\cdot\bflam x}
  \lap \bfY^-_\bfk
\end{eqnarray}

\z is, for any $\mathbf{\tilde{C}}$, 
a solution of (\ref{eqor}) as well. Choosing
$\tilde{C}_1=C_1+S_1;\tilde{C}_i=C_i;\ (i>1)$ all the
terms with $|\bfk|\le 1$ 
in (\ref{sollower}) and (\ref{solupper-}) coincide and thus

\begin{eqnarray}
  \label{difey}
  \bfy^+-\bfy^-=\sum_{|\bfk|>
    1}x^{\bfm\cdot\bfk}\mathrm{e}^{-\bfk\cdot\bflam x}\lap\mathbf{U}_\bfk
  \cr&&
\end{eqnarray}

\z where 
\begin{eqnarray}\label{defU}
\mathbf{U}_\bfk=
\sum_{\bfk'+j\mathbf{e}_1=\bfk}\bfC^{\bfk'}
\bfY^-_{\bfk';j}-\tilde{\bfC}^{\bfk}\bfY^-_\bfk
\end{eqnarray}

\z so that  
applying Proposition~\ref{asymptrick} (ii)

\begin{eqnarray}
  \label{u=0}
 \mathbf{U}_\bfk=0 
\end{eqnarray}

\z Since $C_i\ne 0$ we have, for any $C_1$,
\begin{eqnarray}
  \label{nearfinalresu}
  (C_1+S_1)^{k_1}\mathbf{Y}^-_{\bfk}=
\sum_{\bfk'+j\mathbf{e}_1=\bfk}C_1^{k'_1}
\bfY^-_{\bfk';j}
\end{eqnarray}

\z with arbitrary  $C_1$ so that

\begin{eqnarray}
  \label{primaresu}
 \bfY^-_{\bfk;j}= \binom{k_1+j}{j}S_1^{j}\mathbf{Y}^-_{\bfk+j\mathbf{e}_1}
\end{eqnarray}

\z Combined with the definition
of $\bfY^-_{\bfk;j}$ 
this gives (\ref{thirdresu}).

Solving for $\mathbf{Y}^-_{\bfk+j\mathbf{e}_1}$, (\ref{primaresu})
determines later series in the transseries in terms of earlier ones. The same
arguments work of course with $-/+$ and $+S_1/(-S_1)$ interchanged.

Theorem~\ref{CEQ} part (iii) follows from the following.
 
\begin{Proposition}\label{classic} Any solution $\bfy$ of (\ref{eqor}) so that
  $\bfy\sim\tilde{\bfy}_0$ along some direction $d\subset S_x$ is of
  the form (\ref{solupperlower}), for a unique $\bfC^+$
  (a similar statement holds with $+/-$ interchanged).  Alternatively,
  a solution $\bfy$ of (\ref{eqor}) so that $\bfy\sim\tilde{\bfy}_0$
  along some direction $d\subset S_x$ can be represented as
  (\ref{soleqn}) or more generally as (\ref{soleqnpa}) where Laplace
  integration is along $\RR^+$ (in distributions), for a unique
  $\bfC$.

\end{Proposition}

{\em Proof.} Let $\bfy$ be an arbitrary solution of (\ref{eqor}) so
that $\bfy\sim\tilde{\bfy}_0$ along $d\subset S_x$. Then, by
Proposition~\ref{asymptrick},
$\bfy-\bfy^+=\sum_{j}\overline{\mathbf{C}}_j\mathrm{e}^{-\lambda_{i_j}
  x}x^{-\beta_{i_j}}(\mathbf{e}_{i_j}+o(1))$ for some constant
$\overline{\mathbf{C}}$. Therefore $\bfy_1$ defined as the ``+''
solution in (\ref{solupperlower}) with
$\mathbf{C}={\overline{\mathbf{C}}}$ will have the property
$\bfy_1-\bfy=o(C_j\mathrm{e}^{-\lambda_{i_j} x}x^{-\beta_{i_j}})$ for
all $j$, hence $\bfy=\bfy_1$. Formula (\ref{soleqnpa}). The last part,
as well as the middle formula in (\ref{microsto}), follow through a
straightforward calculation from the first part, (\ref{primaresu})
using (\ref{normVK}) to control convergence.

\Box

{\em Proof of theorem~\ref{Stokestr}}. Let $\bfy^{\pm}$ be defined by
(\ref{solupperlower}) with $\bfC=(\pm\frac{1}{2} S_1+C)\mathbf{e}_1$,
respectively.
The same arguments leading to (\ref{u=0}) show that
$\bfy^+=\bfy^-=:\bfy$.  All the exponentials in the transseries of
$\bfy$ are generated by construction by $e^{-\lambda x}$.  Choosing
$p$ in the path of integration above/below $\RR^+$ and consequently the $+/-$
representation (\ref{solupperlower})
of $\bfy$ we have by Lemma~\ref{higherterms}  that
$\lap\bfY_\bfk^{\pm}\sim\tilde{\bfy}_\bfk$ in (\ref{solupperlower}),
in the half plane $\Re(xp)>0$.
By construction
$\lap\bfY_{\mathbf{e}_1}=x^{-\beta'_1}(\mathbf{e}_1+o(1))$ 
(cf. \ref{eqneuman}) while for $j>1$ we have
$\lap\bfY_{j\mathbf{e}_1}\sim x^{-j\beta'_1}$
by Lemma~\ref{higherterms} (ii).  The condition
$|x^{-\beta_1+1}e^{-x\lambda_1}|\rightarrow 1$ together with
Lemma~\ref{higherterms} guarantee the uniform convergence of the series
(\ref{solupperlower}). The conclusion is immediate.

\Box

\begin{subsubsection}{Local analysis near $p=1$.}
\label{sec:A1}

\z  
We treat (\ref{eqil}) near $p=1$ as a perturbation of a differential equation
having  the same type of singularity. The associated differential
equation splits the singularity, and our convolution equation
is a regular perturbation of it, which is then solved by 
fixed point methods.
Let $\bfY_0$ be the unique solution in $\mathcal{A}_{z,l}$ of (\ref{eqil}) and let
$\epsilon>0$ be small.  Define

\begin{eqnarray}\label{defbfh}
&&\bfH(p):=\left\{\begin{array}{cc}
\bfY_0(p)\ \ \mbox{for $p\in\mathcal{S}_0$\,,$|p|<1-\epsilon$}\cr
0\ \ \ \mbox{otherwise}
\end{array}\right.\cr&&\mbox{and}\ \ \bfW(1-p):=\bfY_0(p)-\bfH(p)
\end{eqnarray}

\z ($\bfY(p)-\bfH(p)=\bfW(p-1)$ would be  more ``natural'',
but would later complicate notations).
In terms of $\bfW$, for real $z=1-p, z<\epsilon$,  (\ref{eqil}) becomes:
\begin{eqnarray}\label{eq002}
-(1-z){\bf W}(z)=\mathbf{F}_1(z)-\hat\Lambda \bfW(z)+
\hat B\int_{\epsilon}^z\bfW(s)\mathrm{d}s+\calnb({\bf H}+{\bf W})
\end{eqnarray}

\z where 

$${\bf F}_1(1-s):=\bffz(s)-{\hat B}\int_0^{1-\epsilon}\bfH(s)\mathrm{d}s
$$

\begin{Proposition}\label{P8}

i) For small $\epsilon$, $\bfH^{*\bfl}(1+z)$  extends to an analytic
function in the disk
$\calv:=\{z:|z|<\epsilon\}$. Furthermore,
for any $\delta$ there is an $\epsilon$ and a
 constant $K_1:=K_1(\delta,\epsilon)$ such that for 
$z\in\calv$ 

\begin{eqnarray}\label{estHl2}
|\bfH^{*\bfl}(1+z)|<K_1\delta^{|\bfl|} 
\end{eqnarray}

ii)
The  equation (\ref{eq002}) can be written as

\begin{gather}\label{formN}
-(1-z){\bf W}(z)={\bf F}(z)-\hat\Lambda \bfW(z)+
\hat B\int_{\epsilon}^z\bfW(s)\mathrm{d}s-\sum_{k=1}^n
\int_\epsilon^zh_j(s)\bfd_j(s-z)\mathrm{d}s
\end{gather}

\z where

\begin{eqnarray}\label{defD1}
{\bf F}(z):=
\calnb(\bfH)(1-z)+{\bf F}_0(z)
\end{eqnarray}

\begin{eqnarray}\label{defderiv}
\bfd_j=
\sum_{|\bfl|\ge 1}l_j{\bf G_{l}}*\bfH^{*\bar\bfl^j}+
\sum_{|\bfl|\ge 2}l_j{\bf g_{\rm 0,\bf l}}\bfH^{*\bar\bfl^j};\ \bar\bfl^j:=(l_1,l_2,..(l_j-1),..l_n)
\end{eqnarray}

\z  extend to analytic functions in $\calv$.
Moreover, if $\bfH$ is a vector in $\lone_\nu(\RR^+)$ then, for large 
$\nu$, $\bfd_j\in\lone_{\nu}(\RR^+)$ and the functions
$ {\bf F}(z)$ and ${\bf D}_j$ extend to analytic functions in
$\calv$.
Furthermore, ${\bf D}_j\in\mathcal{A}_{z,M}$.

iii) Near $p=1$ we have (cf.  Lemma~\ref{analyticase})

\begin{eqnarray}
  \label{Y0p=1}
  \mathcal{P}^{m+1}\bfY_0&=&(1-p)^{\beta'}\mathbf{A}+\mathbf{B}\ \ (\beta\notin\ZZ)\cr
 \mathcal{P}^{m+1}\bfY_0&=&(p-1)\left(\ln(p-1)\mathbf{A}(p)+\mathbf{B}(p)\right)\ \ (\beta\in\ZZ)
\end{eqnarray}

\z where  $\bf A,B$ analytic at $p=1$.
\end{Proposition}
\smallskip

{\em Proof}.

 Parts (i) and (ii), except for the last claim, are proven
in \cite{Costin}, Propositions 18
and 19. To see that $\mathbf{D}_j\in\mathcal{A}_m$ it is enough to remark
again
that $\mathcal{A}_{z,M}$ is a convolution ideal of $\mathcal{A}_{z,0}$ and that
$\mathbf{G}_{\bf j}\in\mathcal{A}_{z,M}$ for
 $|\mathbf{j}|=1$, by (n5).

For (iii), consider again  equation (\ref{formN}).
Let $\ga=\hat\Lambda-(1-z){\hat{1}}$, where
$\hat{1}$ is the identity matrix. By construction
$\ga$ and $\hb$ are diagonal, $\ga_{11}=z$ and $\hb_{11}=\beta_1=:\beta$.
We write this as $\ga=z\oplus\gc(z)$ and similarly,
$\hb=\beta\oplus\bc$, where $\gc$ and $\bc$
are $(n-1)\times(n-1)$ diagonal matrices. $\gc(z)$ and
$\gc^{-1}(z)$ are analytic in $\calv$.

\z Let

\begin{equation}\label{defQ}
\bfQ:=\mathcal{P}_\epsilon^{m+1} \bfW
\end{equation}

\z with $\mathcal{P}_\epsilon:=
\bfW\mapsto\left(z\mapsto\int_\epsilon^z\bfW(s)\mathrm{d}s\right)$. By
Lemma~\ref{analyticase}, (i),
$\bfQ$ is 
analytic in $\calv\cap(\{z:z+1\in\mathcal{S}_0\})$. From (\ref{formN}) we obtain

\begin{multline}\label{difeqbfQ}
\zpp\bfQ^{(m+1)}(z)-\betapp \bfQ^{(m)}(z)\cr={\bf F}(z)-
\sum_{j=1}^n\int_\epsilon^z\bfd_j(s-z)Q_j^{(m+1)}(s)\mathrm{d}s
\end{multline}

\z or, after $m$ integrations by parts in the r.h.s. 
of (\ref{difeqbfQ}), by Proposition~\ref{P8} (ii), we get

\begin{multline}\label{difeqbfQ2}
  \zpp\bfQ^{(m+1)}(z)-\betapp \bfQ^{(m)}(z)\cr={\bf
    F}(z)-(-1)^{m}\sum_{j=1}^n\int_\epsilon^z\bfd_j^{(m)}
  (s-z)Q_j'(s)\mathrm{d}s
\end{multline}

\z so that with $\beta'=\beta'_1, \bcp=m_1+\bc$,

\begin{multline}\label{difeqbfQ2i}
\zpp\bfQ'(z)-(\beta'\oplus\bcp) \bfQ(z)\cr=\mathcal{P}_\epsilon^{m}{\bf
F}(z)-(-\mathcal{P}_\epsilon)^{m}\sum_{j=1}^n\int_\epsilon^z\bfd_j^{(m)}(s-z)Q_j'(s)\mathrm{d}s\cr
=\mathbf{P}(z)+\sum_{j=1}^n\int_\epsilon^z\bfd_j'(s-z)Q_j(s)\mathrm{d}s
\end{multline}

\z where $\mathbf{P}(z)=\mathcal{P}_\epsilon^{m}{\bf F}(z)$.
With the notation
$(Q_1,{\bf Q}_\perp):=(Q_1,Q_2,..,Q_n)$
we write the system in the form

\begin{align}\label{syt1}
(z^{-{\beta'}} Q_1(z))'&=z^{-{\beta'}-1}\left(P_1(z)+
\sum_{j=1}^n\int_\epsilon^z D_{1j}'(s-z)Q_j(s)\mathrm{d}s\right)\cr
(\mathrm{e}^{\hatc(z)}\Qp)'&=\mathrm{e}^{\hatc(z)
}\gc(z)^{-1}\left(\mathbf{P}_\perp+\sum_{j=1}^n\int_\epsilon^z
\dpp'(s-z)Q_j(s)\mathrm{d}s\right)\cr
\hatc(z)&:=-\int_0^z \gc(s)^{-1}\bcp(s)\mathrm{d}s
\cr
\bfQ(\epsilon)&=0
\end{align}

\z After integration we get:

\begin{align}\label{eqperp12}
Q_1(z)&=R_1(z)+J_1(\bfQ)\cr
\Qp(z)&=\Rp(z)+J_\perp(\bfQ)
\end{align}

\z with
\begin{align}\label{eqperp2}
J_1({\bf Q})&=z^{{\beta'}}\int_\epsilon^zt^{-{\beta'}-1}\sum_{j=1}^n\int_\epsilon^tQ_j(s)D_{1j}'(t-s)\mathrm{d}s\mathrm{d}t\cr
J_\perp({\bf Q})(z)&:=\mathrm{e}^{-\hatc(z)}
\int_\epsilon^z \mathrm{e}^{\hatc(t)}\gc(t)^{-1}
\left(\sum_{j=1}^n\int_\epsilon^z
\dpp'(s-z)Q_j(s)\mathrm{d}s\right)\mathrm{d}t\cr
\Rp(z)&:=\mathrm{e}^{-\hatc(z)}
\int_\epsilon^z \mathrm{e}^{\hatc(t)}\gc(t)^{-1}
 {\bf F}_\perp(t)\mathrm{d}t
\cr
R_1(z)&=z^{{\beta'}}\int_\epsilon^zt^{-{\beta'}-1}P_1(t)\mathrm{d}t
\ \ \ \ \ \ \ \ \ \ \ \ \ \ \ \ \ \ \ \ ({\beta'}\ne 1)\cr
R_1(z)&=P_1(0)z\ln z+
z\int_{\epsilon}^z\frac{P_1(s)-P_1(0)}{s}\mathrm{d}s
\ \ ({\beta'} = 1)
\end{align}

\z Consider the space $\mathcal{U}_{\beta'}$ given by

\begin{align}\label{funspace}
{\cal U_{\beta'}}&=
\Big\{{\bf Q}\ \mbox{analytic in }\{z:0<|z|<\epsilon, \arg(z)\ne \pi\}:
{\bf Q}=z^{\beta'}\bfA(z)+\bfB(z)\Big\}\\* \intertext{for   $\beta'\ne 1$ and}
{\cal U}_1&=\Big\{{\bf Q}\ \mbox{analytic in }\{z:0<|z|<\epsilon, \arg(z)\ne \pi\}:
{\bf Q}=z\ln z\bfA(z)+z\bfB(z)\Big\}\notag
\end{align}

\z where $\bfA,\bfB$ are analytic in ${\calv}$.
(The decomposition of $\bf Q$ in (\ref{funspace})
is unambiguous since $z^{\beta'}$ and
$z\ln z$ are not meromorphic in $\calv$.) 

The norm

\begin{equation}\label{normT}
\|{\bf Q}\|=\sup\left\{
|\bfA(z)|,|\bfB(z)|:z\in\calv\right\}
\end{equation}

\z makes $\cal U_{\beta'}$ a Banach space.

\begin{Proposition}\label{eqinte} 
The operator $J:=(\bfQ\mapsto (J_1\bfQ,J_\perp\bfQ)$ has
norm 
$O(\epsilon)$, for small $\epsilon$,
in  ${\cal U_{\beta'}}$ as well as in $\lone[-\epsilon,\epsilon]$.
 Along any segment $d_\epsilon$ originating  at $z=\epsilon$ in the region $|z|<\epsilon, \arg(z)\ne \pi$, Equation
(\ref{difeqbfQ2i})
has a unique solution in $L^1_\nu(d_\epsilon)$. This solution belongs to
$\mathcal{T}_\beta$.

\end{Proposition}

\z The proof uses the following
elementary identities:

\begin{eqnarray}\label{convert}
  &&\int_\epsilon^z A(s)s^r \mathrm{d}s ={\mbox{const.}}+z^{r+1}\int
  _{0}^{1}\!{ A}(zt)t^{r}{\mathrm{d}t}={\mbox{const.}}+
  z^{r+1}Analytic(z)\cr &&\int_0^z s^r\ln
  s\,A(s)\mathrm{d}s=z^{r+1}\ln z\int _{0}^{1}\!{
    A}(zt)t^{r}{\mathrm{d}t}+z^{r+1}\int _{0}^{1 }\!{ A}(zt)t^{r}\ln
  t{\mathrm{d}t}\cr &&
\end{eqnarray}

\z where the second equality is obtained  by differentiating
 with respect to $r$  the first equality.

Using
(\ref{convert}) it is straightforward to check that the r.h.s.
of (\ref{eqperp12}) extends to a linear inhomogeneous
operator
on $\cal U_{\beta'}$ with image in $\cal U_{\beta'}$ and that
the norm of $J$ is $O(\epsilon)$ for small $\epsilon$.
For instance, one of the terms in $J$ for $\beta'=1$,

\begin{multline}\label{arr1}
z\int_0^z t^{-2}
\int_0^t s\ln s\,A(s)D'(t-s)\mathrm{d}s\cr= z^{2}\ln z\int _{0}^{1}
\int _{0}^{1}\sigma A(z\tau\sigma)D'
(z\tau-z\tau
\sigma){\mathrm{d}\sigma}{\mathrm{d}\tau}
\cr +
z^{2}\int _{0}^{1}\mathrm{d}\tau\ln\tau
\int _{0}^{1}\mathrm{d}\sigma\,\sigma(1+\ln\sigma)
A(z\tau\sigma)D'(z\tau-z\tau\sigma)
\end{multline}

\z  manifestly in $\cal U_{\beta'}$ if $A$ is
analytic in $\calv$. Comparing with (\ref{funspace}), the extra power 
of $z$
accounts for a norm $O(\epsilon)$ for this term.

Therefore, in (\ref{syt1}) $(1-J)$ is invertible and the solution $\bf
Q\in\cal U_{\beta'}$. In view of the the
uniqueness of the solution of (\ref{eqil}) in $\mathcal{S}_0$,
(Lemma~\ref{analyticase}, (i)) the rest of the proof of
Proposition~\ref{P8}, (iii) is immediate.

\Box

A short calculation shows that:
\begin{Remark}\label{wrongdirection} i) The equation for $\mathbf{Y}_\bfk$
near $z=0$ where  $z=-p-\bfk\cdot\bflam+\lambda_i$ can be written
in the form (\ref{formN}), for a different $\mathbf{F}$, and with
$\hat{B}+\bfk\cdot\bfm$ instead of $\hat{B}$. Thus
$\mathbf{Y}_\bfk(z)\in\mathcal{T}_{\{\cdot\}}$.
 
ii) The equation for $\mathbf{Y}_\bfk$
near $z=0$ where  $z=-p-\bfk'\cdot\bflam+\lambda_i$ where
$\bfk'\prec\bfk$ can be written as 

\begin{eqnarray}
  \label{wrongdir}
  (1+J_\bfk)\mathbf{Y}_\bfk(z)=\mathbf{R}_\bfk(z)
\end{eqnarray}
\z where

$$J_\bfk\bfY=(\hat{B}+\bfm\cdot\bfk)\hat{M}^{-1}\mathcal{P}\bfY+
\hat{M}^{-1}\sum_{|\bfj|=1}\overline{\mathbf{D}}_\bfj*\bfY^\bfj,$$

\z $\hat{M}=z+\hat{\Lambda}-\lambda_i+(\bfk'-\bfk)\cdot\bflam$,
$\mathbf{R}_\bfk=\hat{M}^{-1}\mathbf{T}_\bfk$ and
$\overline{\mathbf{D}}_\bfj$ analytic for small $z$.  Thus, arguments
virtually identical  to those for (\ref{eqperp12}) imply that
$\mathbf{Y}_\bfk\in\mathcal{T}_{\{\cdot\}}$ near these points.

\end{Remark}

\end{subsubsection}

\begin{subsubsection}{The solutions of (\ref{difeqbfQ2}) on
    $[-\epsilon,\epsilon]$}
\label{sec:-e,e}

Let $\bfQ_0$ be the solution given by Proposition \ref{eqinte}, take
$\epsilon$ small enough and denote by $\cal O_\epsilon$ a neighborhood
in $\CC$ of width $\epsilon$ of the interval $[0,1+\epsilon]$. We look
for solutions of (\ref{difeqbfQ2i}) in $\lone[-\epsilon,\epsilon]$.
The main difference with respect to the previous section is that in
integrating (\ref{difeqbfQ2i}) to the analog of (\ref{eqperp12})for
negative $z$, the constant of integration will now be undetermined leading
to a one-parameter family of solutions. See also Remark~\ref{lincomb}
below.

\begin{Remark}\label{uniformlb}. 
 As $\phi\rightarrow\pm 0$,
$\bfQ_0(z\mathrm{e}^{\mathrm{i}\phi})\rightarrow\bfQ_0^\pm(z)$ in the
sense of $\lone([0,1+\epsilon])$ and also in the sense of
 pointwise convergence for 
$z\ne 0$, where

\begin{align}\label{defsolpm}
\bfQ_0^\pm(z)&:=\left\{\begin{array}{cc}
\bfQ_0(z)
\phantom{\pm0i)^{\beta'}\bfa_1(p)+\bfa_2(p)}&\ \ \ \ {z>0}\cr
|z|^{\beta'}\mathrm{e}^{\mp \mathrm{i}\pi(\beta')}\bfa_1(p)+\bfa_2(p)\
\  \ \ \ \ &{z<0}
\cr
\end{array}\ \ (\beta\ne 1)
\right.\cr
\bfQ_0^\pm&:=\left\{\begin{array}{cc}
\bfQ_0(z)
\phantom{(\ln(|z|)\mp\pi i)\bfa_1(z)+\bfa_2(z)}&{z>0}\cr
z(\ln(|z|)-\pi i)\bfa_1(z)+z\bfa_2(z)&{z<0}
\cr
\end{array}\ \ (\beta'= 1)
\right.
\end{align}

\z Moreover, $\bfQ_0^{\pm}$ are $\lloc$ 
solutions of (\ref{eqinte})
 on the interval $[-\epsilon,\epsilon]$.
\end{Remark}

The proof is immediate from Propositions
\ref{P8} and \ref{eqinte}.

\Box

\begin{Remark}\label{lincomb}
For any $\lambda\in\CC$ the combination
$\bfQ_\lambda=\lambda\bfQ_0^++(1-\lambda)\bfQ_0^-$
is a solution of (\ref{difeqbfQ2}) in $\lone[-\epsilon,\epsilon]$.

\end{Remark}

{\em Proof. } Follows from Remark~\ref{uniformlb} as  (\ref{difeqbfQ2}) is linear.

\Box

\centerline{*}

\z Let now $\bfQ_0$ be any solution of (\ref{difeqbfQ2}) in
$\lone[-\epsilon,\epsilon]$. We search for other solutions in
the form $\bfQ=\bfQ_0+\bfq$. Since (\ref{difeqbfQ2}) is linear and
$\bfQ_0$ is already a solution we have

\begin{eqnarray}\label{difeqbfq2}
&&\zpp\bfq'(z)-(\beta'\oplus\bcp) \bfq(z)=\sum_{j=1}^n\int_\epsilon^z\bfd_j'(s-z)q_j(s)\mathrm{d}s\cr&&
\end{eqnarray}

\z and, by the uniqueness of $\bfQ_0$ for $z>0$ we have $\bfq=0$
for $q<0$ and the equation becomes

\begin{eqnarray}\label{difeqbfq3}
&&\zpp\bfq'(z)-(\beta'\oplus\bcp) \bfq(z)=\sum_{j=1}^n\int_0^z\bfd_j'(s-z)q_j(s)\mathrm{d}s\cr&&
\end{eqnarray}

\z with the initial condition $\bfq(0)=0$. Changing variables to
$z=-p$ ($p>0$ now corresponds to going beyond the singularity) and
$\bfq(z)=\bfY(-p)$ we have

\begin{eqnarray}\label{difeqbfq3+}
&&\ppp\bfY'(p)-(\beta'\oplus\bcp) \bfY(p)+\sum_{j=1}^n\int_0^p\bfd_j'(p-t)Y_j(t)\mathrm{d}t=0\cr&&
\end{eqnarray}

\z We recognize in (\ref{difeqbfq3+}) the equation
for $\mathcal{P}\mathbf{Y}_{\mathbf{e}_1}$:

\begin{Remark}\label{coinci}
Equation (\ref{difeqbfq3+}) is at the same time
the equation for $\mathcal{P}\mathbf{Y}_{\mathbf{e}_1}$ and for
the difference $\mathcal{P}^{m+1}(\mathbf{Y}_0^{[1]}-\mathbf{Y}_0^{[2]})$
where $\mathbf{Y}_0^{[1,2]}$ are any solutions of (\ref{eqil}).
\end{Remark}
\Box

\begin{Proposition}\label{localresult} Let $\epsilon$ be small.  In
$\mathcal{T}_{\beta'}(\{|p|<\epsilon\})$ as well as in
$\mathcal{D}'_{m,\nu}(0,\epsilon \mathrm{e}^{\mathrm{i}\phi})$ for any $\phi$, there is
a unique solution of (\ref{difeqbfq3+}) $\mathbf{W}_0$ such that, for
small $p$, $\mathbf{W}_0=
\Gamma(\beta')^{-1}p^{\beta'}(\mathbf{e}_1+o(1))$. The general solution of
(\ref{difeqbfq3+}) is $\bfY=C\bfW_0$, with $C\in\CC$ arbitrary.

\end{Proposition}

Notes: (1) Modulo relabeling of the spatial directions, 
the statement and proof hold for any of the equations for 
$\bfY_{\mathbf{e}_j}$. 
 
(2) The point $p=0$ is singular, and so
is the ``initial condition'' $\mathbf{W}_0\sim
\Gamma(\beta')^{-1}p^{\beta'}\mathbf{e}_1$.

{\em Proof.} We have

\begin{eqnarray}\label{syt12}
&&(p^{-\beta'} Y_1(z))'=-p^{-\beta'-1}
\sum_{j=1}^n\int_0^p D_{1j}'(p-t)Y_j(t)\mathrm{d}t\cr
&&(\mathrm{e}^{\hat{E}(p)}\bfY_\perp)'=-\mathrm{e}^{\hat{E}(p)
}\gc(-p)^{-1}\sum_{j=1}^n\int_0^p
\dpp'(p-t)Y_j(t)\mathrm{d}t\cr
&&\hat{E}(p):=-\int_0^p \gc(-t)^{-1}\bcp(t)\mathrm{d}s
\cr
&&\bfY(0)=0
\end{eqnarray}

\z After integration we get:

\begin{eqnarray}\label{eqperp122}
(1+J_1)\bfY_1(z)=C R_1(p)\phantom{0}\cr
(1+J_\perp)\bfY_\perp(p)=0\phantom{R_1(z)}
\end{eqnarray}

\z with $C\in\CC$ arbitrary and
\begin{eqnarray}\label{eqperp22}
&&J_1\,{\bf V}=p^{\beta'}\int_0^pt^{-\beta'-1}\sum_{j=1}^n\int_0^tY_j(s)D_{1j}'(t-s)\mathrm{d}s\mathrm{d}t\cr
&& 
J_\perp\,{\bf V}:=\mathrm{e}^{-\hat{E}(p)}
\int_0^p \mathrm{e}^{\hat{E}(t)}\gc(-t)^{-1}
\left(\sum_{j=1}^n\int_0^t
\dpp'(t-s)Y_j(s)\mathrm{d}s\right)\mathrm{d}t\cr
&&R_1(p)=p^{\beta'}\cr
\cr
&&
\end{eqnarray}

In a region $|p|<\epsilon$, for small $\epsilon$,
the norm of the operator $J$ defined 
on $\mathcal{T}_{\beta'}$ is $O(\epsilon)$, exactly as in
Proposition~\ref{eqinte}. Given
$C$ the solution of the system (\ref{eqperp122})
is unique and can be written as 

\begin{equation}\label{eqneuman}
\bfY=C \bfW_0;\ \ \bfW_0:=\Gamma(\beta')^{-1}(\hat 1+J)^{-1}{\bf R}\ne 0
\end{equation}

\z (The prefactor $\Gamma(\beta')^{-1}$ was introduced so that the
coefficient of the leading power in the asymptotic series of
$\lap\bfW_0$ is one).

The proof is essentially the same if we consider (\ref{eqperp122}) in
$L^1_\nu(0,\epsilon \mathrm{e}^{\mathrm{i}\phi}),$ which coincides with
$\mathcal{D}'_{m,\nu}(0,\epsilon \mathrm{e}^{\mathrm{i}\phi})$ for small $\epsilon$.

\begin{Remark}\label{identif}
On $(1,1+\epsilon)$, $(\bfY_0^+-\bfY_0^-)(p)=
S_1\mathbf{Y}_{\mathbf{e}_1}^{(m)}(p-1)$.
\end{Remark}

The existence of some $S_1$  is obvious from Remark~\ref{coinci}
and  Proposition~\ref{localresult}. Its value follows
by comparing (\ref{Y0p=1}), (\ref{primaresu}) and (\ref{eqneuman}).

\end{subsubsection}

\begin{Proposition}\label{medianization} i) Let
  $\mathbf{Y}_\bfk^+,\bfk\ge 0$ solve (\ref{eqil}),
  (\ref{systemformv}) in $\mathcal{T}_{\{\cdot\}}({\mathcal{S}^+}')$.
  Then  $\mathbf{Y}_\bfk^{ba}$, cf.
  (\ref{defmed}) solve (\ref{eqil}), (\ref{systemformv}) in
  $\mathcal{D}'_{m,\nu} (\RR^+)$

ii) For any of the functions $\bfY_\bfk$, interchanging
$+$ with $-$ in (\ref{defmed}) does not change
the balanced average.

\end{Proposition}

{\em Proof} (i) The fact that $\mathbf{Y}_0^{ba}$ is a solution of
(\ref{eqil}) follows from Proposition~\ref{SRY0}. From
Proposition~\ref{SRY0} and Proposition~\ref{medianpropo} we see that
$\mathbf{D}_\bfj^{ba}$ is obtained
by simply replacing $\mathbf{Y}_0^{+}$ by $\mathbf{Y}_0^{ba}$
in (\ref{exprdj}) (notice that on any finite
interval, there are finitely many terms
in the expression of  $\mathbf{D}_\bfj^{ba}-\mathbf{D}_\bfj^{+}$.)
The rest of the proof merely consists in inductively
applying $\mathcal{A}_\alpha$ to the equations (\ref{systemformv}),
noting that each contains finitely many convolutions,
and applying the commutation relation (\ref{assertmed}).

(ii) This is true for $\mathbf{Y}_0$ as an immediate verification
shows that the $+$ and $-$ averages coincide on $(0,2)$ (where they
consist in two terms). Thus by Proposition~\ref{contidis} they have to
coincide on $\RR^+$. With this, for the rest of the $\mathbf{Y}_\bfk$ the
property follows by an obvious induction from
Proposition~\ref{medianpropo}.

\subsection{Appendix}
 
\subsubsection{The $C^*$--algebra $\mathcal{D}'_{m,\nu}$}

\label{starca}

\z Let $\mathcal{D}$ be the space of test functions (compactly
supported $C^\infty$ functions on $(0,\infty)$) and $\mathcal{D}(0,x)$
be the test functions on $(0,x)$.

We say that $f\in\mathcal{D}'$ is a staircase distribution if for any $k=0,1,2,...$ there is an $\lone$ function on
$[0,k+1]$ so that $f=F_k^{(km)}$ (in the sense of distributions) when
restricted to $\mathcal{D}(0,k+1)$ or

\begin{eqnarray}
  \label{imposecond0}
F_k:=\mathcal{P}^{mk}f \in L_1(0,k+1)
\end{eqnarray}

\z (since ${f}\in\lloc[0,1-\epsilon]$ and
by Remark~\ref{density}, $\mathcal{P}f$ is well defined).
With this choice we have

\begin{eqnarray}
  \label{imposecond}
F_{k+1}=\mathcal{P}^m{F_k}\ \mbox{on }[0,k]
\mbox{ and } F_k^{(j)}(0)=0 \ \mbox{for}\ j\le mk-1
\end{eqnarray}

We denote these distributions by $\mathcal{D}'_m$
($\mathcal{D}'_m(0,k)$ respectively, when restricted to
$\mathcal{D}(0,k)$) and observe that
$\bigcup_{m>0}\mathcal{D}'_m\supset S'$, the distributions of slow
growth. The inclusion is strict since any element of $S'$ is of finite
order.

Let $f\in\lone$. Taking $F=\mathcal{P}^j f\in C^j$ we have, by integration
by parts and noting that the boundary terms vanish,

\begin{eqnarray}
  \label{sampleconv}
  (F*F)(p)=\int_0^pF(s)F(p-s)\mathrm{d}s=\int_0^p F^{(j)}(s)\mathcal{P}^jF(p-s)
\end{eqnarray}

\z so that $F*F\in C^{2j}$ and

\begin{eqnarray}
  \label{sample2}
 (F*F)^{(2j)}=f*f
\end{eqnarray}

\z This motivates the following definition:
for $f,\tilde{f}\in\mathcal{D}'_m$ let
\begin{eqnarray}
  \label{defconvd}
  f*\tilde{f}:=(F_k*\tilde{F}_k)^{(2km)}\ \ \mbox{in } \mathcal{D}'(0,k+1)
\end{eqnarray}
\z
We first check that the definition is consistent in the sense that
$$(F_{k+1}*F_{k+1})^{(2m(k+1))}=(F_{k}*F_{k})^{(2mk)}$$ on
$\mathcal{D}(0,k+1)$. For $p<k+1$ integrating by parts and using
(\ref{imposecond}) we obtain

\begin{eqnarray}
  \label{consisit1}
  &&\frac{\mathrm{d}^{2m(k+1)}}{\mathrm{d}p^{2m(k+1)}}\int_0^p F_{k}(s)\mathcal{P}^{2m}
  \tilde{F}_{k}(p-s)\mathrm{d}s=
\frac{\mathrm{d}^{2mk}}{\mathrm{d}p^{2mk}}\int_0^p F_{k}(s)
  \tilde{F}_{k}(p-s)\mathrm{d}s\cr&&
\end{eqnarray}

\z The same argument shows that the definition
is compatible with the embedding of $\mathcal{D}'_{m}$
in $\mathcal{D}'_{m'}$ with $m'>m$.  Convolution is commutative and
associative: with $f,g,h\in\mathcal{D}'_m$ 
and identifying 
 $(f*g)$ and $h$ by the natural inclusion
with elements in $\mathcal{D}'_{2m}$  we obtain
$(f*g)*h=((F*G)*H)^{(4mk)}=f*(g*h)$. 

 The following staircase decomposition
exists in $\mathcal{D}'_m$.

\begin{Lemma}\label{Stcase}. For each $f\in\mathcal{D}'_m$ there is a {\em unique}
sequence 
$\left\{\Delta_i\right\}_{i=0,1,..}$ such that $\Delta_i\in\lone(\RR^+)$,
$\Delta_i=\Delta_i\bchi_{[i,i+1]}$ and 

\begin{eqnarray}
  \label{stdec}
  f=\sum_{i=0}^{\infty}\Delta_i^{(mi)}
\end{eqnarray}

Also (cf. (\ref{imposecond})),

\begin{eqnarray}
  \label{decompik}
  F_i=\sum_{j\le i}\mathcal{P}^{m(i-j)}\Delta_i \ \mbox{on } [0,i+1)
\end{eqnarray}

\end{Lemma}

\z Note that the infinite sum is $\mathcal{D}'-$convergent since
for a given test function only a finite number of
distributions are nonzero.

{\em{Proof}}

\z We start by showing (\ref{decompik}).  For $i=0$ we take
$\Delta_{0}=F_0\bchi[0,1]$ (where $F_0\bchi[0,1]$ $:=$ $ 
\phi\mapsto\int_0^1F_0(s)\phi(s)\mathrm{d}s$).  Assuming
(\ref{decompik}) holds for $i<n$ we simply note that

\begin{multline}
  \label{defDelI}
  \Delta_{n}:=\bchi_{[0,n+1]}\left(F_n-\sum_{j\le
    n-1}\mathcal{P}^{m(n-j)}\Delta_{j}\right)\cr
 =\bchi_{[0,n+1]}
  \left(F_n-\mathcal{P}^m(F_{n-1}\bchi_{[0,n]})\right)=
\bchi_{[n,n+1]}
  \left(F_n-\mathcal{P}^m(F_{n-1}\bchi_{[0,n]})\right)
\end{multline}

\z (with $\bchi_{[n,\infty]}F_n$ defined in the same way as
$F_0\bchi[0,1]$ above) has, by the induction hypothesis and
(\ref{imposecond}) the required properties.  Relation (\ref{stdec}) is
immediate.  It remains to show uniqueness.  Assuming (\ref{stdec})
holds for the sequences $\Delta_i,\tilde{\Delta}_i$ and restricting
$f$ to $\mathcal{D}(0,1)$ we see that $\Delta_0=\tilde{\Delta}_0$.
Assuming $\Delta_i=\tilde{\Delta}_i$ for $i<n$ we then have
$\Delta_n^{(mn)} =\tilde{\Delta}_n^{(mn)}$ on $\mathcal{D}(0,n+1)$. It
follows from Remark~\ref{density} that $\Delta_n(x)
=\tilde{\Delta}_n(x)+P(x)$ on $[0,n+1)$ where $P$ is a polynomial (of
degree $<mn$).  Since by definition
$\Delta_n(x)=\tilde{\Delta}_n(x)=0$ for $x<n$ we have
$\Delta_n=\tilde{\Delta}_n(x)$. \Box

The expression (\ref{defconvd}) hints to  decrease in regularity, but this
is not the case. In fact, we
check that the regularity of convolution is not worse than that of
its arguments.

\begin{Remark}\label{weldefi}
\begin{eqnarray}\label{welldef}
(\cdot\,*\,\cdot):\mathcal{D}_n\mapsto\mathcal{D}_n
\end{eqnarray}

\end{Remark}

\z  Since 

\begin{eqnarray}
  \label{identi0}
\bchi_{[a,b]}*\bchi_{[a',b']}=\left(\bchi_{[a,b]}*\bchi_{[a',b']}\right)
\bchi_{[a+a',b+b']}
\end{eqnarray}

\z we have

\begin{eqnarray}
  \label{verifreg}
  &&F*\tilde{F}=\sum_{j+k\le \lfloor p\rfloor}\mathcal{P}^{m(i-j)}
\Delta_{j}*\mathcal{P}^{m(i-k)}\tilde{\Delta}_{k}
=\sum_{j+k\le \lfloor p\rfloor}
\Delta_{j}*\mathcal{P}^{m(2i-j-k)}\tilde{\Delta}_{k}\cr&&
\end{eqnarray}

\z which is manifestly in $C^{2mi-m(j+k)}[0,p)\subset C^{2mi-m\lfloor
  p\rfloor}[0,p)$.

\Box 

\subsubsection{Norms on $\mathcal{D}'_m$}

\z For $f\in\mathcal{D}'_m$   define 

\begin{eqnarray}
  \label{nord1}
  \|f\|_{\nu ,m}:=c_m\sum_{i=0}^{\infty}\nu ^{im}\|\Delta_i\|_{L^1_\nu }
\end{eqnarray}

\z (the constant $c_m$, immaterial for the moment, is defined in
(\ref{newcm}). When no confusion is possible we will simply write
$\|f\|_\nu $ for $\|f\|_{\nu ,m}$ and $\|\Delta\|_\nu $ for
$\|\Delta_i\|_{L^1_\nu }$ (no other norm is used for the $\Delta$'s).
Let $\mathcal{D'}_{m,\nu }$ be the distributions in $\mathcal{D}'_m$
such that $\|f||_\nu <\infty$.

\begin{Remark}\label{normisnorm}
$\|\cdot\|_\nu $ is a norm on $\mathcal{D'}_{m,\nu }$. 

\end{Remark}

\z If $\|f\|_\nu =0$ for all $i$, then $\Delta_i=0$
whence $f=0$. In view of Lemma~\ref{Stcase} we have 
$\|0\|_\nu =0$. All the other properties are immediate.

\begin{Remark}\label{completen} $\mathcal{D'}_{m,\nu }$ is a Banach space.
The topology
given by $\|\cdot\|_\nu $ on $\mathcal{D'}_{m,\nu }$ is stronger
than the topology inherited from $\mathcal{D}'$. 

\end{Remark}

{\em Proof.}  If we let $\mathcal{D}'_{m,\nu}(k,k+1)$ be the subset
of $\mathcal{D'}_{m,\nu }$ where  all $\Delta_i=0$ except
for $i=k$, with the norm (\ref{nord1}), we have 

\begin{eqnarray}
  \label{directsum}
  \mathcal{D'}_{m,\nu }=\bigoplus_{k=0}^{\infty}\mathcal{D}'_{m,\nu}(k,k+1)
\end{eqnarray}

\z and we only need to check completeness of each
$\mathcal{D}'_{m,\nu}(k,k+1)$ which is immediate: on $\lone[k,k+1]$,
$\|\cdot\|_\nu $ is equivalent to the usual $\lone$ norm and thus if
 $f_n\in \mathcal{D}'_{m,\nu}(k,k+1)$ is a Cauchy sequence
then
$\Delta_{k,n}\stackrel{L_\nu}{\rightarrow}\Delta_k$ (whence weak
convergence) and 
$f_n\stackrel{\mathcal{D}'_{m,\nu}(k,k+1)}{\rightarrow} f$
where $f=\Delta_k^{(mk)}$. \Box

\begin{Lemma}\label{C*}
The space $\mathcal{D'}_{m,\nu}$ is a $C^*$ algebra with
respect to convolution.

\end{Lemma}

{\em Proof.}  Let
$f,\tilde{f}\in\mathcal{D}'_{m,\nu}$ with
$$f=\sum_{i=0}^\infty \Delta_i^{(mi)}\ \ ,\ \
\tilde{f}=\sum_{i=0}^\infty \tilde{\Delta}
_i^{(mi)}$$
Then 
\begin{equation}\label{desf}
f* \tilde{f}=\sum_{i,j=0}^\infty\Delta_i^{(mi)}*\tilde{\Delta}
_j^{(mj)}=\sum_{i,j=0}^\infty\left( \Delta_i*\tilde{\Delta} _j
\right)^{m(i+j)} 
\end{equation}
and the support of $ \Delta_i*\tilde{\Delta} _j$ is in
$[i+j,i+j+2]$ i.e.  $\Delta_i*\tilde{\Delta}
_j=\bchi_{[i+j,i+j+2]}\Delta_i*\tilde{\Delta} _j$.

We  first evaluate the norm in $\mathcal{D}'_{m,\nu}$ of the terms
$\left( \Delta_i*\tilde{\Delta} _j\right)^{m(i+j)} $.

\z {\bf{I. Decomposition formula.}} Let $f=F^{(mk)}\in\mathcal{D}'(\RR_+)$, where
$F\in L^1(\RR_+)$, and $F$ is supported in $[k,k+2]$ i.e.,
$F=\bchi_{[k,k+2]}F $ ($k\geq 0$).  Then $f\in\mathcal{D}'_m$ and the
decomposition of $f$ (cf. (\ref{stdec})) has the terms:
\begin{equation}\label{zerodelta}
\Delta_0=\Delta_1=...=\Delta_{k-1}=0\ \ ,\ \
\Delta_k=\bchi_{[k,k+1]}F
\end{equation}
and 

\begin{equation}\label{decdelta}
\Delta_{k+n}=\bchi_{[k+n,k+n+1]}G_n,\ \mbox{ where }G_n=\mathcal{P}^m\left( \bchi_{[k+n,\infty)}G_{n-1}\right),\ \ G_0=F
\end{equation}

\z {\em{Proof of Decomposition Formula}}. We use first line of (2.98)
of the paper
\begin{equation}\label{deltaform}
\Delta_j=\bchi_{[j,j+1]}\left
( F_j-\sum_{i=0}^{j-1}\mathcal{P}^{m(j-i)}\Delta_i\right)
\end{equation}
where, in our case, $F_k=F,\ 
F_{k+1}=\mathcal{P}^m F,\,  ...,\, F_{k+n}=\mathcal{P}^{mn}F,\, ...$.

The relations (\ref{zerodelta}) follow directly from
(\ref{deltaform}). Formula (\ref{decdelta}) is shown  by induction on $n$.  
For $n=1$ we have
$$\Delta_{k+1}=\bchi_{[k+1,k+2]}\left( \mathcal{P}^m
\, F-\mathcal{P}^m\Delta_{k} \right)$$

$$=\bchi_{[k+1,k+2]}\mathcal{P}^m\left(
  \bchi_{[k,\infty)}F-\bchi_{[k,k+1]}F
\right)=\bchi_{[k+1,k+2]}\mathcal{P}^m\left( \bchi_{[k+1,\infty)}F
\right)$$
Assume (\ref{decdelta}) holds for $\Delta_{k+j}$, $j\le
n-1$. Using (\ref{deltaform}), with $\bchi=\bchi_{[k+n,k+n+1]}$ we
have
$$\Delta_{k+n}=\bchi\left(\mathcal{P}^{mn}F-\sum_{i=k}^{n-1}\mathcal{P}^{m(n-i)}\Delta_i\right)=\bchi
\mathcal{P}^m\left( G_{n-1}-\Delta_{n-1}\right)$$

$$=\bchi\mathcal{P}^m\left(
  \bchi_{[k+n-1,\infty)}G_{n-1}-\bchi_{[k+n-1,k+n]}G_{n-1}\right)=\bchi\mathcal{P}^m\left( \bchi_{[k+n,\infty)}G_{n-1}\right)\Box$$

\z {\bf{II. Estimating $\Delta_{k+n}$.}} For $f$ as in {\bf{I}}, we have 
\begin{equation}\label{normkp1}
||\Delta_{k+1}||_\nu\leq \nu^{-m}|| F||_\nu\ \ ,\ \ 
||\Delta_{k+2}||_\nu\leq \nu^{-2m}|| F||_\nu
\end{equation}
and, for $n\geq 3$
\begin{equation}\label{normkpn}
||\Delta_{k+n}||_\nu\leq e^{2\nu-n\nu}(n-1)^{nm-1} \frac{1}{(nm-1)!}|| F||_\nu
\end{equation}

\z {\em{Proof of estimates of $\Delta_{k+n}$}}. 

\z (A) Case $n=1$. 

\begin{multline}
  \label{fest}
  ||\Delta_{k+1}||_\nu\leq \int_{k+1}^{k+2} \, dt\, e^{-\nu t} \mathcal{P}^m\left( \bchi_{[k+1,\infty)}|F|
\right)(t)\\
= \int_{k+1}^{k+2} \, dt\, e^{-\nu t} \int_{k+1}^t\, ds_1\,
\int_{k+1}^{s_1}\, ds_2\, ...\int_{k+1}^{s_{m-1}}\, ds_m |F(s_m)|\\
\le  \int_{k+1}^{k+2}\, ds_m |F(s_m)|\, \int_{s_m}^{\infty} \,
ds_{m-1}\, ... \int_{s_2}^\infty \, ds_1\, \int_{s_1}^{\infty} \, dt\,
e^{-\nu t} \\
 = \int_{k+1}^{k+2}\, ds_m |F(s_m)|e^{-\nu s_m}\nu^{-m}
\leq \nu^{-m}|| F||_\nu
\end{multline}

\z (B) Case $n=2$: 

$$||\Delta_{k+1}||_\nu\leq \int_{k+2}^{k+3} \, dt\, e^{-\nu t}\mathcal{P}^m
\left( \bchi_{[k+2,\infty)} \mathcal{P}^m\left( \bchi_{[k+1,\infty)}|F|
\right)\right)$$

$$= \int_{k+2}^{k+3} \, dt\, e^{-\nu t} \int_{k+2}^t\, dt_1\,
\int_{k+2}^{t_1}\, dt_2\, ...\int_{k+2}^{t_{m-1}}\, dt_m
\int_{k+1}^{t_m}\, ds_1\,
\int_{k+1}^{s_1}\, ds_2\, ...\int_{k+1}^{s_{m-1}}\, ds_m |F(s_m)|$$

$$\le \int_{k+2}^{k+3}\, ds_m |F(s_m)|\,  \int_{s_m}^{\infty}\, ds_{m-1}\,
...\int_{s_2}^{\infty}\, ds_1 \, \int_{\max\{s_1,k+2\}} ^{\infty}\,
dt_m\, \int_{t_m}^\infty \, dt_{m-1}\, ... \int_{t_1}^\infty \, dt\,
e^{-\nu t}  $$

$$= \int_{k+2}^{k+3}\, ds_m |F(s_m)|\,  \int_{s_m}^{\infty}\, ds_{m-1}\,
...\int_{s_2}^{\infty}\, ds_1
e^{-\nu \max\{s_1,k+2\}}\nu^{-m-1}$$

$$ \le\int_{k+2}^{k+3}\, ds_m |F(s_m)|\,
\int_{s_m}^{\infty}\, ds_{m-1}\,...\int_{s_3}^\infty \, ds_2 e^{-\nu
  s_2}\nu^{-m-2}=\int_{k+2}^{k+3}\, ds_m |F(s_m)|e^{-\nu s_m}\nu^{-2m}$$ 

\z (C) Case $n\geq 3$. We first estimate $G_2,...,G_n$:

$$|G_2(t)|\leq \mathcal{P}^m\left( \bchi_{[k+2,\infty)}
  \mathcal{P}^m\left( \bchi_{[k+1,\infty)}|F| \right)\right) (t)$$

$$= \int_{k+2}^t\, dt_1\, \int_{k+2}^{t_1}\, dt_2\,
...\int_{k+2}^{t_{m-1}}\, dt_m \int_{k+1}^{t_m}\, ds_1\,
\int_{k+1}^{s_1}\, ds_2\, ...\int_{k+1}^{s_{m-1}}\, ds_m |F(s_m)|$$
and using the inequality
$$|F(s_m)|=|F(s_m)|\bchi_{[k,k+2]}(s_m)\le |F(s_m)|e^{-\nu
  s_m}e^{\nu(k+2)}$$
we get
$$|G_2(t)|\le e^{\nu(k+2)} ||F||_\nu \int_{k+1}^t\, dt_1\,
\int_{k+1}^{t_1}\, dt_2\, ...\int_{k+1}^{t_{m-1}}\, dt_m
\int_{k+1}^{t_m}\, ds_1\,
\int_{k+1}^{s_1}\, ds_2\, ...\int_{k+1}^{s_{m-2}}\,
ds_{m-1}$$

$$=e^{\nu(k+2)} ||F||_\nu  (t-k-1)^{2m-1}\frac{1}{(2m-1)!}$$
  
\z The estimate of $G_2$ is used for bounding $G_3$:

$$|G_3(t)|\leq \mathcal{P}^m\left( \bchi_{[k+3,\infty)}
  |G_2|\right)\le \mathcal{P}^m\left( \bchi_{[k+1,\infty)}
  |G_2|\right)\ \ \ \ \ \ \ \ \ \ \ \ \ \ \ \ \ \ \ \ $$

$$\ \ \ \ \ \ \ \ \  \ \  \ \ \ \ \ \ \ \ \ \ \    \ \ \le e^{\nu(k+2)} ||F||_\nu (t-k-1)^{3m-1}\frac{1}{(3m-1)!}$$
and similarly (by induction)
$$|G_n(t)|\le e^{\nu(k+2)} ||F||_\nu  (t-k-1)^{nm-1}\frac{1}{(nm-1)!}$$

\z Then 
$$||\Delta_{k+n}||_\nu\leq  e^{\nu(k+2)} ||F||_\nu \frac{1}{(nm-1)!}
\int_{k+n}^{k+n+1}\, dt\, e^{-\nu t}(t-k-1)^{nm-1} $$
and, for $\nu\geq m$ the integrand is decreasing, and the inequality
(\ref{normkpn}) follows.

\z {\bf{III. Final Estimate}}. Let $\nu_0>m$ be fixed.  For $f$ as in
{\bf{I}}, we have for any $\nu >\nu_0$,
\begin{equation}\label{Fk}
||f||\le c_m \nu^{km}||F||_\nu 
\end{equation}

\z for some $c_m$, if $\nu>\nu_0>m$.

\z {\em{Proof of Final Estimate}}

$$||f||=\sum_{n\geq 0}\nu^{km+kn}||\Delta_{k+n}||_\nu \le
\nu^{km}||F||_\nu \left[ 3+\sum_{n\geq 3}
  \nu^{nm}e^{2\nu-n\nu}\frac{(n-1)^{nm-1}}{(nm-1)!}\right]$$
and,
using $n-1\leq (mn-1)/m$ and a crude Stirling estimate we obtain

\begin{equation}\label{newcm}||f||\le \nu^{km}||F||_\nu \left[ 3+ me^{2\nu-1}\sum_{n\geq 3} \left(
    e^{m-\nu}\nu^m/m^m\right)^{n} \right]\le c_m \nu^{km}||F||_\nu 
\end{equation}
\z   Thus  (\ref{Fk}) is proven for
$\nu>\nu_0>m$.

\z {\bf{End of the proof}}. From (\ref{desf}) and (\ref{Fk}) we get
$$||f*\tilde{f}||\leq \sum_{i,j=0}^\infty||\left(
  \Delta_i*\tilde{\Delta} _j \right)^{m(i+j)}||$$

$$\leq \sum_{i,j=0}^\infty c_m^2 \nu^{m(i+j)}
||\Delta_i*\tilde{\Delta} _j ||_\nu \leq c_m^2 \sum_{i,j=0}^\infty
\nu^{m(i+j)} ||\Delta_i||_\nu \, ||\tilde{\Delta} _j ||_\nu=c_m^2
||f||\, ||\tilde{f}||$$ $\Box$

\begin{Remark}\label{conveB} Let $f\in\mathcal{D}'_{m,\nu}$ for some $\nu >\nu _0$ where
$\nu _0^m=\mathrm{e}^{\nu _0}$. Then
$f\in\mathcal{D}'_{m,\nu'}$ for all $\nu '>\nu $ and furthermore,

\begin{eqnarray}
  \label{tendzero}
  \|f\|_\nu \downarrow 0\ \mbox{as } \nu \uparrow\infty
\end{eqnarray}

\end{Remark}

{\em Proof.} We have

\begin{eqnarray}
  \label{difB}
  \nu ^{mk}\int_k^{k+1}|\Delta_k(s)|\mathrm{e}^{-\nu s}\mathrm{d}s=
(\nu ^m\mathrm{e}^{-\nu })^k\int_0^1 |\Delta_k(s+k)|\mathrm{e}^{-\nu s}\mathrm{d}s
\end{eqnarray}

\z which is decreasing in $\nu $. The rest follows from
the monotone convergence theorem.
\Box


\subsubsection{Embedding of  $L^1_\nu$ in $\mathcal{D}'_m$}

\begin{Lemma}\label{imbeddi}
i) Let  $f\in L^1_{\nu_0}$ (cf. Remark~\ref{conveB}). 
Then $f\in\mathcal{D}'_{m,\nu}$ for all $\nu>\nu_0$.

ii)  $\mathcal{D}(\RR^+\backslash\NN)\cap L^1_\nu(\RR^+) $ is dense in $\mathcal{D}_{m,\nu }$ with
respect to the norm $\|\|_\nu $.

\end{Lemma}

{\em Proof. }

 Note that if for some $\nu_0$ we have $f\in L^1_{\nu_0}(\RR^+)$
then

\begin{eqnarray}
  \label{unifversusint}
  \int_0^p|f(s)|\mathrm{d}s\le \mathrm{e}^{\nu_0 p}\int_0^p|f(s)|\mathrm{e}^{-{\nu_0 s}}\mathrm{d}s
\le \mathrm{e}^{\nu_0 p}\|f\|_{\nu_0}
\end{eqnarray}

\z to which, application of $\mathcal{P}^{k-1}$ yields

\begin{eqnarray}
  \label{genn}
  \mathcal{P}^k|f|\le \nu_0^{-k+1}\mathrm{e}^{\nu_0 p} \|f\|_{\nu_0}
\end{eqnarray}

\z Also, $\mathcal{P}\bchi_{[n,\infty)}\mathrm{e}^{\nu_0 p}\le
{\nu_0}^{-1}\bchi_{[n,\infty)}\mathrm{e}^{\nu_0 p}$ so that

\begin{eqnarray}
  \label{estimexp5}
  \mathcal{P}^m\bchi_{[n,\infty)}\mathrm{e}^{\nu_0 p}\le
{\nu_0}^{-m}\bchi_{[n,\infty)}\mathrm{e}^{\nu_0 p}
\end{eqnarray}

\z so that, by (\ref{defDelI}) (where now $F_n$ and
$\bchi_{[n,\infty]}F_n$ are in $\lloc(0,n+1)$) we have for
$n>1$,

\begin{eqnarray}
  \label{festDel}
  |\Delta_n|\le \|f\|_{\nu_0}\mathrm{e}^{\nu_0 p}{\nu_0}^{1-mn}\bchi_{[n,n+1]}
\end{eqnarray}

\z Let now $\nu$ be large enough.
We have

\begin{multline}
  \label{sumn}
  \sum_{n=2}^{\infty}\nu^{mn}\int_0^{\infty}|\Delta_n|\mathrm{e}^{-\nu p}\mathrm{d}p\le
\nu_0\|f\|_{\nu_0}
\sum_{n=2}^{\infty}\int_n^{n+1}\mathrm{e}^{-(\nu-\nu_0)
  p}\left(\frac{\nu}{\nu_0}\right)^p\mathrm{d}p\cr=
\frac{\mathrm{e}^{-2(\nu-\nu_0-\ln(\nu/\nu_0))}}{\nu-\nu_0-\ln(\nu/\nu_0)}\nu_0\|f\|_{\nu_0}
\end{multline}

\z For $n=0$ we simply have $\|\Delta_0\|\le\|f\|$, while
for $n=1$ we write

\begin{eqnarray}
  \label{case1}
  \|\Delta_1\|_{\nu}\le \|1^{*(m-1)}*|f|\|_{\nu}\le
\nu^{m-1}\|f\|_{\nu}
\end{eqnarray}

Combining the estimates above, the proof of (i) is  complete.
To show (ii), let $f\in\mathcal{D}'_{m,\nu }$ and let 
$k_\epsilon$ be such that
$c_m\sum_{i=k_\epsilon}^{\infty}\nu ^{im}\|\Delta_i\|_\nu <\epsilon$. 
For each $i\le k_\epsilon$ we take a function $\delta_i$ in $\mathcal{D}(i,i+1)$
such that $\|\delta_i-\Delta_i\|_\nu <\epsilon 2^{-i}$. Then
$\|f-\sum_{i=0}^{k_\epsilon}\delta_i^{(mi)}\|_{m,\nu }<2\epsilon$. \Box

{\em Proof} of continuity of $f(p)\mapsto pf(p)$. If $f(p)=
\sum_{k=0}^{\infty}\Delta_k^{(mk)}$ then
$pf=\sum_{k=0}^{\infty}(p\Delta_k)^{(mk)}-
\sum_{k=0}^{\infty}mk\mathcal{P}(\Delta_k^{(mk)})$=
$\sum_{k=0}^{\infty}(p\Delta_k^{(mk)})-1*\sum_{k=0}^{\infty}
(mk\Delta_k)^{(mk)}$.  The rest is obvious from continuity of
convolution, the embedding shown above and the definition of the
norms.

\subsubsection{Laplace transforms}
\label{sec:LT}

{\em Proof of Lemma~\ref{existe} }. Let $\nu >\nu _0$ (cf.
Remark~\ref{conveB}). Equation (\ref{lapcomuta}) follows most easily
from the corresponding well property of Laplace transforms on
$\mathcal{D}$, from the continuity of $\lap$ and Lemma~\ref{imbeddi} (ii).
For the second notice that by the definition of $\mathcal{D}'_m,$
$f'\in\mathcal{D}'_m$ implies $f\in AC(0,1-\epsilon)$ and the property
follows by density from the $\mathcal{D}$ identity
$\mathcal{L}(\mathcal{P}g)= x^{-1}\lap(g)$. The third equality also
follows by density. The rest of the properties, except
injectivity, follow immediately from the definitions and the topology
used.

To show injectivity assume that $\mathcal{L}d(x)=0$ where
$d\in\mathcal{D}'_{m,\nu}$, $\nu<x_0<x\in\RR^+$.  By analyticity,
$\mathcal{L}d(x)=0$ in (say) $S_2:=\{z:|z|>2x_0:|\arg(z)|<\pi/4\}$.
Using dominated convergence, assuming $x_0$ is large enough, we have
$$\left|\sum_{k=1}^{\infty}
x^{m(k-1)}\mathrm{e}^{-(k-1)x}\int_0^{1}\mathrm{e}^{-sx}|\Delta_k(s+k)|\mathrm{d}s\right| \le 1 $$
in $S_2$. Thus $|f(x)|=$
 $\left|\int_0^1 \mathrm{e}^{px}\Delta_0(1-p)\mathrm{d}p\right|$
$\le |x|^m$ in $S_2$. The function $f^{(m)}(x)$ is entire, of
exponential order less than $1+\epsilon$ for any $\epsilon$ and, using
the previous inequality in Cauchy's formula we see that
$|f^{(m)}(x)|<\mbox{const.}$ in $(x_0,\infty)$. Since for
$\phi\in(\pi/2,3\pi/2)$ we obviously also have
$f^{(m)}(r\mathrm{e}^{\mathrm{i}\phi})\rightarrow 0$ as $r\rightarrow\infty$, an
elementary instance of the Phragm\'en-Lindel\"of principle
\cite{Holland} implies that $f^{(m)}$ is bounded in $\CC$, therefore
constant, so $f$ itself is a polynomial that decays in the left half
plane, thus $f=0$. Therefore $\int_0^1
\mathrm{e}^{-px}\Delta_0(p)\mathrm{d}p=\lap\Delta_0=0$ so that $\Delta_0=0$.
Inductively and in the same way, we see that $\Delta_k=0$, $k\in\NN$.

\label{sec:A0}

{\em Proof of Lemma~\ref{cuteCauchy}}.
Take first $r\notin\ZZ$. 
 Choose $a_1,a_2$
so that $0<a_1<a_2<a$ and  consider the closed contour
$C$ going along the upper cut from $\xi =0$ to $\xi =a_2$, continuing towards the lower cut
anticlockwise along the circle $C(a_2)$ of radius $a_2$ centered at
origin, and finally  coming from $\xi =a_2$ back to $\xi =0$
along the lower cut.  For $|\xi|<a_1$ we have, by the assumptions
of the lemma,

\begin{eqnarray}
  \label{intepathf}
  2\pi if(\xi )=\oint_C\frac{f(s)}{s-\xi }\mathrm{d}s=
\oint_{C(a_2)}\frac{f(s)}{s-\xi }\mathrm{d}s+
\int_{0}^{a_2}\frac{s^rA(s)}{s-\xi }\mathrm{d}s\cr&&
\end{eqnarray}

\z On the other hand, defining $z^rA(z)$ in the interior
of $C(a)$ cut along the positive axis (with the usual
convention $\arg(z)=0$ on the upper cut), we have, for
the same contour as above and $\xi \in\mathcal{V}_{a_1}$

\begin{eqnarray}
  \label{intepathfA}
2\pi\,\mathrm{i}\xi ^rA(\xi )=
\oint_{C(a_2)}\frac{f(s)}{s-\xi }\mathrm{d}s+
\left(1-\mathrm{e}^{2\pi i r}\right)\int_{0}^{a_2}\frac{s^rA(s)}{s-\xi }\mathrm{d}s\cr&&
\end{eqnarray}

\z Comparing (\ref{intepathf}) to (\ref{intepathfA}) we get:

\begin{eqnarray}
  \label{resultf}
  &&f(\xi )=\frac{1}{1-\mathrm{e}^{2\pi i r}} \xi^rA(\xi )\cr&&-\frac{1}
{2\pi i(1-\mathrm{e}^{2\pi i r})}\oint_{C(a_2)}\frac{A(s)}{s-\xi}\mathrm{d}s+
\frac{1}{2\pi i}\oint_{C(a_2)}\frac{f(s)}{s-\xi}\mathrm{d}s\cr&&
\end{eqnarray}

\z As integrals of analytic functions with respect to complex
absolutely continuous measures ($A(s)\mathrm{d}s$ and $f(s)\mathrm{d}s$), the last two
terms in (\ref{resultf}) are analytic in $\xi$ for $|\xi|<a_1$.  Since
$a_1$ can be chosen arbitrarily close to $a$, the case $r\notin\ZZ$ is
proven. For $r\in\ZZ$ the argument is essentially the same, in terms
of $A(\xi)\xi^r\ln\xi$ instead of $\xi^r A(\xi)$.  The proof
generalizes immediately to linear combinations
of $\xi^r A(\xi)$.  \Box

{\em Proof of Lemma~\ref{STR1}}. On the interval $(k,k+1)$ we have
$ f^+=\sum_{i=1}^k (f^-_i)^{(mi)} $ or 

\begin{eqnarray}
  \label{deco1}
 \mathcal{P}^{mk+1}f^+=\sum_{i=1}^k\mathcal{P}^{m(k-i)+1} f^-_i
\end{eqnarray}

\z Let $\epsilon$ be small and positive. Since $f(t\mathrm{e}^{\mathrm{i}\phi})$ and
$g_i^{-}(t\mathrm{e}^{\mathrm{i}\phi})$ converge as $\phi\rightarrow 0$
in $\mathcal{D}'_{m,\nu}$ we have
that $\mathcal{P}^{m(k-i)+1}f(t\mathrm{e}^{\mathrm{i}\phi})$ and
$\mathcal{P}^{m(k-i)+1}g_i^{-}(\mathrm{e}^{\mathrm{i}\phi}t)$ converge on
$[0,k+1-\epsilon]$ uniformly to $\mathcal{P}^{mk+1}f^+$ and
$\mathcal{P}^{m(k-i)+1}f_i^{-}$ respectively.  The left side of
(\ref{deco1}) is the limit on $I=[k+\epsilon,k+1-\epsilon]$ of a
function analytic in a neighborhood in the {\em upper} half plane of
$I$ and continuous on $I$ while the right side is the limit of a
function analytic in a neighborhood in the {\em lower} half plane of
$I$ and continuous on $I$. The equality of their continuous limits on
$I$ implies in particular that $\mathcal{P}^{mk+1}f(t\mathrm{e}^{\mathrm{i}\phi})$
extends analytically through $I$ in the lower half plane, and its
continuation is analytic where $\sum_{i=1}^k\mathcal{P}^{m(k-i)+1}
g_i^{-}(t\mathrm{e}^{\mathrm{i}\phi})$ is. A corresponding
statement is true for the upper plane continuation of
$\mathcal{P}^{mk+1}f(t\mathrm{e}^{-\mathrm{i}\phi})$ and (i) follows.  Part (ii) now
follows also, as an immediate application of Lemma~\ref{cuteCauchy}.

\Box

{\em Proof of Proposition~\ref{medianpropo}}.

The fact that multiplication by a bounded analytic function
is well defined on $\mathcal{F}(\mathcal{R}'_1)$ is immediate. 
Since 

\begin{gather}\label{polar}
2f*g = (f+g)*(f+g)-f*f-g*g
\end{gather}

\z we may take $f=g$.  With
$h=\mathcal{P}^{mk+1}f\in\mathcal{F}(\mathcal{R}'_1)$ it suffices to
show for every $k$ that $h*h$ (defined near zero by (\ref{defconv})
and which equals $\mathcal{P}^{2mk+2}(f*f)$ there) extends
analytically to $\mathcal{R}'_1$ for $\Re(x)<k$. Since $f$ 
is analytic in $\mathcal{R}'_1$ so is $h$. 
 In particular $h$ can
be analytically continued along any ray $d\subset\mathcal{R}'_1$ other
than the real line, and we have, by analyticity and with $*_d$ meaning
convolution along $d$,

\begin{equation}\label{forac}
AC(h*h)=AC(h)*_d AC (h)
\end{equation}

\z Also, by (\ref{defindecom2})

\begin{eqnarray}
  \label{eqhh}
  h^-(p)=h^+(p)+\sum_{j=1}^{\infty}(h_j^+(p-j))^{(mj)}
\end{eqnarray}

Let $H_0=h^+$ and $H_j(p)=(h_j^+(p))^{(mj)}$. By construction
$H_j'$ have $\lone$ boundary values on $[0,k-j+1)$ as
$\Re(z)>0,\Im(z)\downarrow 0$ and so $H_j$ extend {\em continuously}
to the strip $0<\Re(z)<k-j+1,\Im(z)\ge 0$.  We
have, by (\ref{defindecom2}) and continuity

\begin{eqnarray}
  \label{defindecom2H}
  h^{-}(z)=\sum_{i=0}^{j}H_i(z-j)
\end{eqnarray}

\z for $\Re(z)\in[0,j)$ and $\Im(z)\ge 0$ since $H_i(x)=0$ in the left
half plane, by definition.  For the same reason we have, with
$p'=p-i-j$ and $i+j\ge 1$,

\begin{eqnarray}
  \label{anahi}
  &&\int_0^p H_i(x-i)H_j(p-x-j)\mathrm{d}x\cr&&=\left\{\begin{array}{cc}
\displaystyle{
\int_0^{p'} H_i(x)H_j(p'-x)\mathrm{d}x=J_{i,j}(p')\  }&(\Re(p)'>i+j)\cr 
0&(\Re(p')<i+j)\end{array}\right.
\end{eqnarray}

\z As both $H_i$ and $H_j$ are analytic in an open strip $\mathcal{S}$
in the first quadrant and continuous on $[0,k+1-\epsilon]$ we see from
(\ref{anahi}) that $J_{ij}(p')$ are also analytic in $\mathcal{S}$ and
continuous on $[0,k+1-\epsilon]$.  For $p\in (0,l+1)$, $l\le k$ we
have

\begin{eqnarray}
  \label{firstequali}
\big(H^{-}*H^{-}\big)(p)=
\sum_{i=0}^{l}\sum_{j=0}^i
J_{ij}(p-i)\cr&&
\end{eqnarray}

\z Now, by (\ref{forac}) and using the continuity of $H$ and of
convolution, we note that the left side of (\ref{firstequali})
represents the continuous limit along $(l,l+1) $ of $(H*H)^-$, a
function analytic in a domain in the lower half plane while the right
side is the limit of a function analytic in the upper half plane
and $(l,l+1)$ is contained in the  common boundary. As in the
proof of Lemma~\ref{STR1} we conclude that $h*h$, thus
$f*f$ extend analytically
in $\mathcal{R}'_1$.

Going back to the definition of $H$ we get on $(0,l+1)$,
 
\begin{eqnarray}
  \label{firstequali2}
\big(f^{-}*f^{-}\big)(p)=\big(f^{+}*f^{+}\big)(p)
+\sum_{i=1}^{l}\sum_{j=0}^i
(f_j^+)^{(mj)}*(f_{i-j}^+)^{m(i-j)}(p-i)\cr&&
\end{eqnarray}

\z where $f_j*f_{i-j}=(H_j*H_{i-j})^{(2mk+2)}$ is the convolution in
$\mathcal{D}'_{m,\nu}$ and in our case gives a function analytic in
the open region $\mathcal{S}$. By comparing
with (\ref{defindecom2}) and (\ref{defindecom2H}) we get by induction
$(f*f)_j=\sum_{s=0}^j f_s*f_{j-s}$ or, using (\ref{polar}) 
we get (\ref{defindecom2}).

Since by assumption $f_s$ and $g_s$  belong to $\mathcal{D}'_{m,\nu}$
and the sum in (\ref{defindecom3}) only contains a finite number of terms, it
follows that all  analytic continuations  
of $(f*g)$ also belong to $\mathcal{D}'_{m,\nu}$. Furthermore,
it follows immediately that $K(f*g,\nu)\le 2K(f,\nu)K(g,\nu)$.

Only the last equality in (\ref{assertmed}) needs a proof; we have

\begin{multline}
  \label{verifymedC}
  (\mathcal{A}_\alpha(f))^{*2}=\left(\sum_{i=0}^{\infty} \alpha^i( f_i(p-i))^{(mi)}\right)^{*2}=
\sum_{k=0}^{\infty} \alpha^k\sum_{j=0}^k (f_j*f_{k-j})^{(mk)}(p-k)\cr
=\sum_{k=0}^{\infty} \alpha^k((f*f)_k)^{(mk)}=\mathcal{A}_\alpha(f*f)
\end{multline}

We have $\|\mathcal{A}_C(f)\|_{m,\nu}\le
\|f\|_{m,\nu}\sum_{j=0}^{\infty}C^j
K(f,\nu)^j=(1-KC)^{-1}\|f\|_{m,\nu}$ so that if
$\bfY\in\mathcal{F}_r$ then $\|\mathcal{A}_C(\bfY^{*\bfl})\|_{m,\nu}
=\|(\mathcal{A}_C\bfY)^{*\bfl}\|_{m,\nu}\le
(\|\bfY\|_{m,\nu}/(1-KC))^{|\bfl|}$ so that if $\nu$ is large
enough the sum involved in the expression of $M$
is uniformly convergent in $\mathcal{D}'_{m,\nu}(\RR^+)$
and (\ref{defindecom3}) follows.

\Box


\begin{Lemma}\label{cinftycase} (i) Let $k_0\ge 0$ and let $\lambda$ be such that 
  $\Re(\lambda)<\alpha_1<k_0$ and
  $\left|\Im(\lambda)\right|<\alpha_2$.  Alternatively, let $k_0\ge 0$
  and $\lambda$ be such that
  $0<\alpha_1<\left|\Im(\lambda)\right|<\alpha_2$.  There exists a
  constant $C(\alpha_1,\alpha_2)$ independent of $k_0,\nu$ and
  $\lambda$ so that

\begin{eqnarray}
  \label{normlemm}
  \|U\|_{\mathcal{D}'_{m,\nu}(k_0,\infty)\mapsto
    \mathcal{D}'_{m,\nu}(k_0,\infty)}\le 
C(\alpha_1,\alpha_2)(1+|\lambda|)^{-1}
\end{eqnarray}

(ii) In both cases in (i), 
$U$ is 
 strongly continuous in $\lambda$. 

\end{Lemma}

{\em Proof}

The impediments in the proof come on the one hand
from  having to estimate
quotients of the form $\int|\mathcal{P}^n U f|/\int|\mathcal{P}^n f|$
and on the other hand from the nonlocal character
of the action of $U$ in our space.

In view of Eq. (\ref{directsum}) it is enough
to find a $k-$independent upper bound
for the norms of the restrictions of $U$ to $\mathcal{D}'_{m,\nu}(k,k+1)$,
$U:\mathcal{D}'_{m,\nu}(k,k+1)\mapsto \mathcal{D}'_{m,\nu}$. We are interested
in $\Re(\lambda)<1$ in
the cases (a) $b>\left|\Im(\lambda)\right|>a>0$, (b) $\lambda<-a<0$  real,
(c) $\lambda\in\RR^+$ or $\lambda$ complex, $|\Im(\lambda)|<b$
but with $\mbox{supp}(f)\in (k_1,\infty)$ 
with $k_1>a>\Re(\lambda)$, $k_1\in(a,a+1)$. We let $k_1=0$ in (a) and (b).

\z We have the following identity

\begin{eqnarray}
  \label{indentino0}
  \frac{f^{(r)}}{p-\lambda}=\left(r(p-\lambda)^{r-1}\int_k^p
  \frac{f(s)}{(s-\lambda)^{r+1}}\mathrm{d}s+\frac{f(p)}{p-\lambda}\right)^{(r)}
\end{eqnarray}

\z which is proved by straightforward 
differentiation of the r.h.s. or by writing 
$f^{(r)}=(p-\lambda)g^{(r)}=(pg-\lambda g)^{(r)}-
rg^{(r-1)}$ so that $f=(p-\lambda)(\mathcal{P}g)'-r\mathcal{P}g$
and solving  for $\mathcal{P}g$. We take $k\in\NN$ with
$k+1>k_1$, a distribution 
 $f$ with $\mbox{supp}(f)\in (k_0,k+1)$ where
$k_0=\max\{k,k_1\}$ and we let

\begin{eqnarray}
  \label{defc}
  c=\int_{k_0}^{k+1}
  \frac{f(s)}{(s-\lambda)^{r+1}}\mathrm{d}s \ \ \ (r:=mk)
\end{eqnarray}

\z For  $\epsilon$ small  (to be made zero
in the end), we write the decomposition

\begin{multline}
  \label{decomf}
  \frac{f^{(r)}}{p-\lambda}\cr=\left(r(p-\lambda)^{r-1}\int_{k_0}^p
  \frac{f(s)}{(s-\lambda)^{r+1}}\mathrm{d}s+\frac{f(p)}{p-\lambda}
-cr(p-\lambda)^{r-1}\bchi_{[k+1-\epsilon,\infty]}\right)^{(r)}\cr
+\left( cr(p-\lambda)^{r-1}\bchi_{[k+1-\epsilon,\infty]}\right)^{(r)}
=f_1+f_2=f_1+c\,r\, g_2^{(r)}
\end{multline}

\z where by construction
 $f_1\in\mathcal{D}'_{m,\nu}(k,k+1)$
whence, for $f\in \mathcal{D}(k_0,k+1)$ we have

\begin{multline}
  \label{est11}
  \left\|\frac{f^{(r)}}{p-\lambda}\right\|_{m,\nu}\le
\|f_1\|_{m,\nu}+\|f_2\|_{m,\nu}=\|f_1\|_{\mathcal{D}'_{m,\nu}(k,k+1)}+\|f_2\|_{m,\nu}
\cr\le 
\|f_1-f_2\|_{\mathcal{D}'_{m,\nu}(k,k+1)}+2\|f_2\|_{m,\nu}
\end{multline}

\z and then, for some $C_1$

\begin{multline}
  \label{term01}
  \|f_1-f_2\|_{\mathcal{D}'_{m,\nu}(k,k+1)}\cr\le
  \nu^{r}\int_{k_0}^{k+1}\left|r(p-\lambda)^{r-1}\int_{k_0}^p
  \frac{f(s)}{(s-\lambda)^{r+1}}\mathrm{d}s+
  \frac{f(p)}{p-\lambda}\right|\mathrm{e}^{-\nu p} \mathrm{d}p\cr \le
  {\sup}_*|p-\lambda|^{-1}\|f\|_{m,\nu}+{\sup}_*
  \left|\frac{(p-\lambda)^{r-1}}{(s-\lambda)^{r+1}}\right|\|f\|_{m,\nu}\cr
    \le C_1{\sup }_*|p-\lambda|^{-1}\|f\|_{m,\nu}
\end{multline}

\z where the supremum is taken over
$\{k\in\NN,p, s\in[k,k+1]\cap(k_0,\infty) \}$.
For the constant $c$ in  (\ref{decomf})  we have, for some
$C_2$ depending on $a$ and otherwise independent
of $\lambda,k$ the estimate

\begin{eqnarray}
  \label{estimc00}
  |c|\le \frac{C_2 \mathrm{e}^{\nu (k+1)}}{|k_0-\lambda|^{r+1}}\int_{k_0}^{k+1}|f(s)|\mathrm{e}^{-\nu
    s}\mathrm{d}s=\nu^{-r}\frac{C_2 \mathrm{e}^{\nu (k+1)}}{|k_0-\lambda|^{r+1}}\|f\|_{m,\nu}
\end{eqnarray}

\z Let $k'=k+1-\epsilon$. For some 
$C_3=C_3(a)\le \exp\left[(k_0+1)|k_0-\lambda|^{-1}\right]$ we have

\begin{multline}
  \label{normk}
  \|g_2\|_{m,\nu,k}=\nu^r \int_{k'}^{k+1}{\mathrm{e}^{-\nu
      x}}{|x-\lambda|^{r-1}}\mathrm{d}x\le {C_3\nu^{r-1}}{|k_0-\lambda|^{r-1}}\mathrm{e}^{-\nu
    k'}\cr
\Rightarrow c\,r\,\|g_2\|_{m,\nu,k}\le
\nu^{-r}{mk}\frac{\mathrm{e}^{\nu (k+1)}}{|k_0-\lambda|^{r+1}}\|f\|_{m,\nu}
{C_2C_3\nu^{r-1}}{|k_0-\lambda|^{r-1}}\mathrm{e}^{-\nu
    k'}\cr=\frac{C_4 mk
    \mathrm{e}^{\nu\epsilon}}{\nu|k_0-\lambda|^2}\|f\|_{m,\nu}
\end{multline}

\z For $n\ge k+1$ we write (\ref{defDelI})  as

\begin{eqnarray}
  \label{deliprel3}
   \Delta_n(g_2^{(r)})=\bchi_{[n,n+1]}\mathcal{P}^{m}\left(\bchi_{[n,\infty)}\mathcal{P}^{m(n-k-1)}g_2\right)
\end{eqnarray}

\z (cf. (\ref{defDelI})). For $\lambda$ complex we take $K_1(a)=\sup_{s\ge k_0} (s-\lambda)^{-1}
(s+1+|\lambda|)$. We let 
 $\tilde{\lambda}=
\lambda$ if $\lambda$ is real and 
$\tilde{\lambda}=-1-|\lambda|$ otherwise and write $q=m(n-k-1)$. Noting that $K_1^{mk_0+m-1}\le C_5(a,b)=K$  we have

\begin{eqnarray}
  \label{primasconv}
  &&\Gamma(q)\mathcal{P}^{q}|g_2|\le K\int_{k'}^{x}{(x-s)^{q-1}}(s-\tilde{\lambda})^{r-1}\mathrm{d}s\cr&&
\le K\int_0^x{(x-s)^{q-1}}(s-\tilde{\lambda})^{r-1}\mathrm{d}s
=\frac{(x-\tilde{\lambda})^{q+r-1}\Gamma(q)\Gamma(r)}{\Gamma(q+r)}
\end{eqnarray}

\z The estimate above is true for $\tilde{\lambda}\le 1$ but
is ``optimal'' only when
the maximum of the integrand is inside the region
of integration i.e. when
$\tilde{\lambda}>-(2+m^{-1})k(n-k)^{-1}+m^{-1}$. If this is
not the case
 we prefer to 
simply estimate the integral in terms of the maximum of the integrand
over the region of integration. So for, say, $\tilde{\lambda}<-3k$
we use the inequality

\begin{eqnarray}
\label{lambdareal6}
&&\Gamma(q)\mathcal{P}^{q}|g_2|\le K(x-k')^{q}(k'-\tilde{\lambda})^{r-1}
\end{eqnarray}

\z Now, for $\tilde{\lambda}>-3k$, using (\ref{primasconv})
and (\ref{deliprel3})

\begin{eqnarray}
  \label{lambdareal3}
  &&\mathcal{P}^{m}\!\left(\bchi_{[n,\infty)}\mathcal{P}^{m(n-k-1)}|g_2|\right)
  \le\frac{\Gamma(r)}{\Gamma(q+r)\Gamma(m)} \int_n^x
  (x-s)^{m-1}(s-\tilde{\lambda})^{q+r-1}\mathrm{d}s\cr&&\le
\frac{\Gamma(r)}{\Gamma(q+r)\Gamma(m)}(x-n)^m(x-\tilde{\lambda})^{q+r-1}
\cr&&
\end{eqnarray}

\z (as $m$ is fixed we do not
lose too much by this evaluation which has the advantage
of preserving the behavior near $x=n$). Further, we have

\begin{multline}
  \label{evalsum0}
  \int_{n}^{n+1}\mathrm{e}^{-\nu x}(x-n)^m(x-\tilde{\lambda})^{q+r-1}\mathrm{d}x
\le (n+1-\tilde{\lambda})^{q+r-1}\int_{n}^{n+1}
\mathrm{e}^{-\nu x}(x-n)^m \mathrm{d}x\cr
\le\frac{\Gamma(m)}{\nu^m} (n+1-\tilde{\lambda})^{q+r-1}
\end{multline}

\z and

\begin{eqnarray}
  \label{evalsum1}
 &&\sum_{n=k+1}^{\infty}\nu^{mn} \|\Delta\|_{m,\nu}\le
\Gamma(mk)K\sum_{n=k+1}^{\infty}\frac{\nu^{m(n-1)} \mathrm{e}^{-\nu
    n}(n+1-\tilde{\lambda})^{mn-1-m}}{\Gamma(m(n-1))}\cr&&
\end{eqnarray}

\z The ratio of two successive terms $s_{n+1}/s_n$ of
the infinite series above is estimated by:

\begin{eqnarray}
  \label{evalsum2}
  \nu^m
  \mathrm{e}^{-\nu}\mathrm{e}^{\frac{mn-1-m}{n+1-\tilde{\lambda}}}\left(\frac{n+2-\tilde{\lambda}}{mn-m}\right)^m
\le \frac{1}{2}
\end{eqnarray}

\z when $\nu>C_1$ for some $C_1$ independent of $k,n,\tilde{\lambda}$ in the
region $k>1,n>k,\tilde{\lambda}\in(-3k,1)$.
This means that

\begin{eqnarray}
  \label{evalsum3}
  &&\sum_{n=k+1}^{\infty}\nu^{mn} \|\Delta\|_{m,\nu}\le
2K{\nu^{mk} \mathrm{e}^{-\nu
    (k+1)}(k+2-\tilde{\lambda})^{mk-1}}\cr&&
\end{eqnarray}

\z and combining with (\ref{estimc00}) and (\ref{est11})
we have

\begin{multline}
  \label{fiest1}
  2\|f_2\|\le 
4K\nu^{-mk}\frac{\mathrm{e}^{\nu (k+1)}}{(k-\tilde{\lambda})^{mk+1}}\|f\|_{m,\nu}
{\nu^{mk} \mathrm{e}^{-\nu
    (k+1)}(k+2-\tilde{\lambda})^{mk-1}}\cr\le
4K\mathrm{e}^{2\frac{(mk+1)}{k-\tilde{\lambda}}}\frac{\|f\|_{m,\nu}}{(k-\tilde{\lambda})^2}\le
4K\mathrm{e}^{\frac{4m}{k_0\tilde{\lambda}}}\frac{\|f\|_{m,\nu}}{(k_0-\tilde{\lambda})^2}\le\frac{C_6(a,b)\|f\|_{m,\nu}}{|k_0-\lambda|^2}
\end{multline}

\z If we started with (\ref{lambdareal6}) we would have obtained in the same way,
for $\tilde{\lambda}<-3k$,
instead of (\ref{lambdareal3}), 

\begin{eqnarray}
  \label{lambdareal10}
  \frac{1}{\Gamma(q)\Gamma(m)}(x-n)^m(x-k')^{q}(k'-\tilde{\lambda})^{r-1}
\end{eqnarray}

\z and the calculations
are similar from this point on. Condition (\ref{evalsum2}) is of the same type, with $k'$ replacing
$-\tilde{\lambda}$ and final estimate is 

\begin{eqnarray}
  \label{fiesta10}
  C_7 {\|f\|_{m,\nu}}
  \frac{(k'-\tilde{\lambda})^{r-1}}{(k-\tilde{\lambda})^{r+1}}
\le C_8 \frac{\|f\|_{m,\nu}}{(k_0-\tilde{\lambda})^2}
\end{eqnarray}

\z Finally we take the limit $\epsilon\rightarrow 0$ 
and noting that $K\rightarrow 0$ as $\lambda\rightarrow\infty$,
(i) is proved.

For (ii), merely notice that $U(\lambda_2)-U(\lambda_1)
=(\lambda_2-\lambda_1)U(\lambda_1)U(\lambda_2)$.

\Box

\begin{Remark}\label{density}
Let $\psi\in\lone[0,1]$ with the property
$\int_0^1\psi(t)\phi^{(m)}(t)\mathrm{d}t=0$ for all $\phi\in\mathcal{D}(0,1)$.
Then $\psi$ is a polynomial of degree at most $m-1$.

\end{Remark}

\z This is a well-known property. We sketch an elementary proof for
$m=1$ (for general $m$ the proof is similar).  Let $x\in(0,1)$, and
consider a sequence $\chi_n$ in $ \mathcal{D}(0,x)$
$\lone$--convergent to $\bchi_{[0,x]}$. Then
$\phi_n(t):=\chi_n(t)-\kappa^{-1}\chi_{n}(\kappa(1-t))$ with
$\kappa=x(1-x)^{-1}$ converges to $\bchi_{[0,x]}-\kappa\bchi_{[x,1]}$.
Furthermore, since $\int_0^1\phi_n(t)\mathrm{d}t=0$ we have
$\Phi_n(t):=\int_0^t\phi_n(s)\mathrm{d}s\in\mathcal{D}(0,1)$. Since
$\Phi_n'+\kappa^{-1}\rightarrow (1+\kappa^{-1})\bchi_{[0,x]}$ it
follows that $\int_0^x(\psi-C)=0$, where $C=\int_0^1\psi(t)\mathrm{d}t$. Thus
$\psi=C\ a.e.$

\subsection{Derivation of the equations for the transseries.}
\label{sec:For}

Consider first the scalar equation

\begin{eqnarray}
  \label{eqscal0}
  y'=f_0(x)-\lambda y-{x}^{-1}By+g(x,y)=
-y+x^{-1}By+\sum_{k=1}^{\infty}g_k(x)y^k
\end{eqnarray} 

\z For $x\rightarrow+\infty$ we take 
\begin{eqnarray}\label{tr1}
{y}=\sum_{k=0}^{\infty}{y}_k \mathrm{e}^{-kx}
\end{eqnarray}

\z where ${y}_k$ will be either formal series
$x^{-s_k}\sum_{n=0}^{\infty}a_{kn}x^{-n}$, with $a_{k,0}\ne 0$  or
actual functions with the condition that (\ref{tr1}) converges uniformly. As a
transseries, (\ref{tr1}) can be also understood as a well ordered
double sequence $t_{kn}=x^{p_{kn}}\mathrm{e}^{-kx}$, with $p_{k\,n+1}<p_{kn}$.
(The order relation  is $x^{p}\mathrm{e}^{-kx}\gg x^{p'}\mathrm{e}^{-k'x}
$ as $x\rightarrow+\infty$ iff $k<k'$ or
$k=k'$ and $p>p'$; thus a strictly {\em increasing} sequence of
terms of a transseries necessarily terminates.) Power series are a
special case of transseries, with $y_1=y_2=\ldots=0$.  Two transseries
$\sum_{k=0}^{\infty}y_k\mathrm{e}^{-kx}$ coincide iff all corresponding
component power series $y_k$ coincide.  Transseries of this type are
closed under addition, multiplication and infinite sums of the form
involved in (\ref{eqscal0}) (this last aspect will become clear in the
calculation leading to (\ref{eqcompl}) below). Note that well-ordering
plays an important part in defining multiplication of transseries; in
contrast, for the unrestricted formal expansion
$S=\sum_{k=-\infty}^{\infty}x^k$, no immediate meaning can be given to
$S^2$.  Let $y_0$ be the first term in (\ref{tr1}) and $\delta=y-y_0$.
We have

\begin{multline}
  \label{part00}
  y^k-y_0^k-ky_0^{k-1}\delta=\sum_{j=2}^k\binom{k}{j}y_0^{k-j}\delta^j=
\sum_{j=2}^k\binom{k}{j}y_0^{k-j}\sum_{i_1,\ldots,i_j=1}^{\infty}
\prod_{s=1}^j\left(y_{i_s}\mathrm{e}^{-i_s x}\right)\cr
=\sum_{m=1}^{\infty}\mathrm{e}^{-mx}\sum_{j=2}^k\binom{k}{j}y_0^{k-j}\sum_{(i_s)}^{(m;j)}
\prod_{s=1}^j y_{i_s}
\end{multline}

\z where $\sum_{(i_s)}^{(m;j)}$ means the sum over
all positive integers $i_1,i_2,\ldots,i_j$ with the restriction
$i_1+i_2+\cdots+i_j=m$. Let $d_1=\sum_{k\ge 1}k g_k y_0^{k-1}$.
Introducing $y=y_0+\delta$ in (\ref{eqscal0})
and equating the coefficients of $\mathrm{e}^{-lx}$ we get, by separating
the terms containing $y_l$ for $l\ge 1$ and 
interchanging the $j,k$ orders of summation,

\begin{multline}
  \label{eqcompl}
  y_l'+(\lambda(1-l)+x^{-1}B-d_1(x))y_l=\sum_{j=2}^{\infty}\sum_{(i_s)}^{(l;j)}
\prod_{s=1}^j y_{i_s}\sum_{k\ge\{2,j\}}\binom{k}{j} g_k
y_0^{k-j}
\cr=\sum_{j=2}^{l}\sum_{(i_s)}^{(l;j)}
\prod_{s=1}^j y_{i_s}\sum_{k\ge\{2,j\}}\binom{k}{ j} g_k
y_0^{k-j}=:\sum_{j=2}^{l}d_j(x)\sum_{(i_s)}^{(l;j)}
\prod_{s=1}^j y_{i_s} 
\end{multline}

\z where for the middle equality we note that the infinite sum
terminates because $i_s\ge 1$ and $\sum_{s=1}^j i_s=l$. The fact
mentioned before that $\sum_{k=1}^{\infty}g_k(x)y_k$ is well defined
when $y_k$ are formal series is now visible: collecting the
coefficient of $x^{p}\mathrm{e}^{-kx}$, only \emph{finite} sums of coefficients
appear.

For a vectorial equation like (\ref{eqor}) we first write

\begin{eqnarray}
  \label{eqsvec0}
  \bfy'=\bff_0(x)
-{\hat{\Lambda}}\bfy+{x}^{-1}{{\hat{B}}}\bfy+\sum_{\bfk\succ 0}{\bf g}_\bfk(x)\bfy^\bfk
\end{eqnarray} 

\z with $\bfy^{\bfk}:=\prod_{i=1}^{n_1} (\bfy)_i^{k_i}$.
 The formal operations
and ordering extend naturally to the vectorial
general transseries (\ref{eqformgen,n}), under
the restriction $\Re(\bfk\cdot\bflam x)>0$
As with (\ref{eqcompl}), we introduce  the transseries
(\ref{eqformgen,n}) in (\ref{eqsvec0}) and equate the
coefficients of $\exp(-\bfk\cdot\bflam x)$. Let
 $\mathbf{v}_\bfk=x^{-\bfk\cdot\bfm}\mathbf{y}_\bfk$ and

\begin{eqnarray}
  \label{defD}
  \bfdd_\bfj(x)=\sum_{\bfl\ge \bfj}\binom{\bfl}{
    \bfj}\bfgg_\bfl(x)\bfv_0^{\bfl-\bfj}
\end{eqnarray}

\z Noting that,
by assumption, $\bfk\cdot\bflam=\bfk'\cdot\bflam\Leftrightarrow
\bfk=\bfk'$ we obtain, for $\bfk\in\NN^{n_1}$, $\bfk\succ 0$

\begin{eqnarray}
  \label{eqmygen}
  &&\bfv_\bfk '+\left({\hat{\Lambda}}-\bfk\cdot\bflam \hat{I}
  +x^{-1}{\hat{B}}\right)\bfv_\bfk+\sum_{|\bfj|=1}\bfdd_{\bfj}(x)(\bfv_\bfk)^{\bfj}\cr&& = \sum_{\stackrel{\scriptstyle \bf \phantom{0}j\le
      k}{|\bfj|\ge 2}}\bfdd_\bfj (x)\sum_{(\bfii_{mp}:\bfk)}
  \prod_{m=1}^n\prod_{p=1}^{j_m}
 \left( \bfv_{\bfii_{mp}}\right)_m=\mathbf{t}_\bfk(\bfv)
\end{eqnarray}
 
\z where $\binom{\bf l}{j}= \prod_{j=1}^n\binom{l_i}{j_i}$, $(\bfv)_m$
means the component $m$ of $\bfv$, and $\sum_{(\bfii_{mp}:\bfk)}$
stands for the sum over all vectors $\bfii_{mp}\in\NN^n$, with $p\le
j_m,m\le n$, such that $\bfii_{mp}\succ 0$ and
$\sum_{m=1}^n\sum_{p=1}^{j_m}\bfii_{mp}=\bfk$.  We use the convention
$\prod_{\emptyset} =1, \sum_{\emptyset}=0$.  With
$m_i=1-\lfloor\beta_i\rfloor$ we obtain for $\bfy_\bfk$

\begin{gather}
  \label{homogeqv}
  \bfy_\bfk '+\left({\hat{\Lambda}}-\bfk\cdot\bflam \hat{I}
  +x^{-1}({\hat{B}}+\bfk\cdot\bfm)\right)\bfy_\bfk+\sum_{|\bfj|=1}\bfdd_{\bfj}(x)(\bfy_\bfk)^{\bfj}
  = \bft_\bfk(\bfy)
\end{gather}

There are clearly finitely many 
terms in $\mathbf{t}_\bfk(\bfy)$. To find a (not too unrealistic) upper
bound for this number of terms,  we  compare with
 $\sum_{(\bfii_{mp})'}$ which stands for the same as
$\sum_{(\bfii_{mp})}$ except with $\bfii\ge 0$ instead of
$\bfii\succ 0$. Noting that $\binom{k+s-1}{ s-1}=\sum_{a_1+\ldots+a_s=k} 1$
is the number of ways $k$ can be written as a sum of $s$ integers,
we have

\begin{gather}
  \label{combineq}
  \sum_{(\bfii_{mp})}1\le \sum_{(\bfii_{mp})'} 1
=\prod_{l=1}^{n_1}\sum_{(\bfii_{mp})_l}1=
\prod_{l=1}^{n_1}\binom{k_l+|\bfj|-1}{|\bfj|-1}\le \binom{|\bfk|+|\bfj|-1}{
|\bfj|-1}^{n_1}
\end{gather}

\begin{Remark}\label{homogstruct}
Equation  (\ref{eqmygen})  can be written in the form (\ref{homogeq})

\end{Remark}

{\em Proof.}  The fact that only predecessors of $\bfk$ are involved
in $\bft(\bfy_0,\cdot)$ and the homogeneity property of
$\bft(\bfy_0,\cdot)$ follow immediately by combining the conditions
$\sum {\bfii_{mp}}=\bfk$ and $\bfii_{mp}\succ 0$. \Box

\z The formal inverse Laplace transform of (\ref{homogeqv}) is 
then

\begin{eqnarray}
  \label{invlapvk11}
  &&\left(-p+\hat{\Lambda}-\bfk\cdot\bflam\right)\bfY_\bfk
+\left(\hat{B}+\bfk\cdot\bfm\right)
\mathcal{P}\bfY_\bfk+\sum_{|\bfj|=1}\bfd_\bfj*\left(\bfY_\bfk\right)^\bfj
=\bfT_\bfk(\bfY)\cr&&
\end{eqnarray}

\z with 

\begin{equation}
  \label{defT}
  \bfT_\bfk(\bfY)=\bfT\left(\bfY_0,\{\bfY_{\bfk'}\}_{0\prec\bfk'\prec\bfk}\right)=
  \sum_{\bfj\le k;\ |\bfj|>1}\mathbf{D}_\bfj(p)*\sum_{(\bfii_{mp};\bfk)}
  \sideset{^*}{}\prod_{m=1}^{n_1}\sideset{^*}{}\prod_{p=1}^{j_m}\left(\bfY_{\bfii_{mp}}\right)_m
\end{equation}

\z and

\begin{gather}
  \label{definitionDj}
  \mathbf{D}_\mathbf{j}=\sum_{\bfl\ge\bfm}\binom{\bfl}{\bfm}\mathbf{G}_\bfl*
\mathbf{Y}_0^{*(\bfl-\bfm)}+\sum_{\bfl\ge\bfm;|\bfl|\ge
  2}
\binom{\bfl}{\bfm}\mathbf{g}_{0,\bfl}\mathbf{Y}_0^{*(\bfl-\bfm)}
\end{gather}

\begin{subsection}{Useful formulas}
\label{usefulfor}

   A straightforward computation shows that

\begin{eqnarray}\label{u1}{\cal B}(\frac{1}{x^n})=\frac{p^{n-1}}{\Gamma(n)}
\ \mbox{or}\ {\cal L}(p^{n})=\frac{\Gamma(n+1)}{x^{n+1}}\end{eqnarray}

\begin{eqnarray}\label{u3}p^q*p^r=\frac{\Gamma(q+1)\Gamma(r+1)}{\Gamma(q+r+2)}
p^{q+r+1}\end{eqnarray}
Also,  with $f_{1,2}(p):=p\mapsto\heav(p-k_{1,2})g_{1,2}(p-k_{1,2})$
we have
\begin{eqnarray}\label{u4}\Big(f_1*f_2\Big)(p)=\heav(p-k_1-k_2)\Big(g_1*g_2\Big)(p-k_1-k_2)\end{eqnarray}
\end{subsection}

\begin{section}{Acknowledgments}
  
  Special thanks are due to Professors Martin Kruskal and Jean
  \'Ecalle for their many in-depth comments, and to Professor B. L. J.
  Braaksma for carefully reading the manuscript and pointing out an
  error in a formula. The author is grateful to Professors Michael
  Berry and Percy Deift for very interesting discussions.

\end{section}


\begin{thebibliography}{99}
\bibitem{Ecalle-book} J. \'Ecalle {\em Fonctions 
Resurgentes, Publications Mathematiques D'Orsay, 1981}
\bibitem{Ecalle}
 J. \'Ecalle {\em in Bifurcations and periodic orbits of 
vector fields NATO ASI Series, Vol. 408,  1993}
\bibitem{Ecalle2}  J. \'Ecalle {\em Finitude des cycles limites.., Preprint
90-36 of Universite de Paris-Sud, 1990} 
\bibitem{Ecalle3} J. \'Ecalle, F. Menous {Well behaved averages
and the non-accumulation theorem..} Preprint
\bibitem{BRBS} W. Balser, B.L.J. Braaksma, J-P. Ramis, Y. Sibuya
{\em Asymptotic Anal. {\bf 5}, no. 1 (1991), 27-45}
\bibitem{Braaksma} B. L. J. Braaksma {\em Ann. Inst. Fourier,
Grenoble,{\bf 42}, 3 (1992), 517-540}
\bibitem{Balser} Balser, W. {\em From divergent power series to
    analytic functions, Springer-Verlag, (1994).}
\bibitem{Borel} Borel, E. {\em Lecons sur les series divergentes,
Gauthier-Villars, 1901}
\bibitem{Hardy} Hardy, C. G. {\em Divergent series}

\bibitem{Stokes} G. G. Stokes {\em Trans. Camb. Phil. Soc {\bf 10}
106-128}. Reprinted in {\em Mathematical and Physical papers by late
sir George Gabriel Stokes. Cambridge University Press 1904, vol. IV,
77-109}
\bibitem{Wasow} W. Wasow {\em  Asymptotic expansions
for ordinary differential equations, Interscience Publishers 1968 }
\bibitem{Sibuya} Y. Sibuya {\em Global theory of a second order linear
ordinary differential equation with a polynomial coefficient , North-Holland 1975}
\bibitem{Costin} O. Costin {\em IMRN 8, 377-417 (1995)}
\bibitem{CK1} O. Costin, M.D. Kruskal {\em Proc. R. Soc. Lond. A 452,
    1057-1085 (1996)}
\bibitem{CK2} O. Costin, M.D. Kruskal {\em in preparation}
\bibitem{Costin3} O. Costin {\em in preparation}
\bibitem{Cope} F. T. Cope {\em Amer. J. Math. vol. 56 pp 411-437 (1934)}
\bibitem{Ritt} J.F. Ritt   {\em Differential algebra,
 American Mathematical Society, New York 1950}
\bibitem{To1} A. Tovbis{ Linear
   Algebra and Applications, { 162-164}, 389-407 (1992)}.
\bibitem{Fabry} C. E. Fabry {\em Th\`ese (Facult\'e des Sciences), Paris, 1885}
\bibitem{Iwano}M. Iwano {\em Ann. Mat. Pura Appl. (4) {\bf 44} 1957, 261-292}
\bibitem{Berry} M.V. Berry {\em  Proc. R. Soc. Lond. A 422, 7-21,
1989}
\bibitem{Berry-hyp} M.V. Berry {\em  Proc. R. Soc. Lond. A 430, 653-668,
1990}

\bibitem{Berry-Howls} M.V. Berry, C.J. Howls {\em
    Proc. Roy. Soc. London Ser. A 443 no. 1917, 107--126 (1993)}
\bibitem{Berry-gamma} M.V. Berry {\em  Proc. Roy. Soc.
London Ser. A 434  no. 1891, 465--472. (1991)}
\bibitem{Confe}H. Segur, S.
Tanveer and H. Levine, ed. {\em Asymptotics Beyond all
Orders,  Plenum Press 1991}
\bibitem{Kruskal} M.D. Kruskal, H. Segur {\em Studies
in Applied Mathematics 85:129-181, 1991}
\bibitem{Holland} A.S.B. Holland,{\em Introduction
to the theory of entire functions, Academic Press, 1973}
\end{thebibliography}
\end{document}